\documentclass{amsart}
\usepackage{color}

\usepackage{amsmath} 
\usepackage{amssymb}
\usepackage{mathrsfs}

\newtheorem{theorem}{Theorem}[section] 
\newtheorem{claim}[theorem]{Claim}

\newtheorem{observation}[theorem]{Observation} 
 
\newtheorem{conclusion}[theorem]{Conclusion}

\theoremstyle{definition}
\newtheorem{definition}[theorem]{Definition}
\newtheorem{dc}[theorem]{Definition/Claim}

\newtheorem{discussion}[theorem]{Discussion}
\newtheorem{convention}[theorem]{Convention}

\newtheorem{hypothesis}[theorem]{Hypothesis}

\theoremstyle{remark}
\newtheorem{remark}[theorem]{Remark}
\newtheorem{question}[theorem]{Question}
\newtheorem{notation}[theorem]{Notation}
\newtheorem{context}[theorem]{Context}
\newtheorem{thesis}[theorem]{Thesis}

\newcommand{\rest}{{\restriction}}
\newcommand{\dom}{{\rm dom}} 
\newcommand{\hrtg}{{\rm hrtg}} 
\newcommand{\Ax}{{\rm Ax}} 
\newcommand{\DC}{{\rm DC}} 
 
\newcommand{\PC}{{\rm PC}}

\newcommand{\Card}{{\rm Card}} 
\newcommand{\ZF}{{\rm ZF}} 
\newcommand{\Ord}{{\rm Ord}} 
\newcommand{\Hom}{{\rm Hom}} 
\newcommand{\Min}{{\rm Min}} 
 
\newcommand{\Fil}{{\rm Fil}}

\newcommand{\Car}{{\rm Car}} 
\newcommand{\otp}{{\rm otp}}

\newcommand{\tcf}{{\rm tcf}}

\newcommand{\fil}{{\rm fil}} 
\newcommand{\dual}{{\rm dual}}
\newcommand{\Depth}{{\rm Depth}}

\newcommand{\purely}{{\rm purely}} 
 
\newcommand{\Range}{{\rm Range}} 
\newcommand{\pcf}{{\rm pcf}} 
\newcommand{\cov}{{\rm cov}} 
\newcommand{\pp}{{\rm pp}}

\newcommand{\eq}{{\rm eq}} 
 
\newcommand{\rk}{{\rm rk}} 
\newcommand{\qu}{{\rm qu}} 
\newcommand{\AC}{{\rm AC}} 
\newcommand{\ZFC}{{\rm ZFC}} 
\newcommand{\FIL}{{\rm FIL}} 
\newcommand{\Rang}{{\rm Rang}} 
\newcommand{\eub}{{\rm eub}} 
\newcommand{\Dom}{{\rm Dom}} 
\newcommand{\wlor}{{\rm wlor}} 
 
\newcommand{\supp}{{\rm supp}} 

\newcommand{\set}{{\rm set}}

\newcommand{\card}{{\rm card}}
\newcommand{\wilog}{{\rm without loss of generality}}
\newcommand{\Wilog}{{\rm Without loss of generality}}

\newcommand{\then}{{\underline{then}}}
\newcommand{\when}{{\underline{when}}}
\newcommand{\Then}{{\underline{Then}}}

\newcommand{\Iff}{{\underline{iff}}}
\newcommand{\mn}{{\medskip\noindent}}
\newcommand{\sn}{{\smallskip\noindent}}

\newcommand{\cA}{{\mathscr A}}

\newcommand{\cE}{{\mathscr E}}

\newcommand{\gb}{{\mathfrak b}}
\newcommand{\ga}{{\mathfrak a}}

\newcommand{\varp}{{\varepsilon}}

\newcommand{\cH}{{\mathcal H}}

\newcommand{\cF}{{\mathscr F}}
\newcommand{\cG}{{\mathscr G}}

\newcommand{\bbL}{{\mathbb L}}

\newcommand{\bbN}{{\mathbb N}}
\newcommand{\cP}{{\mathscr P}}

\newcommand{\bbQ}{{\mathbb Q}}
\newcommand{\bbR}{{\mathbb R}}
\newcommand{\bbZ}{{\mathbb Z}}

\newcommand{\cS}{{\mathscr S}}

\newcommand{\gy}{{\mathfrak y}} 
\newcommand{\gz}{{\mathfrak z}} 
\newcommand{\gY}{{\mathfrak Y}} 
\newcommand{\cU}{{\mathscr U}}

\newcommand{\cW}{{\mathscr W}}

\newcommand{\cZ}{{\mathscr Z}}

\newcommand{\cf}{{\rm cf}}
\newcommand{\id}{{\rm id}}
\newcommand{\pr}{{\rm pr}}

\newcount\skewfactor
\def\mathunderaccent#1#2 {\let\theaccent#1\skewfactor#2
\mathpalette\putaccentunder}
\def\putaccentunder#1#2{\oalign{$#1#2$\crcr\hidewidth
\vbox to.2ex{\hbox{$#1\skew\skewfactor\theaccent{}$}\vss}\hidewidth}}

\newenvironment{PROOF}[2][\proofname.]
   {\begin{proof}[#1]\renewcommand{\qedsymbol}{\textsquare$_{\rm #2}$}}
   {\end{proof}}

\usepackage[hidelinks]{hyperref}

\begin{document}
\makeatletter\def\shfiuwefootnote{\gdef\@thefnmark{}\@footnotetext}\makeatother\shfiuwefootnote{Version 2021-09-15. See \url{https://shelah.logic.at/papers/1005/} for possible updates.}

\title {ZF+DC+AX$_4$ \\
Sh1005}
\author {Saharon Shelah}
\address{Einstein Institute of Mathematics\\
Edmond J. Safra Campus, Givat Ram\\
The Hebrew University of Jerusalem\\
Jerusalem, 9190401, Israel\\
 and \\
 Department of Mathematics\\
Hill Center - Busch Campus \\ 
 Rutgers, The State University of New Jersey \\
 110 Frelinghuysen Road \\
 Piscataway, NJ 08854-8019 USA}
\email{shelah@math.huji.ac.il}
\urladdr{http://shelah.logic.at}
\thanks{The author thanks Alice Leonhardt for the beautiful typing.
  The author would like to thank the ISF for partial support of this
  research (Grant No. 1053/11).
Publication 1005.  First typed August 12, 2010.}

\subjclass[2010]{Primary 03E04, 03E25; Secondary: }

\keywords {set theory, weak axiom of choice, pcf, abelian groups,
  group. References to outside papers like \cite[2.13=Ls.2]{Sh:835}
  means to 2.13 where s.2 is the label used there, so intended only to
  help the author if more is added to \cite{Sh:835}}



\date{February 3, 2016}

\dedicatory {Dedicated to the memory of Richard Laver}

\begin{abstract}
We consider mainly the following version of set theory:
``ZF + DC and for every $\lambda,\lambda^{\aleph_0}$ 
is well ordered", our thesis is that this is a reasonable set theory,
e.g. on the one hand it is much weaker than full choice, and on the
other hand much can be said or at least this is what the present work
tries to indicate.  In particular, we prove that
for a sequence $\bar\delta = \langle \delta_s:s \in
Y\rangle,\cf(\delta_s)$ large enough compared to $Y$, we can prove the
pcf theorem with minor changes (in particular, using true cofinalities not the
pseudo ones).  We then deduce the existence of covering numbers and
define and prove existence of a class of true successor cardinals.  
Using this we give some
diagonalization arguments (more specifically some black boxes and
consequences) on Abelian groups, chosen as a characteristic case.  
We end by showing that some such
consequences hold even in ZF above.
\end{abstract}

\maketitle
\numberwithin{equation}{section}
\setcounter{section}{-1}

\newpage

\centerline {Anotated Content}
\bigskip

\noindent
\S0 \quad Introduction, (labels z -), pg.\pageref{0} 

\S(0A) \quad Background and results, pg.\pageref{0A}
\mn
\begin{enumerate}
\item[${{}}$]  [We investigate $\ZF + \DC + \Ax_4$ asserting it is
  quite strong, not like the chaos usually related to universes
  without choices.  We consider using weaker versions and relatives
of $\Ax_4$ but not in the Anotated Content.]
\end{enumerate}

\S(0B) \quad Preliminaries, pg.\pageref{0B}
\mn
\begin{enumerate}
\item[${{}}$]  [We define Ax$_4$, Ax$_{4,\delta}$, prove that a suitable
 closure operation  $c \ell$ exists, and define ``$\partial$-uniformly
 definable".  We also define ``$\gy$-eub", tcf and $A \le_{\qu} B$.]
\end{enumerate}
\bigskip

\noindent
\S1 \quad The pcf theorem again, (labels c-), pg.\pageref{1}
\mn
\begin{enumerate}
\item[${{}}$]  [The version of the pcf theorem proved here is quite strong.
  Assume $\bar\delta = \langle \delta_s:s \in Y \rangle,\cf(\delta_s)$
  large enough compared to $Y$; we do not demand 
``$\delta_s$ regular cardinals". 
We prove first the existence of scales for $\aleph_1$-complete filters;
  note that we have said ``for any $Y$" in spite of our having 
$\Ax_4$ (or less) only.  Then we prove that we have $\langle
  (D_\varepsilon,A_\varepsilon/D_\varepsilon,\bar
  f_\varepsilon):\varepsilon \le \varepsilon_*\rangle$ as usual (so
$D_\varepsilon$ not necessarily $\aleph_1$-complete) but
\sn
\begin{enumerate}
\item[$(a)$]   $\ell g(\bar f_\varepsilon)$ is not necessarily 
a regular cardinal,
\sn
\item[$(b)$]  the cofinality of $\ell g(\bar f_\varepsilon)$ is
not necessarily increasing
\sn
\item[$(c)$]   as generators, for the time being we have only
  $A_\varepsilon/D_\varepsilon$ not $A_\varepsilon$.
\end{enumerate}
\sn
However, here there is a gain compared to the ZFC version 
because of a new phenomena: the 
results apply also when many (even all) $\delta_s$ have 
small cofinality but $\bar\delta$ does not;
expressed by $\cf-\id_{< \theta}(\bar\delta)$.  Of course, an
additional gain is that the objects above are definable (from a well
ordering of some $[\lambda]^{\aleph_0}$).]
\end{enumerate}
\bigskip

\noindent
\S2 \quad More on the pcf theorem, pg.\pageref{2}

\S(2A) \quad When Cofinalities are smaller, pg.\pageref{2A}
\mn
\begin{enumerate}
\item[${{}}$]  [A drawback of \S1 is that we need 
$\cf-\id_{< \theta}(\bar\delta)$ where $\theta >
  \hrtg(\cP(\cP(Y)))$.  We weaken the assumption to
$\Ax_{4,\infty,\partial,\kappa}$ with possibly $\kappa > \aleph_1$.
If $Y$ is countable we can weaken the large cofinality 
demand to $> \aleph_1$.  Moreover, there is a pcf analysis of $(\Pi
\bar\delta,<_{\cf-\id_{< \theta}(\bar\delta)})$ iff there is a well
 orderable $\cF \subseteq \Pi \bar\delta$ which is 
$<_{\id-\cf_{<\theta}(\Pi \bar \delta)}$-cofinal, and we can choose
 generators $A_\varepsilon$ under reasonable conditions.]
\end{enumerate}

\S(2B) \quad Elaborations, pg.\pageref{2B}
\mn
\begin{enumerate}
\item[${{}}$]  [We revisit some points.  We give a sharper version of
  the result of \cite{Sh:835} that ${}^\kappa \lambda$ can be divided
  to few (really $X_\kappa = {}^\omega(\Fil^4_{\aleph_1}(Y))$ well
  ordered subsets (in Theorem \ref{d29}).  We also reconsider the eub-existence
  (in \ref{d19}), existence of $\langle e_\alpha:\alpha\rangle$ and
  existence of $u$ with a minimal $c \ell(u)$ such that $u$ includes a
  club of $\delta_s$ for $s \in Y$ (in \ref{d17}).  We finish getting
  essential equality in $\hrtg({}^\kappa \mu)$, \ wilog \, $({}^\kappa
  \mu)$ and so called o-$\Depth^0_\kappa({}^\kappa \mu)$ in \ref{d31}.
See \ref{d19}, \ref{d29}, \ref{d34}.]
\end{enumerate}

\S(2C) \quad True successor cardinal, pg.\pageref{2C}
\mn
\begin{enumerate}
\item[${{}}$]  [See \ref{p16}.  We say that $\lambda$ is true
  successor cardinal when $\lambda = \mu^+$ and there is $\bar f =
  \langle f_\alpha:\alpha \in [\mu,\lambda\rangle,f_\alpha$ a
  one-to-one function from $\mu$ onto $\alpha$.  We investigate this
  notion in particular proving many successor cardinals are true
  successor cardinals.]
\end{enumerate}

\S(2D) \quad Covering numbers, pg.\pageref{2D}
\mn
\begin{enumerate}
\item[${{}}$]  [We prove that covering number exists.  Note that we can
  present the results: if $\mathbf L[X]$ contain
  $[\lambda_*]^{\aleph_0}$ then ``enough below $\lambda_*$", $\mathbf
  L[X]$ is closed enough to $\mathbf V$ by covering lemmas, singulars
  being true successor, etc.]
\end{enumerate}
\bigskip

\noindent
\S3 \quad Black boxes, pg. \pageref{3}
\mn
\begin{enumerate}
\item[${{}}$]  [Normally theorems using diagonalization used choice
  quite heavily.  We show that at least for one way (one kind of 
black boxes), $\Ax_4$
  suffice.]
\end{enumerate}

\S(3A) \quad Existence proof, pg.\pageref{3A}
\mn
\begin{enumerate}
\item[${{}}$]  [We show that using $\ZF + \Ax_4$, we can prove a Black
 Box which has been used not a few times, e.g. in the book of
 Eklof-Mekler \cite{EM02} and in the book of G\"obel-Trlifaj
 \cite{GbTl12}.  We then as an example prove one such theorem: the existence
  of an $\aleph_1$-free Abelian group with trivial dual.]
\end{enumerate}

\S(3B) \quad Black boxes with no choice, pg.\pageref{3B}
\mn
\begin{enumerate}
\item[${{}}$]  [Here we go in another direction: we try to build
 examples on sets which are are not well ordered, working in $\ZF$
 only.]
\end{enumerate}
\newpage

\section {Introduction} \label{0}

\subsection {Background and Results} \label{0A}\

Everyone knows that the issue of weakening AC, the axiom of 
choice issue, is dead,
settled, as naturally the axiom of choice is true, and its weakenings
lead to bizarre universes on which there is not much to be proved, 
\underline{or} assuming $\AC$ is irrelevant (as in inner models).

The works on determinacy are not a real exception: it e.g. replace
Borel sets and projective sets by sets in $\bbL[\bbR]$, so have much
to say on this inner model, for which the only choice missing is a well
ordering of $\cP(\bbN)$.  In \cite{Sh:835} we suggest to consider
several related axioms, the strongest of them being $\Ax_4$, 
assuming ZF + DC of course.  It is 
in a sense an anti-thesis to considering $\bbL[\bbR]$: 
it says we can well order (not all the subsets just) the
countable subsets of any ordinal. This was continued in \cite{Sh:937},
\cite{Sh:955} and in Larson-Shelah \cite{Sh:925}. We may wonder how
to get natural models of ZF + DC + Ax$_4$.  Such a natural model 
is gotten starting with $\mathbf V \models$ G.C.H. and
forcing by the choiceless version of Easton forcing except for
$\aleph_0$.

While \cite{Sh:497} claims to prove that ``the theory of
$\pcf$ with weak choice is
non-empty", \cite{Sh:835} seems to us the true beginning of such set
theory, proving (in $\ZFC + \DC + \Ax_4$ or so): there is a class of 
successor regular cardinals, and for any set
$Y,{}^Y \lambda$ can, in a suitable sense, 
be decomposed to ``few" well order sets (see \cite[0.3]{Sh:835} and
more here in \ref{d29}).

Much attention there was given to trying to get the results from
weaker relatives of Ax$_4$.  A major aim of this work is to try to justify:
\begin{thesis}
\label{y4}
ZF + DC + Ax$_4$ is a reasonable set theory, for which much of
combinatorial set theory can be generalized, but many times in a
challenging way and even discover new phenomena.
\end{thesis}

\noindent
In particular we consider diagonalization arguments, including in $\ZF$ alone.
Returning to the original issue, i.e. the position that ``set theory
with weak choice is dead", which we had wholeheartedly 
supported, the paper's position here is that:
\mn
\begin{enumerate}
\item[$(a)$]  AC is obviously true
\sn
\item[$(b)$]  general set theory in ZF + DC + Ax$_4$ is a worthwhile
  endeavor
\sn
\item[$(c)$]   an important reason for not adopting ZF + DC was the
  lack of something like (b), hence intellectual honesty urges you to 
investigate this direction
\sn
\item[$(d)$]  this is just a way to look at strengthening existence
  results to existence by nicely definable sets.
\end{enumerate}
\mn
Let us try to explain the results.

We assume ZF + DC.  Consider a sequence $\bar\delta = \langle \delta_s:s \in
Y\rangle$ of limit ordinals, 
when can we get a cofinal $<_I$-increasing sequence in
$(\Pi \bar\delta,<_I)$ for $I$ on ideal on $Y$?  When can we get a
parallel to the pcf-theorem?

In \cite[\S5]{Sh:938},\cite{Sh:955} we use AC$_{\cP(Y)}$ (and DC) to deal with
true pseudo cofinality, but here instead we continue \cite{Sh:835} assuming
Ax$_4$.  In \cite[1.8=L6.1]{Sh:835} we generalize the pcf-theorem
(i.e. existence of $\langle \gb_{\ga,\theta},\bar
f_{\ga,\theta}:\theta \in \text{ pcf}(\ga)\rangle$) for countable
index set $Y$.
What about large $Y$, with each $\delta_s$ having cofinality large
compared to $Y$?  Here first we deal with $D$ an $\aleph_1$-complete
filter in \ref{c13}; this continues the ideas of \cite[1.2=Lr.2]{Sh:835}.  
We then can\footnote{We temporarily cheat a little, only
  $A_\varepsilon/I_\varepsilon$ is defined.}
 choose $\langle A_\varepsilon,J_\varepsilon,\bar
f:\varepsilon < \varepsilon(*)\rangle, J_\varepsilon$ the
$\aleph_1$-complete ideal on $Y$ generated by $\{A_\zeta:\zeta <
\varepsilon\},\bar f$ cofinal in $(\Pi(\bar\delta \rest
A_\varepsilon),<_{I_\varepsilon})$ and $<_{I_\varp}$-increasing. Can we waive
``$\aleph_1$-complete"?  For this in \ref{c17} we combine the above with a
generalization of \cite[1.6=1p.4]{Sh:835}, i.e. above $I_\varepsilon$
is the ideal on $Y$ generated by $\{A_\zeta:\zeta < \varepsilon\}$.
If $I_\varepsilon$ is not $\aleph_1$-complete we deal essentially
with all quotients of $I_\varepsilon$ which are ideals on countable
sets.  

But in Theorem \ref{c17}, what about $\Pi \bar\delta$ when $s \in
Y \Rightarrow \cf(\delta_s)$ small?  With choice, recalling \cite{Sh:673}
we cannot generalize
the pcf theorem\footnote{Still by \cite{Sh:506}, in ZFC, we can deal
  with $(\Pi \bar\lambda,<_I)$ if $\lambda_s > \theta$ and a relative
  of ``$\cP(I)/I$ satisfies the $\theta$-c.c." hold.}, but here, 
even if each $\delta_s$ has countable
cofinality this is not necessarily the case.
This motivates the definition of the ideal $\cf-\id_{<
  \theta}(\bar\delta)$ noting that in general it may well be
that $s \in Y \Rightarrow \cf(\delta_s) = \aleph_0$ but $\cf(\Pi \bar
\delta)$ is large.

In our context, the set ${}^\kappa \lambda$ does not in general have a
cardinality, i.e. its power is not a cardinal, i.e. an $\aleph$, 
equivalently the set is not
well orderable.  But surprisingly, by Theorem \ref{s14} in \S(2D),
relevant covering numbers exist, i.e.
$\cov(\lambda,\theta_3(\kappa),\kappa,\sigma)$ is a well defined $\aleph$
when the cardinality of the sets by which we cover
$(< \theta_3(\kappa))$ is large enough compared to the ones we cover
$(< \kappa)$.  This is an additional witness for the covering number's  
naturality.  This follows by moreover proving when $\kappa = 
\sigma = \aleph_1$, there is a cofinal subset which is well orderable.
In particular here it
gives us a way to circumvent the non-existence of well orders of ${}^\kappa
\lambda$.  

In \S(2A), \S(2B) we deal with relatives of \S1: pcf system, eub and
more.  Also in \ref{d29} 
we give an improvement of the result of \cite[\S1]{Sh:835}.

Another issue is the ``successor of a singular cardinal is regular" in \S(2C).
Recall that the consistency strength of two successive singular
cardinal is large, but not for ``a successor cardinal is singular".
So a posteriori (i.e. after \cite[\S1]{Sh:835}) it 
is natural to hope that if $\mu$ is singular large enough
then $\mu^+$ is regular.  In \cite[2.13=Ls.2]{Sh:835} we show that
for many $\mu$ the answer is yes; here we get a stronger
conclusion: $\mu^+$ is a true successor cardinal; in fact $\alpha <
\mu \Rightarrow |\alpha|^{\aleph_0} < \mu$ suffice; see \ref{p21}(2).

Many proofs rely on diagonalizing so seemingly inherently use strong
choice.  Still we succeed to save some, see \S3.  As a 
test problem, we deal with constructing Abelian groups 
and with Black Boxes.  We also note that
\cite{Sh:460} applies even in $\ZF + \AC_{\aleph_0}$ in \ref{p38}.

A natural question is:
\mn
\begin{enumerate}
\item[$(*)$]  assume $\cf(\mu) = \aleph_0,(\forall \alpha <
\mu)(|\alpha|^{\aleph_0} < \mu)$
\sn
\begin{enumerate}
\item[$(a)$]   if $\mu \le \lambda < \mu^{\aleph_0}$ and $\lambda$ is
singular, is $\lambda^+$ a true successor? or at least
\sn
\item[$(b)$]   if $\mu \le \lambda < \pp(\mu)$ and $\lambda$ is
singular is $\lambda^+$ is regular?
\end{enumerate}
\end{enumerate}
\mn
We may try to use a closure operation 
$c \ell$ which is only $\aleph_1$-well founded, hence
have to use DC$_{\aleph_1}$.

How can we try to prove?
We may try to prove that if $\mu > 2^{\aleph_0}$ is singular
then $\lambda = \mu^+$ is regular improving \cite[2.13=Ls.2]{Sh:835}, where
there are further restrictions on $\mu$.  A natural approach is letting
$\chi \le \mu$ be minimal such that $\chi^{\aleph_0} \ge \mu$, so
$\chi > 2^{\aleph_0}$, so as there we can find $\bar C_1 = \langle
C_\alpha:\alpha \in S^\lambda_{< \chi}\rangle,C_\alpha \subseteq
\alpha = \sup(C_\alpha)$ and $|C_\alpha| < \chi$.  
But what about $S^\lambda_{\ge \chi}$?
Assume $\lambda = \pp(\chi)$ so we can find $\langle
\lambda_n:n < \omega\rangle$, each $\lambda_n$ is $< \chi,J$ ideal on 
$\omega,\tcf(\Pi \lambda_n,<_J) = \lambda$ and 
$\bar f = \langle f_\alpha:\alpha <
\lambda\rangle$ is $<_J$-increasing cofinal in $(\Pi \lambda_n,<_J)$.
Without loss of generality $\cf(\alpha) > 2^{\aleph_0} \Rightarrow
f_\alpha$ a $<_J$-eub of $\bar f \rest \alpha$.

Another approach is to build an AD family $\cA \subseteq
[\lambda]^{\aleph_0}$ which induces a ``good" function $c
\ell_{\cA}:\cP(\lambda) \rightarrow \cP(\lambda)$: where $c \ell_{\cA}(u) =
\cup\{A \in \cA:A \cap u$ infinite$\}$, maybe let $\cA_0$ be induced
by $\bar f$.

\noindent
Naturally we may ask (and deal with some, as mentioned).
\begin{question}
\label{d15}
1) Can we bound $\hrtg(\cP(\mu))$ for $\mu$ singular?
(recall Gitik-Koepke \cite[pg.2]{GiKo10pr}).

\noindent
2) Can we deduce \wlor$({}^Y \mu) = \hrtg({}^Y \mu)$ when $\mu$ is
singular large enough?  Maybe see \cite[Ld21]{Sh:F1303}.

\noindent
3) In \S1 we may replace $\theta$ by several $\theta_\ell$, defined by
the proof (i.e. $\theta_\ell$ is minimal satisfying some demands
involving $\theta_0,\dotsc,\theta_{\ell-1}$ and the pcf problem); 
but seemingly this does not make a serious gain, maybe see
on this in \cite[5.2=Le4]{Sh:F1303}. 

\noindent
4) Can we generalize RGCH (see \cite{Sh:460}, \cite[\S1]{Sh:829}), see
 \ref{p38}, \ref{s21}.  Maybe see more in \cite{Sh:F1438}.

We thank the referee for checking the paper very carefully discovering
many things which should be mended much above the call of duty.
\end{question}
\bigskip

\subsection {Preliminaries} \label{0B}\
\bigskip

\begin{hypothesis}
\label{z2}
1) We work in $\ZF + \DC$.

\noindent
2) Usually we assume $\Ax_{4,\partial}$, see Definition \ref{z6}(5)
relying on \ref{z4}(3), \ref{z6}(4),
so a reader may assume it throughout; or even assume $\Ax_4$, see
\ref{z4}(2),(1).  
Many times we use weaker
relatives so we try to mention the case of 
$\Ax_{4,\lambda,\theta,\partial}$ actually used.
So the case $\theta = \partial = \aleph_1$ means $\Ax_{4,\lambda}$
holds and note $\Ax_4$ is stronger than $\Ax_{4,\aleph_1}$.

\noindent
3) So no such assumption means $\ZF + \DC$ \underline{but} still
   $\partial$ is a fixed cardinal $\ge \aleph_1$.
\end{hypothesis}

\begin{definition}
\label{z6}
1) \hrtg$(A) = \Min\{\alpha$: there is no function from $A$ onto $\alpha\}$.

\noindent
2) $\wlor(A) = \Min\{\alpha$: there is no one-to-one function from
$\alpha$ into $A$ or $\alpha=0 \wedge A = \emptyset\}$ so
$\wlor(A) \le \hrtg(A)$.
\end{definition}

\begin{definition}
\label{z4}
1) Ax$^4_\lambda$ means $[\lambda]^{\aleph_0}$ can be well ordered 
so $\lambda^{\aleph_0}$ is a well defined cardinal.

\noindent
2) Ax$_4$ means Ax$^4_\lambda$ for every cardinality $\lambda$.

\noindent
3) Ax$_{4,\lambda,\partial,\theta}$ means that $(\lambda \ge \partial
\ge \theta \ge \aleph_1$ and): there is a witness
   $\cS$ which means:
\mn
\begin{enumerate}
\item[$(a)$]  $\cS \subseteq ([\lambda]^{< \partial},\subseteq)$
\sn
\item[$(b)$]  for every $u_1 \in [\lambda]^{< \theta}$ there is $u_2
\in \cS$ such that $u_1 \subseteq u_2$
\sn
\item[$(c)$]  $\cS$ is well-orderable
\sn
\item[$(d)$]   for notational simplicity: $\cS$ of minimal
 cardinality.
\end{enumerate}
\mn
3A) But we may use an ordinal $\beta$ instead of $\lambda$ above.
So trivially $\Ax^4_\lambda \Rightarrow
\Ax_{4,\lambda,\aleph_1,\aleph_1}$ because we can choose
$\cS = [\lambda]^{\le \aleph_0}$.

\noindent
3B) If $\Ax_{4,\lambda,\partial,\theta}$ then we let
$\cov(\lambda,\partial,\theta,2)$ be the minimal $|\cS|$ for $\cS$ as
in \ref{z4}(3); necessarily it is $< \wlor([\lambda]^{< \partial})$ 
which is $\le \hrtg([\lambda]^{< \partial})$; so if $\neg
\Ax_{4,\lambda,\partial,\theta}$ then it is not well defined.

\noindent
3C) We say $(\cS_*,<_*)$ witness $\Ax_{4,\lambda,\partial,\theta}$
  \when \, $\cS_*$ is as in part (3) and $<_*$ is a well ordering of
  $\cS_*$.

\noindent
3D) We say $(\cS_*,<_*)$ witness $\Ax^4_\lambda$ when $\cS_* =
[\lambda]^{\le \aleph_0}$ and $<_*$ is a well ordering of
$[\lambda]^{\le \aleph_0}$.

\noindent
4) Let Ax$_{4,\lambda,\partial}$ mean
   Ax$_{4,\lambda,\partial,\aleph_1}$; note that even if $\partial =
   \aleph_1,\Ax_{4,\lambda,\partial}$ is not $\Ax^4_\lambda$.

\noindent
5) Let Ax$_{4,\partial}$ mean Ax$_{4,\lambda,\partial}$ for every
$\lambda$, so $\Ax_{4,\partial}$ is not the same as $\Ax^4_\partial$.

\noindent
6) We may write $\le \theta$ instead of $\theta^+$, and writing an
   ordinal $\alpha$ instead of $\partial$ means $\otp(u_1) < \alpha$
   in clause (b) of part (3); similarly for the other parameters. 
\end{definition}

\noindent
We try to make the paper reasonably self-contained.  Still we assume
knowledge of \cite[\S(0B)]{Sh:835}, the preliminaries, in particular,
recall:
\begin{claim}
\label{z8}
1) For every $\lambda,\partial$ satisfying 
$\Ax_{4,\lambda,\partial}$ there is a function
$c \ell$, moreover one which is (we may use $\alpha$ instead of $\lambda$)
definable from $(\cS_*,<_*)$ where $(\cS_*,<_*)$ 
witness $\Ax_{4,\lambda,\partial}$, see \ref{z4}(3),(3B), even
uniformly such that:
\mn
\begin{enumerate}
\item[$(a)$]  $c \ell:\cP(\lambda) \rightarrow \cP(\lambda)$
\sn
\item[$(b)$]  $u \subseteq c \ell(u) \subseteq \lambda$, (but we do not
  require $c \ell(c \ell(u)) = c \ell(u))$
\sn
\item[$(c)$]  $|c \ell(u)| < \hrtg([u]^{\aleph_0} \times \partial)$
\sn
\item[$(c)'$]  if $\Ax_4$ and $(\cS_*,<_*)$ witness it then $|c \ell(u)|
\le |u|^{\aleph_0}$ for $u \subseteq \lambda$
\sn
\item[$(d)$]  there is no sequence $\langle u_n:n < \omega\rangle$
such that $u_{n+1} \subseteq u_n \nsubseteq c \ell(u_{n+1})$.
\end{enumerate}
\mn
2) We can above replace $\Ax_{4,\lambda,\partial}$ by: there is a well
orderable $\cS_* \subseteq [\lambda]^{< \partial}$ such that there is no
$u \in [\lambda]^{\aleph_0}$ satisfying $v \in \cS_* \Rightarrow
\aleph_0 > |v \cap u|$.
\end{claim}

\begin{PROOF}{\ref{z8}}
1) Recall $\cS_* \subseteq [\lambda]^{< \partial}$ and $u_1 \in
[\lambda]^{\le \aleph_0} \Rightarrow (\exists u_2 \in \cS_*)(u_1
\subseteq u_2)$ and $<_*$ is a well ordering of $\cS_*$ and let
$\langle w^*_i:i < \otp(\cS_*,<_*)\rangle$ list $\cS_*$ in
$<_*$-increasing order; if $\Ax_4$ we can use $\cS_* =
[\lambda]^{\aleph_0}$.  For $v \in [\lambda]^{\le \aleph_0}$ let
$\mathbf i(v) = \mathbf i(v,\cS_*,<_*) = \min\{i:v \backslash w^*_i$ is
finite$\}$.

For $u \subseteq \lambda$ let $c \ell(u) = \cup\{w^*_i$: for some $v
\in [u]^{\aleph_0}$ we have $i = \mathbf i(v)\} \cup u \cup\{0\}$.

So clearly clauses (a),(b) of the conclusion hold.

For clause (c) define $F:[u]^{\aleph_0} \times \partial 
\rightarrow \lambda$ by
$F(v,\alpha) =$ the $\alpha$-th member of $w^*_{\mathbf i(v)}$ when
$\otp(w^*_{\mathbf i(v)}) > \alpha$, and 0 otherwise; clearly $F$ is a
function from $[u]^{\aleph_0} \times \partial$ to $\lambda$
and its range is included in $c \ell(u)$ and includes $c \ell(u)
\backslash u$; we like $F$ to be onto $c \ell(u)$, but clearly $u
\backslash \Rang(F)$ is finite, hence
this last part can be corrected easily
 hence $c \ell(u)$ has cardinality $< \hrtg([u]^{\aleph_0} \times \partial)$
so we are done with clause (c).

Lastly, to prove clause (d), toward contradiction assume 
$\bar u = \langle u_n:n < \omega\rangle$ and $u_{n+1}
\subseteq u_n \nsubseteq c \ell(u_{n+1})$ for every $n$; by $\DC$ or
just $\AC_{\aleph_0}$ choose
$\bar\alpha = \langle \alpha_n:n < \omega\rangle$ such that $\alpha_n
\in u_n \backslash c \ell(u_{n+1})$.  Now let $v = \{\alpha_n:n < \omega\}$
and $i = \mathbf i(v)$, so for every $n,v \backslash (v \cap u_n)$ is
finite hence $\mathbf i(v) = \mathbf i(v \cap u_n)$ and let $n$ be such
that $v \backslash w^*_i \subseteq \{\alpha_0,\dotsc,\alpha_{n-1}\}$,
so $\alpha_n \in w^*_i \subseteq c \ell(u_{n+1})$, contradicting the
choice of $\alpha_n$.

\noindent
2) Similarly but first for any infinite $v \subseteq \lambda$ let
$\mathbf i(v) = \mathbf i(v,\cS_*,<_*) := \min\{i:v
\cap w^*_i$ is infinite$\}$.  Second, $F(v,\alpha)$ is:
\mn
\begin{enumerate}
\item[$\bullet$]  the $\alpha$-th member of $w^*_{\mathbf i(v)}$ if
$\alpha < \otp(w^*_{\mathbf i(v)})$
\sn
\item[$\bullet$]  $0$ otherwise.
\end{enumerate}
\mn
Third, note:
\mn
\begin{enumerate}
\item[$\bullet$]  if $u \subseteq \lambda$ then $u \backslash
  \{F(v,\alpha):v \in [u]^{\aleph_0}$ and $\alpha < \partial\}$ is
  finite.
\end{enumerate}
\mn
[Why?  If not, let the difference be $v^*$ and let $v = \{\alpha \in
v^*:v^* \cap \alpha$ is finite$\}$ so $v$ is a subset of 
the difference of cardinality $\aleph_0$, (infinite by our assumption), hence 
$\{F(v,\alpha):\alpha < \lambda\}$ is not disjoint to $v$, contradiction.]

Fourth, in the end, instead of ``let $n$ be such that 
$v \backslash w^*_i \subseteq \{\alpha_0,\dotsc,\alpha_{n-1}\}"$ we choose
$n$ such that $\alpha_n \in w^*_{\mathbf i(v)} \cap v$; possible as
$w^*_{\mathbf i(v)} \cap v = w^*_{\mathbf i(v)} \cap \{\alpha_n:n <
\omega\}$ is infinite and $n < \omega \Rightarrow \mathbf i(v) = \mathbf
i(\{\alpha_k:k > n\})$.
\end{PROOF}

\begin{observation}
\label{z5}
1) For any set $Y$, if $\mu$ a cardinal and $\theta := \hrtg(Y)$ \then \,
$\hrtg(Y \times \mu) \le (\theta + \mu)^+$.

\noindent
2) In \ref{z8} we can replace clause (c) by:
\mn
\begin{enumerate}
\item[$(c)'$]  $|c \ell(u)| < \max\{\partial^+,\hrtg([u]^{\aleph_0})\}$.
\end{enumerate}
\end{observation}

\begin{PROOF}{\ref{z5}}
1) Assume $F$ is a function from $Y \times \mu$ onto an ordinal
   $\gamma$.

For $\beta < \mu$ let $v_\beta = \{F(y,\beta):y \in Y\}$, so $\langle
v_\beta:\beta < \mu\rangle$ is a well defined sequence of subsets of
the ordinal $\gamma$ with union $\gamma$, and clearly $\beta < \mu
\Rightarrow |v_\beta| < \hrtg(Y) =
\theta$.  Really we can use $v'_\beta = v_\beta \backslash
\cup\{v_\alpha:\alpha < \beta\}$, in this case 
clearly $\langle v'_\beta:\beta <
\mu\rangle$ is a partition of $\gamma$.
Hence easily $|\gamma| = |\bigcup\limits_{\beta < \mu} v_\beta| =
|\bigcup\limits_{\beta < \mu} v'_\beta| \le \theta + \mu$, 
so the desired result follows.

\noindent
2) Let $\theta = \hrtg([u]^{\aleph_0})$, if $\theta \le \partial$ then
   applying part (1), $\hrtg([u]^{\aleph_0} \times \partial) \le (\theta
+ \partial)^+ = \partial^+$ so we are done.  If $\theta >
   \partial$, then $\hrtg([u]^{\aleph_0} \times \partial) \le
   \hrtg([u]^{\aleph_0} \times [u]^{\aleph_0})$ and if $|u| \ge
\aleph_0$ we have $|[u]^{\aleph_0} \times [u]^{\aleph_0}| =
   |u|^{\aleph_0}$ hence we are done.

Lastly, if $\neg(|u| \ge \aleph_0)$ then (as $u \subseteq \lambda$)
necessarily $u$ is finite and so $c \ell(u) = u \cup \{0\}$ hence $|c
\ell(u)| < \partial$, so having covered all cases we are done.
\end{PROOF}

\begin{convention}
\label{z9}
1) Let ``there is $y$ satisfying $\psi(y,a),\partial$-uniformly definable
(or uniformly $\partial$-definable) for $a \in A$" means that there 
is a formula $\varphi(x,y,z)$ such that:
\mn
\begin{enumerate}
\item[$\bullet$]  for every $\mu$ large enough
if $a \in A$ and $\Ax_{4,\mu,\partial}$ holds and $<_*$ well orders 
some $\cS_* \subseteq [\mu]^{< \partial}$ as in \ref{z4}(3)
then $(\exists!y)[\varphi(y,a,<_*) \wedge \psi(y,a)]$.
\end{enumerate}
\mn
1A) Note that it follows that there is a definable
  function $A \mapsto \mu_A \in \card$ such that above, $\mu \ge
  \mu_A$ suffice.

2) Similarly with $(\partial,\theta)$-uniformly definable \when \, we use
$\Ax_{4,\mu,\partial,\theta}$; and $(\mu,\partial,\theta)$-uniformly
definable when we fix $\mu$.

\noindent
3) If the parameter $(\partial)$ or $(\partial,\theta)$ or
$(\mu,\partial,\theta)$ is clear from the context 
we may omit it.  We may not always remember to state this.

\noindent
4) $\delta$ denotes an ordinal, limit one if not said otherwise.
\end{convention}

\begin{definition}
\label{z10}
Let $D$ be a filter on a set $Y$.

\noindent
1) For $\bar\delta \in {}^Y\Ord$ let
$\lambda = \tcf(\Pi \bar\delta,<_D)$ means that $(\Pi \bar\delta,<_D)$
has true cofinality $\lambda$, i.e. $\lambda$ is a regular cardinal and
there is a witness that is a $<_D$-increasing
sequence $\langle f_\alpha:\alpha <\lambda\rangle$ of members of $\Pi
\bar\delta$ which is cofinal in $(\Pi \bar\delta,<_D)$; but sometimes
we allow $\lambda$ to be an ordinal so not unique.  (Why helpful?  See
part (2)).

\noindent
2) We say that $\bigwedge\limits_{i \in I} \lambda_i = \tcf(\Pi
\bar\delta_i,<_D)$ \when \, $\bar\delta_i \in {}^Y \Ord$ for $i \in I$
and there is a sequence $\big< \langle f^i_\alpha:\alpha < \lambda_i \rangle:i
\in I \big >$ such that $\langle f^i_\alpha:\alpha < \lambda_i\rangle$
is as above for $\lambda_i = \tcf(\Pi \bar\delta_i,<_D)$, 
but $\lambda_i$ may be any ordinal hence is not unique; so
$\bigwedge\limits_{i \in I} \lambda_i = \tcf(\Pi \bar\delta_2,<_D)$
and $i \in I \Rightarrow \lambda_i = \tcf(\Pi \bar\delta_i,<_D)$ 
has a different meaning.

\noindent
3) Assume $\bar f = \langle f_\alpha:\alpha < \delta\rangle$ and
$\alpha < \delta \Rightarrow f_\alpha \in {}^Y \Ord$ and $D$ is a
filter on $Y$.  We say $f \in {}^Y\Ord$ is a $<_D-\eub$ of $\bar f$ 
\when \,:
\mn
\begin{enumerate}
\item[$(a)$]   $\alpha < \delta \Rightarrow f_\alpha \le f \mod D$
\sn
\item[$(b)$]   if $g \in {}^Y\Ord$ and
$(\forall s \in Y)(g(s) < f(s) \vee g(s)=0)$ then $(\exists \alpha <
 \delta)(g \le f_\alpha \mod D)$.
\end{enumerate}
\end{definition}

\begin{definition}
\label{z11}
1) Let $Y$ be the set and let $\kappa$ be an infinite cardinal.
\mn
\begin{enumerate}
\item[$(a)$]    $\Fil^1_\kappa(Y)$ is the set of $\kappa$-complete
filters on $Y$, (so $Y$ is defined from $D$ as $\cup\{X:X \in D\}$)
\sn
\item[$(b)$]   $\Fil^2_\kappa(Y) = \{(D_1,D_2):D_1 \subseteq D_2$ are
$\kappa$-complete filters on $Y$, ($\emptyset \notin D_2$, of course)$\}$;
 in this context $Z \in \bar D$ means $Z \in D_2$
\sn
\item[$(c)$]   $\Fil^3_\kappa(Y,\mu) = \{(D_1,D_2,h):(D_1,D_2) \in
\Fil^2_\kappa(Y)$ and $h:Y \rightarrow \alpha$ for some
$\alpha < \mu\}$, if we omit $\mu$ we mean $\mu = \hrtg([Y]^{\le
  \aleph_0} \times \partial) \cup \omega$, recalling \ref{z2}
\sn
\item[$(d)$]    $\Fil^4_\kappa(Y,\mu) = \{(D_1,D_2,h,Z):
(D_1,D_2,h) \in \Fil^3_\kappa(Y,\mu)$ and $Z \in D_2\}$; omitting
$\mu$ means as above. 
\end{enumerate}
\mn
2) For ${\gy} \in \Fil^4_\kappa(Y,\mu)$ let 
$Y = Y^{\gy} = Y_{\gy},{\gy} = (D^{\gy}_1,D^{\gy}_2,h^{\gy},Z^{\gy}) =
(D_{\gy,1},D_{\gy,2},h_{\gy},Z_{\gy}) 
= (D_1[{\gy}],D_2[{\gy}],h[{\gy}],Z[{\gy}])$; similarly for the
others and let $D^{\gy} = D[\gy]$ be $D^{\gy}_1 + Z^{\gy}$ recalling
$D+Z$ is the filter generated by $D \cup \{Z\}$.

\noindent
3)If $\kappa = \aleph_1$ we may omit it.

\noindent
4) For $D$ a filter on $Y$ and $f \in {}^Y \Ord$ we define $\rk_D(f)
\in \Ord \cup \{\infty\}$ by $\rk_D(f) = \sup\{\rk_D(g)+1:g \in {}^Y
\Ord$ and $g <f \mod D\}$, (the Galvin-Hajnal rank).
\end{definition}

\noindent
We now repeat to a large extent \cite{Sh:835}, \cite{Sh:938}
\begin{dc}
\label{z13}
Assume $\delta$ is a limit ordinal (or zero for some parts), 
$D = D_1 \in \text{ Fil}^1_{\aleph_1}(Y),\bar f = \langle
f_\alpha:\alpha < \delta\rangle$ is a sequence of members\footnote{We
  can use any index set instead of $\delta$ (in particular the empty
  one), \underline{except} in part (5); this
  applies also to Definition \ref{z10}.} of ${}^Y \Ord$, usually
 $<_{D_1}$-increasing in
${}^Y\Ord,f$ is a $\le_D$-upper bound of $\bar f$ but there
is no such $g <_D f$; necessarily there is such $f$ (using DC).

\noindent
1) [Definition]  Let $J = J[f,\bar f,D] := 
\{A \subseteq Y$: either $A = \emptyset$ mod $D$ 
\underline{or} $A \in D^+$ but there is a $\le_{D+A}$-upper bound $g
 <_{D+A} f$ of $\bar f\}$.

\noindent
2) $J[f,\bar f,D]$ is an $\aleph_1$-complete ideal on $Y$  
disjoint to $D$.

\noindent
3) [Definition]  Recalling $D_1=D$, let 
$D_2 = D_2(f,\bar f,D_1) = \text{ dual}(J[f,\bar f,D_1]) :=
\{A \subseteq Y:Y \backslash A \in J[f,\bar f,D_1]\}$; note that, e.g. as
$D_1$ is $\aleph_1$-complete then $D_2$ is an 
$\aleph_1$-complete filter on $Y$ extending $D_1$.

\noindent
4) In (3), $f$ is a unique modulo $D_2$, i.e. if also $g \in
   {}^Y\Ord$, is a $<_{D_1}$-upper bound of $\bar f$ and
$J[g,\bar f,D_1] = J[f,\bar f,D_1]$ \then \, $g = f \mod D_2$,
   equivalently $\mod J[f,\bar f,D_1]$.

\noindent
5) If ($\bar f$ is $\le_{D_1}$-increasing, and)
$\cf(\delta) \ge \hrtg(\cP(Y))$ \then \, $f$ from above is a $<_{D_2}$-eub of
$\bar f$, see Definition \ref{z10}(3).
\end{dc}

\begin{definition}
\label{z15}
Assume $f \in {}^Y \Ord,D_2 \supseteq D_1$ are $\aleph_1$-complete
filters on $Y,c \ell$ is as in \ref{z8} for $\alpha(*)$
 and $\Rang(f) \subseteq \alpha(*)$.

\noindent
0) For some $\gy \in \Fil^4_{\aleph_1}(Y),D^{\gy}_1 = D_1,D^{\gy}_2 =
D_2$ and the function $f$ satisfies $\gy$, see below.

\noindent
1) We say $f:Y \rightarrow \Ord$ satisfies $\gy \in
\Fil^4_{\aleph_0}(Y)$ \when \,:
\mn
\begin{enumerate}
\item[$(a)$]  if $Z \in D^{\gy}_2$ and $Z \subseteq Z_{\gy}$ then $c
  \ell(\{f(t):t \in Z\}) = c \ell(\{f(t):t \in Z_{\gy}\}$
\sn
\item[$(b)$]   $y \in Z_{\gy} \Rightarrow h_{\gy}(y) = \otp(f(y) \cap c
 \ell(\Rang(f \rest Z_{\gy})))$
\sn
\item[$(c)$]   if $t \in Y$ and $f(t) \in c \ell\{f(s):s \in
  Z_{\gy}\}$ then $t \in Z_{\gy}$
\sn
\item[$(d)$]  $y \in Y \backslash Z_{\gy} \Rightarrow f(y) = 0$.
\end{enumerate}
\mn
2) ``Semi satisfies" mean we omit clause (d).

\noindent
3) Let ``weakly satisfies" means we omit clauses (c),(d).
\end{definition}

\begin{definition}
\label{z17}
Let $Y,f,\bar f,D$ be as in \ref{z13} and $Y,\alpha(*),c \ell$ as in
\ref{z15}.

\noindent
1) We say $f$ is the $(\gy,c \ell)$-eub of $\bar f$ \underline{or}
canonical $\bar f$-eub for $\gy$ and $c \ell$ or for $(\gy,c \ell)$
\when \,:
\mn
\begin{enumerate}
\item[$(a)$]  $\gy \in \Fil^4_{\aleph_1}(Y)$
\sn
\item[$(b)$]  $\bar f = \langle f_\alpha:\alpha < \alpha_*\rangle$
\sn
\item[$(c)$]  $f_\alpha,f$ are from ${}^Y \alpha(*)$
\sn
\item[$(d)$]  $f_\alpha \le_{D_{\gy,1}} f$
\sn
\item[$(e)$]  $D_{\gy,1}=D$ and $D_{\gy,2} \supseteq \dual(J[f,\bar
  f,D_{\gy,q}])$
\sn
\item[$(f)$]  $f$ satisfies $\gy$ (for $c \ell$).
\end{enumerate}
\end{definition}

\begin{claim}
\label{z19}
Let $Y,f,\bar f,D$ as in \ref{z13}, $f,\alpha(*),c \ell$ as in \ref{z15}.

\noindent
1) The ``the" is \ref{z17} is justified, that is, $f$ is unique given
$c \ell$ (so $\alpha(*),\bar f,\gy$).

\noindent
2) There is one and only one $\gy$ such that
\mn
\begin{enumerate}
\item[$(a)$]  $\gy \in \Fil^4_{\aleph_1}(Y)$
\sn
\item[$(b)$]  $D_{\gy,1} = D$
\sn
\item[$(c)$]  $D_{\gy,2} = \dual(J[f,\bar f,D])$
\sn
\item[$(d)$]  $f$ semi satisfies $\gy$.
\end{enumerate}
\mn
3) For the $\gy$ from part (2), letting $g = (f \rest Z_{\gy}) \cup
(0_{Y \backslash Z_{\gy}})$ we have $g$ is the canonical $\bar f-\eub$
for $\gy$ (and $c \ell$), in particular it satisfies $y$.
\end{claim}

\begin{PROOF}{\ref{z19}}
Should be clear.
\end{PROOF}

\noindent
Recall the related (not really used)
\begin{dc}
\label{z21}
Assume $D \in F^1_{\aleph_1}(Y)$ and $f:Y \rightarrow \Ord$.

\noindent
1) [Definition] $J[f,D] = \{A \subseteq Y:A = \emptyset \mod D$ or $A
\in D^+$ and $\rk_{D+A}(f) > \rk_D(f)\}$.

\noindent
2) $J$ is an $\aleph_1$-complete filter disjoint to $D$.

\noindent
3) If $f_1,f_1:Y \rightarrow \Ord$ and $J[f_1,D] = J[f_2,D]$.

\noindent
4) There is one and only $\gy \in \Fil^4_{\aleph_1}(Y)$ such that $f$
semi satisfies $\gy,D_{\gy,1} = D$ and $D_{\gy,2} = \dual(J[f,D])$.

\noindent
5) In (4) there is a unique $f'$ which satisfies $\gy$ and $f' \rest
Z_{\gy} = f \rest Z_{\gy}$.
\end{dc}

\begin{notation}
\label{z22}
Let $A \le_{\qu} B$ means that $A = \emptyset$ or there is a function
from $B$ onto $A$.
\end{notation}

\begin{observation}
\label{z24}
Assume $\partial \le |Y|$ and even $\partial \subseteq Y$ 
for transparency.

\noindent
1)  $\Fil^4_{\aleph_1}(Y) \le_{\qu}|\cP(\cP(3 \times Y))|$.

\noindent
2) Also ${}^\omega(\Fil^4_{\aleph_1}(Y)) \le_{\qu} \cP(\cP(Y))$.

\noindent
3) If $\theta = \hrtg(\cP(\cP(Y))$ \then \, $\theta$ satisfies:
\mn
\begin{enumerate}
\item[$\bullet$]  if $\alpha < \theta$ then $\hrtg(\cP([\alpha]^{\aleph_0}
  \times \partial)) \le \theta$
\sn
\item[$\bullet$]  so if $\Ax_4$ then $|\alpha|^{\aleph_0} \times
  \partial < \theta.$
\end{enumerate}
\mn
4) Assume $\Ax_4$. If $\alpha < \hrtg(\cP(Y))$ then
$|\alpha|^{\aleph_0} < \hrtg(\cP(Y))$; hence if $\partial \le |Y|$ and
$\alpha < \hrtg(\cP(Y))$ then $|\alpha|^{\aleph_0} \times \partial <
\hrtg(\cP(Y))$. 
\end{observation}

\begin{remark}
If $Y$ is a set of ordinals, infinite to avoid trivialities \then \, $|Y
\times 3| = |Y|$, justifying this see \ref{c59}.
\end{remark}

\begin{PROOF}{\ref{z24}}
1) Let $Y_0 = Y,Y_{\ell +1} = \cP(Y_\ell)$ for $\ell=0,1$ and let $Y^*_1 = 
[Y_1]^{\le \aleph_0},Y^*_2 = \cP(Y^*_1),Y'_0 = 3 \times Y$ and $Y'_{\ell +1}
= \cP(Y'_\ell)$ for $\ell=0,1$
\mn
\begin{enumerate}
\item[$(*)_1$]  $|Y_0| +1 = |Y_0|$ and even $|Y_0| + \partial = |Y_0|$.
\end{enumerate}
\mn
[Why?  As $\partial \le |Y|$ is an infinite cardinal.]
\mn
\begin{enumerate}
\item[$(*)_2$]  $|Y_1| = \partial \times |Y_1|$ and 
$\partial \times |Y^*_1| = |Y^*_1|$ and $|Y'_1| = |Y'_1
\times \partial| = \partial \times |Y'_1|$.
\end{enumerate}
\mn
[Why?  Both follow by $(*)_1$.]
\mn
\begin{enumerate}
\item[$(*)_3$]  $|Y_2| \times |Y_2| = |Y_2|$ and $|Y_0| \le |Y_1| \le
  |Y_2|$ and $|Y'_2| \times |Y'_2| = |Y'_2|$; moreover (for part (2))
  $|{}^\omega(Y_2)| = |Y_2|$ and $|{}^\omega(Y'_2)| = |Y'_2|$
\end{enumerate}
\mn
[Why?  Follows by $(*)_2$.]
\mn
\begin{enumerate}
\item[$(*)_4$]  $\{D_{\gy,\ell}:\gy \in \Fil^4_{\aleph_1}(Y)\}$ has power
  $\le |Y_2|$ for $\ell=1,2$.
\end{enumerate}
\mn
[Why?  By the definition each $D_{\gy,\ell}$ is a subset of $\cP(Y)
  =\cP(Y_0) = Y_1$.]
\mn
\begin{enumerate}
\item[$(*)_5$]  $\{Z_{\gy}:\gy \in \Fil^4_{\aleph_1}(Y)\}$ has power 
$\le |Y_1|$.
\end{enumerate}
\mn
[Why?  As $Z_{\gy} \subseteq Y = Y_0$ so $Z_{\gy} \in Y_1$.]
\mn
\begin{enumerate}
\item[$(*)_6$]  $[Y]^{\aleph_0} \times \partial$ has the same power as
$[Y]^{\le \aleph_0}$.
\end{enumerate}
\mn
[Why?  Let $Z$ be a set of ordinals disjoint to $Y$ of order type
  $\partial$; by $(*)_1$ we have $|Y| = |Y \cup Z|$ hence
$|[Y]^{\le \aleph_0} = |[Y \cup Z]^{\le \aleph_0}| 
\ge |[Y]^{\le \aleph_0} \times 
[\partial]^{\le \aleph_0}| \ge |[Y]^{\le \aleph_0} \times \partial|
\ge [Y]^{\le \aleph_0}$.]
\mn
\begin{enumerate}
\item[$(*)_7$]  $|Y \times [Y]^{\aleph_0} \times [Y]^{\aleph_0}| \le
 |\cP(3 \times Y)| \le |Y_2|$.
\end{enumerate}
\mn
[Why?  The mapping $(y,u_1,u_2) \mapsto \{(0,y),(1,z_1),(2,z_2):z_1 \in
u_1,z_2 \in u_2\}$ from $Y \times [Y]^{\aleph_0} \times
[Y]^{\aleph_0}$ into $\cP(3 \times Y)$ prove the first inequality, the
second inequality follows from $|3 \times Y| = |3 \times Y_0| \le 
|Y'_1| = |Y_1|$.]
\mn
\begin{enumerate}
\item[$(*)_8$]  $\cH := \{h_{\gy}:\gy \in \Fil^4_{\aleph_1}(Y)\} 
\le_{\qu} |Y'_2|$.
\end{enumerate}
\mn
[Why?  Recalling $(*)_6$ clearly $|\cH| \le |\{h:h$ a function, $\Dom(h)=Y$ 
and $\Rang(h)$ a bounded subset of 
$\hrtg([Y]^{\le \aleph_0} \times \partial)\}| 
\le |\{h:h$ a function from $Y$ into some 
$\alpha < \hrtg([Y]^{\le \aleph_0})\}| \le_{\qu} |X_1|$ where 

\begin{equation*}
\begin{array}{clcr}
X_1 := \{(h,g): &\text{ for some ordinal } \alpha,g 
\text{ is a partial function from } [Y]^{\le \aleph_0} \text{ onto }
\alpha,\\
  &\text{ so necessarily } \alpha < \hrtg([Y]^{\le \aleph_0})\text{ and }
h \text{ is a function from } Y \text{ into } \alpha\}.
\end{array}
\end{equation*}

\mn
Clearly $|\cH| \le |X_1|$.  Let $t \notin Y$ and for 
$(h,g) \in X_1$ let set$(h,g) := \{(y,u_1,u_2):y=t \wedge \{u_1,u_2\}
\subseteq [Y]^{\le \aleph_0} \wedge g(u_1) \le
g(u_2)$ \underline{or} $y \in Y$ and
$u_1,u_2 \in [Y]^{\le \aleph_0}$ satisfies $h(y) = g(u_1)$ and
$g(u_2) = g(u_1)\}$.  Easily $(h,g) \mapsto \set(h,g)$ is a one-to-one
function from $X_1$ into $X_3 := \cP(X_2)$ where $X_2 := (Y \cup \{t\})
 \times [Y]^{\le \aleph_0} \times [Y]^{\le \aleph_0}$ and by $(*)_7$
 we have $|X_2| = |\cP(3 \times Y)|$.  Hence $|X_1| \le |X_3| =
|\cP(X_2)| \le |\cP(\cP(3 \times Y))|$.  Recalling $|\cH| \le |X_1|$ 
we are done proving $(*)_8$.]

Now $|\Fil^4_{\aleph_1}(Y)| \le |\Fil^1_{\aleph_1}(Y) \times
\Fil^1_{\aleph_1}(Y) \times \cH \times \cP(Y)|$ by the definition of
$\Fil^4_{\aleph_1}$ and this is, by the inequalities above 
$\le_{\qu} |Y'_2| \times |Y'_2| \times |Y'_2| \times |Y'_1| 
\le_{\qu} |Y'_2|^4 = |Y'_2|$.

\noindent
2),3),4) Should be clear.
\end{PROOF}

\noindent
Note also we may wonder about the RGCH, see \cite{Sh:460}, we note
(not using any version of $\Ax_4$), that we can get such a result using
only $\AC_{\aleph_0}$.  From the results of \S1 we can deduce
more. see \ref{s21}.
\begin{theorem}
\label{p38}
[$\ZF + \AC_{\aleph_0}$]  Assume that $\mu > \aleph_0$ and $\chi < \mu
\Rightarrow \hrtg(\cP(\chi)) < \mu$.
\Then \, for every $\lambda > \mu$ for some
$\kappa < \mu$ we have:
\mn
\begin{enumerate}
\item[$(*)_{\lambda,\mu,\kappa}$]  if $\theta \in (\kappa,\mu)$ and
  $D$ is a $\kappa$-complete filter on $\theta$ \then \, there is no
  $<_D$-increasing sequence $\langle f_\alpha:\alpha <
  \lambda^+\rangle$ of members of ${}^\theta \lambda$.
\end{enumerate}
\end{theorem}

\begin{remark}
In \ref{p38} we can replace ``$\chi < \mu \Rightarrow \hrtg(\cP(\chi))
< \mu$" by $\chi < \mu \Rightarrow \wlor(\cP(\chi)) < \mu$; this holds
by the proof.
\end{remark}

\begin{PROOF}{\ref{p38}}
Assume that this fails for a given $\lambda$.  
We choose $\kappa_n < \theta_n < \mu$ by induction on
$n$.  Let $\kappa_0 = \aleph_0$, so $\kappa_0 = \aleph_0 < \mu$ as
required.  Assume $\kappa_n < \mu$ has been chosen, note that it 
cannot be as required so there is $\theta \in [\kappa_n,\mu)$ 
such that it exemplifies $\neg(*)_{\lambda,\mu,\kappa_n}$ and let
$\theta_n$ be the first such $\theta$.

Given $\theta_n$ let 
$\kappa_{n+1} := \wlor(\cP(\theta_n))$ so $\kappa_{n+1} \in
(\theta_n,\mu) \subseteq (\kappa_n,\mu)$.  So $\langle \kappa_n:n <
\omega\rangle$ is well defined increasing and 
$\mu_* = \sum\limits_{n} \kappa_n \le \mu$.  
Let $X_n = \{(\theta,D,\bar f):\theta \in
   [\kappa_n,\kappa_{n+1}),D$ is a $\kappa_n$-complete filter on
     $\theta,\bar f = \langle f_\alpha:\alpha < \lambda^+\rangle$ is a
$<_D$-increasing sequence of members of ${}^\theta \lambda\}$, so
by the construction we have $X_n \ne \emptyset$ and $\langle X_n:n <
\omega\rangle$ exist being well defined.  As we are
assuming $\AC_{\aleph_0}$ there is a sequence 
$\langle (\theta_n,D_n,\bar f_n):n < \omega\rangle$ from 
$\prod\limits_{n} X_n$.

We can consider $\bar f = \langle \bar f_n:n < \omega\rangle$ (and
also $\bar\kappa = \langle \kappa_n:n < \omega\rangle$) as a set
of ordinals (using a pairing function on the ordinals) hence $\mathbf V_* =
\mathbf L[\bar f,\bar\kappa]$ is a model of ZFC and a transitive class.
In $\mathbf V_*$ we can define $D'_n$ as the minimal $\kappa_n$-complete
filter on $\theta _n$ such that $\bar f_n$ is $<_{D'_n}$-increasing.
Clearly $(2^{\theta_n})^{\mathbf V_*} < \wlor(\cP^{\mathbf V}(\theta_n)) <
\mu$ hence $\mathbf V_* \models ``\mu_*$ is strong limit".  By
\cite{Sh:460} or see \cite[\S1,1.13=Lg.8]{Sh:829} where
$\lambda^{[\partial,\theta]}$ is defined we get a contradiction.
\end{PROOF}
\newpage

\section {The pcf theorem again} \label{1}

We prove a version of the pcf theorem; weaker than
\cite[Ch.I,II]{Sh:g} as we do not assume just min$\{\cf(\alpha_y):y
\in Y\} > \hrtg(Y)$ but a stronger inequality.  Still we gain in a
point which disappears under $\AC$: dealing with a sequence of possibly
singular ordinals (and the ideal $\cf-\id_{< \theta}(\bar\delta)$, see below).
In addition we gain in having the scales being uniformly definable.
Also the result is stronger than in
\cite{Sh:955}, as we use functions rather than sets of functions;
(i.e. true cofinality rather than pseudo true cofinality; of
course, the axioms of set theory used are different accordingly;
full choice in \cite{Sh:g}, ZF + DC + AC$_{\cP(Y)}$ in \cite{Sh:955}
and ZF + DC + Ax$_4$ here).
\bigskip

\centerline {$* \qquad * \qquad *$}
\bigskip

It seems natural in our context instead of looking 
at $\{\cf(\delta_s):s \in Y\}$ we should look at:
\begin{definition}
\label{c2}
1) For a sequence $\bar\delta = \langle \delta_s:s \in Y\rangle$ of limit
ordinals and a cardinal $\theta$ let cf-id$_{< \theta}(\bar \delta)
= \{X \subseteq Y$: there is a sequence $\bar u = \langle u_s:s \in
Y\rangle$ such that $s \in X \Rightarrow u_s \subseteq \delta_s =
\sup(u_s)$ and $s \in X \Rightarrow \otp(u_s) <  \theta\}$.

\noindent
2) Let $\cf-\fil_{<\theta}(\bar\delta)$ be the filter dual to the ideal 
$\cf-\id_{< \theta}(\bar\delta)$.

\noindent
3) We may replace $\bar\delta$ by a set of ordinals, i.e. instead of
   $\langle \alpha:\alpha \in u\rangle$ we may write $u$.

\noindent
4) For $\bar\delta = \langle \delta_s:s \in Y\rangle$ and $\bar\theta
   = \langle \theta_s:s \in Y\rangle$ we define 
$\cf-\id_{< \bar\theta}(\bar\delta)$ similarly to 
part (1); similarly in the other cases.

\noindent
5) For $\bar\theta$ a sequence of infinite cardinals, let 
$\cf-\fil_{< \bar\theta}(\bar\delta)$ be the dual filter;
   similarly in the other cases.
\end{definition}

\begin{observation}
\label{c5}
1) In \ref{c2}, $\cf-\id_{< \theta}(\bar\delta),
\cf-\id_{< \bar\theta}(\bar\delta)$ are ideals on $Y$ or equal to $\cP(Y)$.

\noindent
1A) Moreover $\aleph_1$-complete ideals.

\noindent
2) Similarly for the filters.
\end{observation}

\begin{PROOF}{\ref{c5}}
Should be clear, e.g. use the definitions recalling we are assuming
$\AC_{\aleph_0}$.
\end{PROOF} 

\begin{observation}
\label{c10}
Assume
\mn
\begin{enumerate}
\item[$(a)$]  $D = \cf-\fil_{< \bar\theta}(\bar\delta)$ is a well
  defined filter (that is $\emptyset \notin D$), so 
$\bar\delta \in {}^Y\Ord$ is a sequence of limit
ordinals, $\bar\theta = \langle \theta_s:s \in Y\rangle \in {}^Y\Car$,
e.g. $\bigwedge\limits_{s} \theta_s = \theta$
\sn 
\item[$(b)$]  $\bar\cU = \langle \cU_s:s \in Y\rangle$ satisfies
$\cU_s \subseteq \delta_s,\otp(\cU_s) < \theta_s$ for $s \in Y$,
\sn
\item[$(c)$]  $g \in \Pi \bar\delta$ is defined by
\sn
\begin{enumerate}
\item[$\bullet$]  $g(s)$ is $\sup\{\alpha +1:\alpha \in \cU_s\}$ if
this value is $< \delta_s$
\sn
\item[$\bullet$]  $g(s)$ is zero otherwise.
\end{enumerate}
\end{enumerate}
\mn
\Then
\mn
\begin{enumerate}
\item[$(\alpha)$]   $g$ belongs to $\Pi \bar\delta$ indeed
\sn
\item[$(\beta)$]  if $f \in \prod\limits_{s \in Y} \cU_s \subseteq \Pi
\bar\delta$ then $f < g \mod D$.
\end{enumerate}
\end{observation}

\begin{remark}
Clause (b) of \ref{c10} holds, e.g. if $\cU \subseteq \Ord,\otp(\cU) < 
\min\{\theta_s:s \in Y\},\cU_s = \cU \cap \delta_s$.
\end{remark}

\begin{PROOF}{\ref{c10}}
Clause $(\alpha)$ is obvious by the choice of the function $g$; 
for clause $(\beta)$ let $f \in
\prod\limits_{s \in Y} \cU_s$ and let $X = \{s \in Y:f(s) \ge
g(s)\}$.  Necessarily $s \in X$ implies (by the assumption on $f$ and
the definition of $X$) that 
$(\exists \alpha)(\alpha \in \cU_s \wedge g(s) \le \alpha)$ which
implies (by clause (c), the definition of $g$) that $g(s) = 0 \wedge
\sup(u_s) = \delta_s$.  So by the definition of $\cf-\fil_{<
\bar\theta}(\bar\delta)$ we have $X \in 
\cf-\fil_{<\bar\theta}(\bar\delta)$ hence we are done.
\end{PROOF}

\begin{claim}  
\label{c13}
Assume $\Ax_{4,\partial}$, see Definition \ref{z4}(3); if 
(A) then (B) \underline{where}:
\mn
\begin{enumerate}
\item[$(A)$]  we are given $Y$, an arbitrary set, $\bar\delta$, a
  sequence of limit ordinals and $\mu$, an infinite cardinal (or just
  a limit ordinal) such that:
\sn
\begin{enumerate}
\item[$(a)$]  $\bar\delta = \langle \delta_s:s \in Y\rangle$
  and $\mu = \sup\{\delta_s:s \in Y\}$
\sn
\item[$(b)$]  $D_*$ is an $\aleph_1$-complete filter on
$Y$, it may be $\{Y\}$
\sn
\item[$(c)$]  $\theta$ is any cardinal satisfying:
\sn
\item[${{}}$]  $(\alpha) \quad \cf-\id_{< \theta}(\bar\delta)
  \subseteq \dual(D_*)$, 
\sn
\item[${{}}$]  $(\beta) \quad \alpha < \theta \Rightarrow
\hrtg([\alpha]^{\aleph_0} \times \partial) \le \theta$ so 
$\partial < \theta$ 
\sn
\item[${{}}$]  $(\gamma) \quad \hrtg(\cP(Y)) \le \theta$
\sn
\item[${{}}$]  $(\delta) \quad \hrtg(\Fil^4_{\aleph_1}(Y))\le
  \theta$
\end{enumerate}
\sn
\item[$(B)$]   there are $\alpha_*,f,\bar f,A_*/D_*$ 
$\partial$-uniformly defined 
from the triple $(Y,\bar \delta,D_*)$, see \ref{z9} such that (see more in
the proof):
\mn
\begin{enumerate}
\item[$(a)$]  $\alpha_*$ is a limit ordinal of cofinality $\ge \theta$
\sn
\item[$(b)$]  $\bar f = \langle f_\alpha:\alpha < \alpha_*\rangle$
\sn
\item[$(c)$]  $f_\alpha \in \Pi \bar\delta$ and $f \in \Pi \bar\delta$
\sn
\item[$(d)$]  $\bar f$ is $<_{D_*}$-increasing
\sn
\item[$(e)$]  $A_* \in D^+_*$
\sn
\item[$(f)$]  $\bar f$ is cofinal in $(\Pi \bar\delta,<_{D_*+A_*})$
\sn
\item[$(g)$]   if $Y \backslash A_* \in D^+_*$ then $f$ is
$a <_{D_* + (Y \backslash A_*)}$-ub of the sequence $\bar f$.
\end{enumerate}
\end{enumerate}
\end{claim}

\begin{remark}
\label{c14}
1) Note that we do not use $\AC_{\cP(Y)}$ and even not
$\AC_Y$ which would simplify.

\noindent
2) Note that $\theta$ is not necessarily regular.

\noindent
3) In $(A)(c)(\delta)$, we can restrict ourselves to
   $\aleph_1$-complete filters on $Y$ extending $D_*$.

\noindent
4) Originally we use several $\theta$'s to get best results but not
   clear if worth it.

\noindent
5) Why for a given $Y$ there is $\theta$ as in
\ref{c13}(A)(c)$(\beta),(\gamma),(\delta)$? see \ref{z24}(3).

\noindent
6) In \ref{c13} we can replace the assumption $\Ax_{4,\partial}$ by
 $\Ax_{4,\hrtg({}^Y \mu),\partial}$, see \ref{z4}(4),(5).

\noindent
7) Concernig $(A)(c)(\alpha)$ note that this
holds when each $\delta_s$ is an
ordinal $\le \mu$ of cofinality $\ge \theta$.

\noindent
7A) In $(A)(c)(\beta)$, if $\Ax_4$ then 
the demand is equivalent to ``$\partial < \theta$ and
$\alpha < \theta \Rightarrow |\alpha|^{\aleph_0} < \theta$", see \ref{z24}(4).
\end{remark}

\begin{PROOF}{\ref{c13}}
We can define $\mu$ by clause (A)(a) and \wilog \, $\theta$ is minimal
such that (A)(c) holds and recall $\partial$ is given so fixed.

Let
\mn
\begin{enumerate}
\item[$(*)_1$]  $(a) \quad \lambda_* = \hrtg({}^Y \mu)$
\sn
\item[${{}}$]  $(b) \quad \cS_{\lambda_*} \subseteq [\lambda_*]^{<
\partial}$ is as in \ref{z4}(3) 
\sn
\item[${{}}$]  $(c) \quad <_{\lambda_*}$ be a well ordering of
  $\cS_{\lambda_*}$ 
\sn
\item[${{}}$]  $(d) \quad \bar w^* = \langle w^*_i:i <
\otp(\cS_{\lambda_*},<_{\lambda_*})\rangle$ list $\cS_{\lambda_*}$ in
$<_{\lambda_*}$-increasing order
\sn
\item[$(*)_2$]  $c \ell$ be as in \ref{z8} for $\lambda_*$
\sn
\item[$(*)_3$]   $\Omega = \{\alpha < \lambda_*:\aleph_0 \le 
\cf(\alpha) < \theta\}$.
\sn
\item[$(*)_4$]  There is a sequence $\bar e$ (in fact, $\partial$-uniformly 
definable one) such that:
\sn
\begin{enumerate}
\item[$(a)$]   $\bar e = \langle e_\alpha:\alpha \in \Omega\rangle$
\sn
\item[$(b)$]  $e_\alpha \subseteq \alpha = \sup(e_\alpha)$
\sn
\item[$(c)$]  $e_\alpha$ has order type $< \theta$; 
\end{enumerate}
\end{enumerate}
\mn
and we can add
\mn
\begin{enumerate}
\item[${{}}$]  $(c)_1 \quad e_\alpha$ has order type $< \partial$ 
if \cf$(\alpha) = \aleph_0$
\sn
\item[${{}}$]  $(c)_2 \quad e_\alpha$ has cardinality $<
\hrtg([\cf(\alpha)]^{\aleph_0} \times \partial)$.
\end{enumerate}
\mn
[How?  
\mn
\begin{enumerate}
\item[$\bullet$]  If $\cf(\alpha) = \aleph_0$ let $\mathbf i(\alpha) =
\min\{i:w^*_i \cap \alpha$ is unbounded in $\alpha\}$ and $e_\alpha =
w^*_{\mathbf i(\alpha)} \cap \alpha$.
\sn
\item[$\bullet$]  If $\cf(\alpha) > \aleph_0$ let $e_\alpha = c
\ell(e)$ where $e$ is any club of $\alpha$ of order type $\cf(\alpha)$
such that $(\forall e')[e' \subseteq e$ a club of $\alpha \Rightarrow
c \ell(e') = c \ell(e)]$.
\end{enumerate}
\mn
[Why?  Such $e$ exists by the choice of $c \ell$ in \ref{z8} and if
$e'_*,e''_*$ are two such clubs then $e'_* \cap e''_*$ is a club of
$\alpha$ of orer type $\cf(\alpha)$ and $c \ell(e') = c \ell(e' \cap
e'') = c \ell(e'')$ by the assumption on $e'$ and on $e''$
respectively, so $e_\alpha$ is well defined.]

Lastly, the cardinality is as required by the clause $(A)(e)(\beta)$ 
and \ref{z8}(c); similarly to \cite[2.11=Lr.9]{Sh:835}.

So $(*)_4$ holds indeed.]

Now we try to choose $f_\alpha \in \Pi \bar\delta$ by induction on
$\alpha$ such that $\beta < \alpha \Rightarrow f_\beta < f_\alpha \mod D_*$.
\bigskip

\noindent
\underline{Case 1}:  $\alpha = 0$

Let $f_\alpha$ be constantly zero, i.e. $s \in Y \Rightarrow
f_\alpha(s) = 0$, clearly $f_\alpha \in \Pi \bar\delta$ as each
$\delta_s$ is a limit ordinal.
\bigskip

\noindent
\underline{Case 2}:  $\alpha = \beta +1$

Let $f_\alpha(s) = f_\beta(s) +1$ for $s \in Y$, so $f_\alpha \in \Pi
\bar\delta$ as $f_\beta \in \Pi \bar\delta$ and each $\delta_s$ is a
limit ordinal and $\gamma < \alpha \Rightarrow f_\gamma < f_\alpha$
mod $D_*$ as $f_\gamma \le f_\beta < f_\alpha \mod D_*$.
\bigskip

\noindent
\underline{Case 3}: $\alpha$ is a limit ordinal of cofinality $< \theta$.

So $e_\alpha$ is well defined and we define $f_\alpha:Y \rightarrow
\Ord$ as follows:  $f_\alpha(s)$ is equal to $\sup\{f_\beta(s)+1:
\beta \in e_\alpha\}$ if this is $< \delta_s$ and is zero otherwise.
\mn
\begin{enumerate}
\item[$(*)_5$]  $f_\alpha \in \Pi \bar\delta$.
\end{enumerate}
\mn
[Why?  Obvious.]

Let $\cU_{\alpha,s} = \{f_\beta(s)+1:\beta \in e_\alpha\}$, so 
clearly $\langle \cU_{\alpha,s}:s \in Y\rangle$ is well defined and
$\sup(\cU_{\alpha,s})$ is an ordinal, it is $\le \delta_s$ as $\beta
\in e_\alpha \Rightarrow f_\beta \in \Pi \bar\delta$.  Let $X = \{s
\in Y:f_\alpha(s) > 0$ equivalently $\delta_s > \sup(\cU_{\alpha,s})\}$ 
\mn
\begin{enumerate}
\item[$(*)_6$]  $X \in D_*$, i.e. $X = Y \mod D_*$.
\end{enumerate}
\mn
[Why?  For $s \in Y \backslash X$ note that $|\cU_{\alpha,s}|
  \le_{\qu} |e_\alpha|$ and $|e_\alpha| < \theta$ by $(*)_4(c)$, hence 
  $|\cU_{\alpha,s}| < \theta$.  By the choice of $X$ and Definition
\ref{c2} we have $Y \backslash X \in \cf-\id_{< \theta}(\bar\delta)$
  hence the clause (A)(c)$(\alpha)$ of the assumption of the claim, $X
  = Y \mod D_*$ as promised.]
\mn
\begin{enumerate}
\item[$(*)_7$]  if $\beta < \alpha$ then $f_\beta < f_\alpha \mod D_*$.
\end{enumerate}
\mn
[Why?  Clearly $e_\alpha$ has no last element so we can choose
$\gamma \in e_\alpha \backslash (\beta +1)$ and let $X' = \{s \in
Y:f_\beta(s) < f_\gamma(s)\}$.  Necessarily $X' \in D_*$ hence $X' \cap
X \in D_*$ but clearly $s \in X' \cap X \Rightarrow f_\beta(s) <
  f_\gamma(s) < f_\alpha(s)$ so $(*)_7$ holds.]
\bigskip

\noindent
We arrive to the main case.
\medskip

\noindent
\underline{Case 4}:  $\alpha$ a limit ordinal of cofinality $\ge \theta$

Let 
\mn
\begin{enumerate}
\item[$\bullet$]  $\bar f^\alpha = \langle f_\beta:\beta <
\alpha\rangle$
\sn
\item[$\bullet$]  $\mathbf D = \{D:D \text{ is an } 
\aleph_1 \text{-complete filter on } Y \text{ extending } D_*\}$
\sn
\item[$\bullet$]  $\mathbf D^1_\alpha = \{D \in \mathbf D:\bar f^\alpha$
is not cofinal in $(\Pi \bar\delta,<_D)\}$
\sn
\item[$\bullet$]  $\mathbf D^2_\alpha = \{D \in \mathbf D^1_\alpha:\bar
f^\alpha$  has a $<_D \text{-upper bound } f \in \Pi \bar\delta\}$
\sn
\item[$\bullet$]  $\mathbf D^3_\alpha = \{D \in \mathbf D^2_\alpha:\bar f^\alpha$
has a $<_D \text{-eub } f \in \Pi \bar\delta\}$.
\end{enumerate}
\mn
For every $D \in \mathbf D^3_\alpha$ let
\mn
\begin{enumerate}
\item[$\bullet$]  $\cF^3_{\alpha,D} = \{f \in \Pi \bar\delta:
f \text{ is a } <_D\text{-eub of } \langle f_\beta:\beta < \alpha \rangle\}$.
\end{enumerate}
\mn
Note
\mn
\begin{enumerate}
\item[$\odot_1$]  if $D_1 \in \mathbf D^1_\alpha$ and $f$ exemplifies this
\then \, for some $D_2,D_1 \subseteq D_2 \in \mathbf D$ and 
$f$ is a $<_{D_2}$-upper bound of $\bar f$, i.e. $f$ exemplifies $D_2
\in \mathbf D^2_\alpha$; in fact $D_2$ is uniformly definable from $f$
(and $\bar f^\alpha,D_1$).
\end{enumerate}
\mn
[Why?  Let $\bar A = \langle A_\gamma:\gamma < \alpha\rangle$ be
defined by $A_\gamma := \{s \in Y:f(s) \le f_\gamma(s)\}$.  So $\langle
A_\gamma/D_1:\gamma < \alpha\rangle$ is increasing (in the Boolean
algebra $\cP(Y)/D_1$, of course), but clearly
$|\{A/D_1:A \subseteq Y\}| \le_{\text{qu}} |\cP(Y)|$ and $\hrtg(\cP(Y))
\le \theta$ by clause $(A)(c)(\gamma)$ of the assumption. 
Let $\cU = \{\gamma < \alpha$: for no $\beta < \gamma$
do we have $A_\gamma = A_\beta \mod D\}$, so clearly $|\cU| <
\hrtg(\cP(Y)) \le \theta$ by $(A)(c)(\gamma)$
but by the present case assumption,
$\cf(\alpha) \ge \theta$ so $\langle A_\gamma/D_1:\gamma <
\alpha\rangle$ is necessarily eventually  constant.  Let $\alpha(*) =
\min\{\gamma$: if $\beta \in (\gamma,\alpha)$ then $A_\beta =
A_\gamma$ mod $D_1\}$; it is well defined (and $< \alpha$).  Now
$A_{\alpha(*)} \notin D_1$ as otherwise $f \le f_{\alpha(*)} <
f_{\alpha(*)+1} \mod D_1$ contradicting the assumption on $f$.  Let
$D_2 := D_1 + (Y \backslash A_{\alpha(*)})$.  Clearly $D_2$ is as
required.]
\mn
\begin{enumerate}
\item[$\odot_2$]   if $D \in \mathbf D^2_\alpha$ and $f$ exemplifies it
\then \, for some $g$ we have:
\sn
\begin{enumerate}
\item[$(a)$]  $g \in \Pi \bar\delta$
\sn
\item[$(b)$]  $g \le_D f$
\sn
\item[$(c)$]  $g$ is a $<_D$-upper
bound of $\langle f_\gamma:\gamma < \alpha\rangle$
\sn
\item[$(d)$]   there is no $h \in \Pi \bar\delta$ which is an $<_D$-upper
  bound of $\langle f_\gamma:\gamma < \alpha\rangle$ such that $h <_D g$.
\end{enumerate}
\end{enumerate}
\mn
[Why?  Use DC and $D$ being $\aleph_1$-complete.]
\mn
\begin{enumerate}
\item[$\odot_3$]   if $D_1 \in \mathbf D^2_\alpha$ and $g$ is as in
$\odot_2$  \then \, for a unique pair $(\gy,f)$ we have
\sn
\begin{enumerate}
\item[$(a)$]   $\gy \in \Fil^4_{\aleph_1}(Y)$
\sn
\item[$(b)$]   $D_{\gy,1} = D_1$
\sn
\item[$(c)$]  $D_{\gy,2} = \dual(J[g,\bar f^\alpha,D_1])$
from \ref{z13}(1)
\sn
\item[$(d)$]  $Z_{\gy}$ satisfies:
\sn
\item[${{}}$]  $(\alpha) \quad Z_{\gy} \in D_{\gy,2}$
\sn
\item[${{}}$]  $(\beta) \quad Z \in D_{\gy,2} \wedge Z 
\subseteq Z_{\gy} \Rightarrow c \ell((\Rang(g \rest Z_{\gy}) =
c \ell(\Rang(g \rest Z)$,
\sn
\item[${{}}$]  $(\gamma) \quad$ if $t \in Y$ and $g(t) \in c
  \ell(\Rang(g \rest Z_{\gy})$ then $t \in Z_{\gy}$
\sn
\item[$(e)$]  $h_{\gy}:Z_{\gy} \rightarrow \Ord$ (really into some
  $\alpha < \hrtg(\cP(Y))$ is defined
by $g(s) =$ the $h_{\gy}(s)$-th member of 
$c \ell(\Rang(g \rest Z_{\gy}))$ if $s \in Z_{\gy}$ and
\sn
\item[$(f)$]  $f:Y \rightarrow \Ord$ is defined by $f \rest Z_{\gy} =
  g \rest Z_{\gy}$ and $f(s)=0$ for $s \in Y \backslash Z_{\gy}$.
\end{enumerate}
\end{enumerate}
\mn
[Why?  We apply \ref{z19}(2) with $g,\langle f_\gamma:\gamma <
\alpha\rangle$ here standing for $f,\bar f$ there to define $\gy$ and
then let $f = (g \rest Z_{\gy},0_{Y \backslash Z_{\gy}})$ as in
  \ref{z19}(3).]

In particular, the ``the" in $\odot_3(c)$ is justified by:
\mn
\begin{enumerate}
\item[$\odot'_3$]  if $\gy \in \Fil^4_{\aleph_1}(Y)$ and $f',f''$ are
$(\gy,c \ell)-\eub$ of $\bar f^\alpha$ \then \, $f' = f''$, i.e. \ref{z19}(3).
\end{enumerate}
\mn
Also, (recalling $\dom(f') = \dom(f'') = Z_{\gy}$ by
$\odot_3,\delta)$, see \ref{z13}(4))
\mn
\begin{enumerate}
\item[$\odot''_3$]   if $\gy \in \Fil^4_{\aleph_1}(Y)$ and
  $f',f''$ satisfy $\odot_3(e)$ \then \, $f' = f'' \mod D_{\gy,2}$.
\end{enumerate}
\mn
Recalling \ref{z13}(5), let
\mn
\begin{enumerate}
\item[$\odot_4$]  $\gY^2_\alpha = \{\gy \in \Fil^4_{\aleph_1}(Y):
D_* \subseteq D_{\gy,1} \text{ and some } f \in
{}^{Z[\gy]}\Ord$ is a $\gy$-eub of $\bar f^\alpha\}$ 
\sn
\item[$\odot_5$]  for each $\gy \in \gY^2_\alpha$, let $f_{\gy} =
f^2_{\alpha,\gy}$ be the unique function $f \in \Pi(\bar\delta 
\rest Z_{\gy})$ which is the canonical $\gy$-eub of $\langle f_\gamma:\gamma <
  \alpha\rangle$.
\end{enumerate}
\mn
Now let
\mn
\begin{enumerate}
\item[$\odot_6$]  for $s \in Y$ let $\cU^*_{\alpha,s} = 
\{f_{\gy}(s):\gy \in \gY^2_\alpha\}$.
\end{enumerate}
\mn
Clearly
\mn
\begin{enumerate}
\item[$\odot_7$]  $(a) \quad \langle \cU^*_{\alpha,s}:s \in Y\rangle$
  is well defined
\sn
\item[${{}}$]  $(b) \quad \cU^*_{\alpha,s} \subseteq \delta_s$
\sn
\item[${{}}$]  $(c) \quad$ if $s \in Y$ then $|\cU^*_{\alpha,s}| < \theta$.
\end{enumerate}
\mn
[Why?  Clause (a) holds by $\odot_6$ and clause (b) by $\odot_5 +
  \odot_6$.  As for clause (c) by
  $\odot_6, \cU^*_{\alpha,s}$ is the range of the function
$\gy \mapsto f_{\gy}(s)$ for $\gy \in \gY^2_\alpha,s \in
Z_{\gy}$, so clearly $|\cU^*_{\alpha,s}| \le_{\qu}
  |\gY^2_\alpha| \le_{\qu} |\Fil^4_{\aleph_1}(Y)|$ hence
$|\cU^*_{\alpha,s}| < \hrtg(\Fil^4_{\aleph_1}(Y))$ which is $\le
  \theta$ by (A)(c)$(\delta)$ of the claim.]
\mn
\begin{enumerate}
\item[$\odot_8$]  $X := \{s \in Y:\sup(\cU^*_{\alpha,s}) < \delta_s\} =
  Y \mod \cf-\id_{< \theta}(\bar\delta)$ hence $X \in D_*$. 
\end{enumerate}
\mn
[Why?  By $\odot_7(a),(b),(c)$ and Definition \ref{c10} we have $X
  = Y \mod \cf-\id_{< \theta}(\bar\delta)$ but by (A)(c)$(\alpha)$, this
  implies $X \in D_*$.]

So define $f_\alpha \in \Pi \bar\delta$ by:
\mn
\begin{enumerate}
\item[$\odot_9$]  $f_\alpha(s) \begin{cases} \text{ is }
  \sup(\cU^*_{\alpha,s}) \quad &\text{ if } s \in X \\
  \text{ is } 0 \quad &\text{ if } s \in Y \backslash X \end{cases}$
\end{enumerate}
\mn
Also clearly
\mn
\begin{enumerate}
\item[$\odot_{10}$]  $f_\alpha \in \Pi \bar\delta$
\end{enumerate}
\mn
and also
\mn
\begin{enumerate}
\item[$\odot_{11}$]  if $\gy \in \gY^2_\alpha$ and
  $\beta < \alpha$ \then \, $f_\beta < f_\alpha \mod D_{\gy,2}$.
\end{enumerate}
\mn
For $\beta < \alpha$ let $A^\alpha_\beta = \{s \in
Y:f_\beta(s) < f_\alpha(s)\}$ so $\bar A^\alpha = \langle
A^\alpha_\beta:\beta < \alpha\rangle$ is well defined and $\langle
A^\alpha_\beta/D_*:\beta < \alpha\rangle$ is decreasing (in the
Boolean Algebra $\cP(Y)/D_*$) and is eventually constant as
$\hrtg(\cP(Y)/D_*) \le \hrtg(\cP(Y)) \le \theta$ by clause
(A)(c)$(\gamma)$ of the assumption so let $\gamma(\alpha) =
\min\{\gamma < \alpha$: for every $\beta \in (\gamma,\alpha)$ we have
$A^\alpha_\beta/D_* = A^\alpha_\gamma/D_*\}$.

If $A^\alpha_{\gamma(\alpha)} \in D_*$ then $\beta < \alpha
\Rightarrow A^\alpha_\beta \supseteq
A^\alpha_{\max\{\beta,\gamma(\alpha)\}} = A^\alpha_{\gamma(\alpha)}
\mod D_* \Rightarrow f_\beta < f_\alpha \mod D_*$, so 
$f_\alpha$ is as required. 
Otherwise, $A^\alpha_{\gamma(\alpha)} \notin D_*$, so $A_* :=
Y \backslash A^\alpha_{\gamma(\alpha)} \in D^+_*$ so $D_1 = D_* + A_* \in
\mathbf D$.  Now if $D_1 \in \mathbf D^1_\alpha$ then by $\odot_1$ there is
$D_2$ such that $D_1 \subseteq D_2 \in \mathbf D^2_\alpha$ hence there
is $g \in \Pi \bar\delta$ as in $\odot_2$ for $D_2$ hence there is
$\gy \in \Fil^4_{\aleph_1}(Y)$ as in $\odot_3$ hence $f_{\gy} \in
\Pi(\bar\delta)$ as in $\odot_5$, so $Z_{\gy} \in
D_{\gy,2}$, and by the choice of $\cU^*_{\alpha,s}(s \in Y)$ and
$f_\alpha$ we have $f_{\gy} \le f_\alpha \mod D_{\gy,2}$ hence $\beta
< \alpha \Rightarrow f_\beta < f_\alpha \mod D_{\gy,2}$ so
$f_{\gamma(\alpha)} < f_\alpha \mod D_{\gy,2}$.  But $A_* \in D_1 =
D_{\gy,1} \subseteq D_{\gy,2}$ and by the choice of
$A^\alpha_{\gamma(\alpha)}$ and $A_*$ we have $f_\alpha \rest A_* \le
f_{\gamma(\alpha)} \rest A_*$ contradicting the previous sentence.

So necessarily $(A_* \in D^+_*$ and) $D_1 = D_* + A_* \in \mathbf D$
does not belong to $\mathbf D^1_\alpha$ which means $\bar f^\alpha$ is
cofinal in $(\Pi\bar\delta,<_{D_* + A_*})$ hence
letting the desired $(\alpha_*,f,\bar f,A_*/D_*)$ in (B) of \ref{c13}  
be $(\alpha,f_\alpha,\bar f^\alpha,A_*/D_*)$ we are done.
\end{PROOF}

\begin{theorem}
\label{c17}
\underline{The pcf Theorem}: 
[$\Ax^4_{4,\theta,\partial},\theta =\hrtg({}^Y \mu) + \DC$]

If $(A)$ then $(B)^+$ where:
\mn
\begin{enumerate}
\item[$(A)$]  we\footnote{Clause (A) here is as 
in \ref{c13}(A) but $D_*$ is just 
a filter on $Y$, not necessarily $\aleph_1$-complete filter on 
$Y$ (i.e. we weaken clause (b)), noting that possibly 
$D_* = \{Y\}$, still we require $\cf-\fil_{< \theta}(\bar\delta) 
\subseteq D_*$.}
 are given $Y$, an arbitrary set, $\bar\delta$, a
 sequence of limit ordinals and $\mu$, an infinite cardinal (or just
 a limit ordinal) and $D_*,\theta$ such that
\sn
\begin{enumerate}
\item[$(a)$]  $\bar\delta = \langle \delta_s:s \in Y\rangle$
  and $\mu = \sup\{\delta_s:s \in Y\}$
\sn
\item[$(b)$]  $D_*$ is an $\aleph_1$-complete filter\footnote{This is
    reasonable as we normally use $D_* = \dual(\cf-\id_{<
      \theta}(\bar\delta))$ which is $\aleph_1$-complete by \ref{c10}(1A).}
 on $Y$, it may be $\{Y\}$
\sn
\item[$(c)$]  $\theta$ is any cardinal satisfying:
\sn
\item[${{}}$]  $(\alpha) \quad \cf-\id_{< \theta}(\bar\delta)
  \subseteq \dual(D_*)$, note that this
holds when each $\delta_s$ is an

\hskip25pt  ordinal $\le \mu$ of cofinality $\ge \theta$, see below
\sn
\item[${{}}$]  $(\beta) \quad \alpha < \theta \Rightarrow
\hrtg([\alpha]^{\aleph_0} \times \partial) \le \theta$ so 
$\partial < \theta$ 
\sn
\item[${{}}$]  $(\gamma) \quad \hrtg(\cP(Y)) \le \theta$
\sn
\item[${{}}$]  $(\delta) \quad \hrtg(\Fil^4_{\aleph_1}(Y))\le
  \theta$
\end{enumerate}
\sn
\item[$(B)^+$]  there are $\varepsilon(*),\bar D^*,\bar A^*,\bar E^*,
\bar\alpha^*,\bar g$, in fact $\partial$-uniformly definable from 
$(Y,\bar\delta,D_*)$ such that:
\sn
\begin{enumerate}
\item[$(a)$]  $\varepsilon(*) < \hrtg(\cP(Y))$
\sn
\item[$(b)$]  $\bar D^* = \langle D^*_\varepsilon:\varepsilon \le
\varepsilon(*)\rangle$ and $\bar E^* = \langle
E^*_\varepsilon:\varepsilon < \varepsilon(*)\rangle$
\sn
\item[$(c)$]  $\bar D^*$  is a $\subset$-increasing continuous sequence
(of filters on $Y$, but see $(f)$) 
\sn
\item[$(d)$]  if $\varepsilon = \zeta +1$ \then \,
$D^*_\varepsilon$ is a filter on $Y$ generated by $D_\zeta \cup
\{A\}$ for some $A \subseteq Y$ such that $A \in D^+_\zeta$
\sn
\item[$(e)$]  $D^*_0 = D_*$ 
\sn
\item[$(f)$]  $D^*_\varepsilon$ is a filter on $Y$ for $\varepsilon <
  \varepsilon(*)$ but $D^*_{\varepsilon(*)} = \cP(Y)$, 
\sn
\item[$(g)$]  $(\alpha) \quad \bar\alpha^* = 
\langle \alpha^*_\varepsilon:\varepsilon \le \varepsilon(*)\rangle$
\sn
\item[${{}}$]  $(\beta) \quad \bar\alpha^*$ is an increasing
continuous sequence of ordinals
\sn
\item[${{}}$]  $(\gamma) \quad \alpha^*_0 = 0,
\cf(\alpha^*_{\varepsilon +1}) \ge \theta$
\sn
\item[${{}}$]  $(\delta) \quad \varepsilon(*)$ is a
  successor ordinal
\sn 
\item[$(h)$]  $\bar g = \langle g_\alpha:\alpha <
\alpha^*_{\varepsilon(*)}\rangle$ is a sequence of members of $\Pi
\bar\delta$, so of functions from $Y$ into the ordinals
\sn
\item[$(i)$]  if $\beta < \alpha < \alpha_{\varepsilon +1}$ \then \,
$g_\beta < g_\alpha \mod D^*_\varepsilon$ 
\sn
\item[$(j)$]  $\bar A^* = \langle
  A^*_\varepsilon/D^*_\varepsilon:\varepsilon < \varepsilon(*)\rangle$
  where $A^*_\varepsilon \subseteq Y$, 
so only $A^*_\varepsilon/D^*_\varepsilon$ is computed\footnote{But see
  \ref{d12}.}  not $A^*_\varepsilon$, still 
$(Y \backslash A^*_\varepsilon)/D^*_\varepsilon$ and
$D^*_\varepsilon + (Y \backslash A^*_\varepsilon)$ are well defined
\sn
\item[$(k)$]  $D^*_{\varepsilon +1} = D^*_\varepsilon +
A^*_\varepsilon$ and $E^*_\varepsilon = D^*_\varepsilon 
+ (Y \backslash A^*_\varepsilon)$ if $\varepsilon$ is a successor
ordinal and $D_\varepsilon$ if otherwise
\sn
\item[$(l)$]  $\langle g_\alpha:\alpha \in
[\alpha^*_\varepsilon,\alpha^*_{\varepsilon +1})\rangle$ is increasing
and cofinal in $(\Pi \bar\delta,<_{E_\varepsilon})$ 
so also $\bar g \rest \alpha^*_{\varepsilon +1}$ is.
\end{enumerate}
\end{enumerate}
\end{theorem}

\begin{remark}
\label{c19}
1) Note that unlike the ZFC case, the $\alpha^*_{\varepsilon +1}$'s
(and even $\alpha^*_{\varepsilon +1} - \alpha^*_\varepsilon$) are
ordinals rather than regular cardinals and we do not exclude here
$\varepsilon < \zeta \wedge \cf(\alpha^*_{\varepsilon +1}) =
\cf(\alpha^*_{\zeta +1})$.  Also we do not know that $\langle
\cf(\alpha^*_\varepsilon):\varepsilon < \varepsilon(*)\rangle$ 
is increasing or even non-decreasing.

\noindent
2) We may get $\langle \alpha^*_{\varepsilon +1} -
   \alpha^*_\varepsilon:\varepsilon < \varepsilon(*)\rangle$
   non-decreasing but this is of unclear value.  [For this we proceed
as below but when we arrive to $\varepsilon +1$ and there is $\zeta <
   \varepsilon$ such that $\alpha^*_{\varepsilon +1} -
\alpha^*_\varepsilon < \alpha^*_\zeta - \alpha^*_{\zeta +1}$, choose
   the first one, we go back, retaining only $\bar g \rest
   \alpha^*_\zeta$.  Now we try again to choose $g'_\alpha$ for
 $\alpha \ge \alpha^*_\zeta$ but demanding $g'_{\alpha^*_\zeta +
   \beta} \ge g_{\alpha^*_\varepsilon + \beta},g_{\alpha^*_\zeta +
   \beta}$.  This process converges.]

\noindent
3) However \ref{c54}(5) below is a simpler way.  Working harder we get
   $\langle \alpha^*_{\varepsilon +1} -
   \alpha^*_\varepsilon:\varepsilon < \varepsilon(*)\rangle$ is
   (strictly) increasing (using increasing rectangles of functions).

\noindent
4) As in $(*)_4$ of the proof of \ref{c13}, \wilog \,
$\alpha^*_{\varepsilon +1} - \alpha^*_\varepsilon <
\hrtg([\cf(\alpha^*_{\varepsilon +1})]^{\aleph_0}) =
([\cf(\alpha^*_{\varepsilon +1}))^{\aleph_0})^+$.

\noindent
[Why?  As have first chosen $\langle g'_\alpha:\alpha \in
(\alpha^*_\varepsilon,\alpha'_{\varepsilon +1}]\rangle$ and just as
$\langle g_\alpha:\alpha \in
(\alpha^*_\varepsilon,\alpha^*_{\varepsilon +1}]\rangle$ was chosen
before we choose $\langle g_\alpha:\alpha \in (\alpha^*_\varepsilon,
\alpha^*_{\varepsilon +1}]\rangle$ by ($e_\alpha$ as in $(*)_4$ of the
proof even if $\alpha \in \lambda_* \backslash \Omega$)
\mn
\begin{enumerate}
\item[$\bullet$]  $\alpha^*_{\varepsilon +1} = \alpha^*_\varepsilon +
\otp(e_{\alpha'_{\varepsilon +1}} \backslash (\alpha^*_\varepsilon + 1))$
\sn
\item[$\bullet$]  if $\beta \in e_{\alpha'_{\varepsilon +1}}
  \backslash (\alpha^*_\varepsilon + 1)$ and $\gamma = \otp(e_\alpha
  \cap \beta) \backslash (\alpha^*_\varepsilon +1))$ then $g_\gamma = g'_\beta$
\sn
\item[$\bullet$]  if $\beta = \alpha^*_{\varepsilon +1}$ then $g_\beta
  = g'_{\alpha'_{\varepsilon + 1}}$.
\end{enumerate}
\mn
So we are done.

\noindent
5) Concerning $(\beta)$ of \ref{c17}$(B)^+(e)$, recall that 
$D_*$ include $\cf-\fil_{< \theta}(\bar\delta)$ by $(A)(c)(\alpha)$.

\noindent
6) Concerning \ref{c17}$(B)^+(f)$, if $D^*_{\varepsilon(*)}
   = \cP(Y)$ then it is not really a filter.

\noindent
7) Concerning \ref{c17} $(B)^+(i)$, note that
 using this clause in Definition \ref{d2}(2) we mean only $\le$!, that is
 we may have
\mn
\begin{enumerate}
\item[$(B)^+$]  $(i)' \quad$ if $\beta < \alpha < \alpha_{\varepsilon
    +1}$ then $g_\beta \le g_\alpha \mod D^*_\varepsilon$.
\end{enumerate}
\end{remark}

\begin{PROOF}{\ref{c17}}
Let $\cS_{\lambda_*},<_{\lambda_*} \langle w^*_i:i <
\otp(\cS_{\lambda_*},<_{\lambda_*})\rangle$ as well as $\bar e$ 
be as in the proof of \ref{c13}.

We try to choose $(\alpha^*_\varepsilon,\bar g \rest (\alpha^*_\varepsilon
+1),\bar D^\varepsilon),\bar D^\varepsilon = \langle D^*_\xi:\xi <
\varepsilon\rangle$ by induction on $\varepsilon < \hrtg(\cP(Y))$ 
such that the relevant parts of $(B)^+$ holds, but if $\emptyset \in
D^*_\varepsilon$ then $g_{\alpha^*_\varepsilon}$ is not well defined, 
so $\bar g^\varepsilon = \bar g \rest \alpha^*_\varepsilon = \langle 
g_\alpha:\alpha < \alpha^*_\varepsilon\rangle$ and $\langle
A^*_\zeta/D^*_\zeta:\zeta < \varepsilon\rangle$ are determined.
Clearly the induction has to stop before $\hrtg(\cP(Y))$, otherwise
the sequence $\langle A_\zeta/D^*_\zeta:\zeta < \hrtg(\cP(Y))\rangle$ 
gives a contradiction to the definition of $\hrtg(\cP(Y))$.
\bigskip

\noindent
\underline{Case A}:  $\varepsilon = 0$

Let $\alpha^*_\varepsilon = 0,D^*_\varepsilon = D_*$ and $g_0$ is
constantly zero.
\bigskip

\noindent
\underline{Case B}:  $\varepsilon$ a limit ordinal 

Let $\alpha^*_\varepsilon = \cup\{\alpha^*_\zeta:\zeta <
\varepsilon\},D^*_\varepsilon = \cup\{D^*_\zeta:\zeta < \varepsilon\}$
and $\bar g \rest \alpha^*_\varepsilon$ is naturally defined and
define $g_{\alpha^*_\varepsilon} \in \Pi \bar\delta$ by, for $s \in Y$
letting $g_{\alpha^*_\varepsilon}(s) = 
\cup\{g_{\alpha^*_\zeta}(s)+1:\zeta <
\varepsilon\}$ if it is $< \delta_s$ and 0 otherwise.  As in Case 3 of
the proof of \ref{c13}, clause $(B)^+(i)$ is satisfied, 
because $\hrtg(\cP(Y)) > \varepsilon$. 
\bigskip

\noindent
\underline{Case C}:  $\varepsilon = \zeta +1$ and $\emptyset \notin D^*_\zeta$.

Let (note that 
$A_{\mathbf a,n}$ in (b) below is almost equal to $Y \backslash A^*_{\xi_n}$
but we know only $A^*_{\xi_n}/D^*_{\xi_n})$:
\mn
\begin{enumerate}
\item[$(*)_1$]  $(a) \quad \mathbf J_{\zeta,1} = \{A \subseteq Y:A \in
(D^*_\zeta)^+$ and $D^*_\zeta + A$ is $\aleph_1$-complete$\}$
\sn
\item[${{}}$]  $(b) \quad \mathbf U_\zeta = \{\mathbf a:\mathbf a = \langle
(A_n,\xi_n):n < \omega\rangle = \langle (A_{\mathbf a,n},
\xi_{\mathbf a,n}):n < \omega\rangle$,

\hskip35pt  for every $n < \omega$ we have

\hskip35pt $\xi_n < \zeta$ and $D^*_{\xi_n+1} =
D^*_{\xi_n} + (Y \backslash A_n)$ and 

\hskip35pt $A_{\mathbf a} := \cup\{A_n:n < \omega\} \ne \emptyset 
\mod D^*_\zeta\}$;

\hskip35pt so this concerns witnesses to $D^*_\zeta$ being not
$\aleph_1$-complete and 

\hskip35pt $A_{\mathbf a} \in D^+_\zeta \subseteq \cP(Y)$
\sn
\item[${{}}$]  $(c) \quad \mathbf J_{\zeta,2} = \{A \subseteq Y:A \in
(D^*_\zeta)^+$ and for some $\mathbf a \in \mathbf U_\zeta$ we have $A
  \subseteq A_{\mathbf a}\}$.
\end{enumerate}
\mn
Note
\mn
\begin{enumerate}
\item[$(*)_2$]  $(a) \quad \mathbf J_{\zeta,1} \cup \mathbf J_{\zeta,2} \subseteq
(D^*_\zeta)^+$ is dense, i.e. if $A \in (D^*_\zeta)^+$ then for some
$B \subseteq A$,

\hskip25pt we have $B \in \mathbf J_{\zeta,1} \cup \mathbf J_{\zeta,2}$
\sn
\item[${{}}$]  $(b) \quad$ if $\ell \in \{1,2\},A \in \mathbf
  J_{\zeta,1},B \subseteq A$ and $B \in D^+_\zeta$ \then \, $B \in
  \mathbf J_{\zeta,\ell}$.
\end{enumerate}
\mn
[Why Clause (a)?  Because we are assuming that $D_*$ is
$\aleph_1$-complete in (A)(b).  For clause (b), just read the
definition of $\mathbf J_{\zeta,\ell}$.]

Now we try to choose $f_\alpha$ (or pedantically $f^\varepsilon_\alpha$ if you
like) by induction on $\alpha$ such that:
\mn
\begin{enumerate}
\item[$(*)_3$]  $(a) \quad f_\alpha \in \Pi \bar\delta$
\sn
\item[${{}}$]  $(b) \quad \beta < \alpha^*_\zeta \Rightarrow g_\beta <
f_\alpha \mod D^*_\zeta$; follows by (c) + (d)
\sn
\item[${{}}$]  $(c) \quad \beta < \alpha \Rightarrow f_\beta <
f_\alpha \mod D^*_\zeta$
\sn
\item[${{}}$]  $(d) \quad f_0 = g_{\alpha^*_\zeta}$.
\end{enumerate}
\mn
Arriving to $\alpha,\bar f = \langle f_\beta:\beta < \alpha \rangle$
has been defined.  Let $\mathbf J^*_{\zeta,\alpha} 
= \{A \subseteq Y:A \in (D^*_\zeta)^+$ and $\bar f$ has
an upper bound in $(\Pi \bar\delta,<_{D^*_\zeta + A})\}$.
\medskip

\noindent
\underline{Sub-case C1}:  $(\mathbf J_{\zeta,1} \cup 
\mathbf J_{\zeta,2}) \cap \mathbf J^*_{\zeta,\alpha}$ is
dense in $((D^*_\zeta)^+,\supseteq)$.

First, as in the proof of \ref{c13}, (that is, choosing $f_\alpha$ in
the inductive step in the proof) we can define $\bar f^1_{\zeta,\alpha}$
such that:
\mn
\begin{enumerate}
\item[$(*)_4$]  $(a) \quad \bar f^1_{\zeta,\alpha} = 
\langle f^1_{\zeta,\alpha,A}:A \in \mathbf J_{\zeta,1} 
\cap \mathbf J^*_{\zeta,\alpha}\rangle$
\sn
\item[${{}}$]  $(b) \quad f^1_{\zeta,\alpha,A} \in \Pi \bar\delta$
\sn
\item[${{}}$]  $(c) \quad f^1_{\zeta,\alpha,A}$ is a $<_{D^*_\zeta +A}$-upper
bound of $\{g_{\alpha^*_\zeta}\} \cup \{f_\beta:\beta < \alpha\}$.
\end{enumerate}
\mn
Second, we consider $\mathbf a \in \mathbf U_\zeta$ hence $A_{\mathbf a}
\in \mathbf J_{\zeta,2}$.  

Let
\mn
\begin{enumerate}
\item[$\bullet$]  for $u \subseteq \alpha^*_\zeta$ let 
$g^{[u]} \in \Pi \bar\delta$ be defined by
$g^{[u]}(s) = \sup\{g_\beta(s)+1:\beta \in u\}$ 
if this supremum is $< \delta_s$ and 0 otherwise.
\end{enumerate}
\mn
Note that 
\mn
\begin{enumerate}
\item[$(*)_5$]   $(a) \quad$ if $A \subseteq Y,A = 
\emptyset \mod D^*_\zeta$ \then \, for every $f \in \Pi \bar\delta$ 
for some finite 

\hskip25pt $v \subseteq \alpha^*_\zeta$ we have 
$\{s \in A: \neg(\exists \beta \in v)(f(s) < g_\beta(s))\} 
= \emptyset \mod D_*$
\sn
\item[${{}}$]  $(b) \quad$ if $u_1 \subseteq u_2$ are from
$[\alpha^*_\zeta]^{< \partial}$ then $g^{[u_1]} \le g^{[u_2]} \mod D_*$.
\end{enumerate}
\mn
[Why?  By induction on $\zeta$ using $(B)^+(k),(l)$ 
recalling $D^*_0 = D_*$ or see the proof of \ref{c59}.   
Clause (b) is proved by $\cf-\id_{< \partial}(\bar\delta) \subseteq
\cf - \id_{< \theta}(\bar\delta) \subseteq \dual(D_*)$ recalling
$\aleph_0 < \theta$ by $(A)(c)(\beta)$ of the claim.]  
\mn
\begin{enumerate}
\item[$(*)_6$]  if $f \in \Pi \bar\delta$ \then \, $f$ 
has a $<_{D^*_\zeta + A_{\mathbf a}}$-upper bound and even a
$<_{D_* +A_{\mathbf a}}$-upper bound of the form $g^{[u]}$ 
for some countable $u \subseteq \alpha^*_\zeta$.
\end{enumerate}
\mn
[Why?  Let $f \in \Pi \bar\delta$, now for each $n$ there is $\alpha_n
< \alpha^*_{\xi_{\mathbf a,n}+1}$ such that $f < g_{\alpha_n}
\mod(D^*_{\xi_{\mathbf a,n}} + A_{\mathbf a,n})$, moreover, see
$(*)_5(a)$, there is a finite set $v_n \subseteq \alpha^*_{\xi_{\mathbf
a,n}+1}$ such that $(\forall s \in A_{\mathbf a,n})(\exists \beta \in
v)(f(s) < g_\beta(s))$.  As those are finite sets of ordinals (or use
$\AC_{\aleph_0}$) there is such a sequence
$\langle v_n:n < \omega\rangle$, so $u = \cup\{v_n:n <\omega\}$ is as
required, recalling $\cf-\id_{< \aleph_1}(\bar\delta) \subseteq
\dual(D_*)$ as in earlier cases so we have proved (a) of $(*)_5$.]

Lastly (well defined by $(*)_5(b) + (*)_6$ recalling our sub-case assumption):
\mn
\begin{enumerate}
\item[$(*)_7$]  let $\bar f^2_{\zeta,\alpha} = \langle 
f^2_{\zeta,\alpha,\mathbf a}:\mathbf a \in \mathbf U_\zeta\rangle$ be defined by: $f^2_{\zeta,\alpha,\mathbf a}$ 
is $g^{[u]}$ where $u = u_{\mathbf a} \in \cS_{\lambda_*}$ is the
$<_{\lambda_*}$-first $u \in \cS_{\lambda_*}$ for which $g^{[u]}
\in \Pi \bar\delta$ is a $<_{D^*_\zeta + A_{\mathbf a}}$-common upper
bound of $\{g_{\alpha^*_\zeta}\} \cup \{f_\beta:\beta < \alpha\}$.
\end{enumerate}
\mn
Note that
\mn
\begin{enumerate}
\item[$(*)_8$]  if $\mathbf a_1,\mathbf a_2 \in \mathbf U_\zeta$ and 
if $A_{\mathbf a_1}/D^*_\zeta = 
A_{\mathbf a_2}/D^*_\zeta$ then $f^2_{\zeta,\mathbf a_1,\alpha} 
= f^2_{\zeta,\mathbf a_2,\alpha}$.
\end{enumerate}
\mn
Having defined $\langle f^1_{\zeta,\alpha,A}:A \in \mathbf J_{\zeta,1} \cap
\mathbf J^*_{\zeta,\alpha}\rangle$ and $\langle f^2_{\zeta,\alpha,\mathbf a}:
\mathbf a \in \mathbf U_\zeta \cap \mathbf J^*_{\zeta,\alpha}\rangle$, of course, 
they all depend on $\zeta$; we define $f_\alpha \in {}^Y\Ord$ by
\mn
\begin{enumerate}
\item[$(*)_9$]  $f_\alpha(s)$ is: the supremum below if it is
$< \delta_s$ and zero otherwise. where the supremum is
$\sup(\{f^1_{\alpha,\zeta,A}(s)+1:A \in \mathbf J_{\zeta,1} \cap 
\mathbf J^*_{\zeta,\alpha}\} \cup 
\{f^2_{\zeta,\alpha,\mathbf a}(s)+1:\mathbf a \in \mathbf U_\zeta\})$.
\end{enumerate}
\mn
So indeed $f_\alpha \in \Pi \bar\delta$ as in the end of the proof of \ref{c13}
and is as required for $\alpha$ as
$\hrtg(\mathbf J_{\zeta,1} \cap \mathbf J^*_{\zeta,\alpha}) \le
\hrtg(\cP(Y)/D_*) \le \theta$ and $\hrtg(\{f_{\zeta,\alpha,\mathbf a}:\mathbf a
\in \mathbf U_\zeta\}) \le \hrtg(\{A_{\mathbf a}:\mathbf a \in \mathbf
U_\zeta\}) \le \hrtg(\cP(Y)) \le \theta$ because of $(*)_8$ (so even
$\hrtg(\cP(Y)/D^*_\zeta)$ suffice); note that we have used $(A)(c)(\beta)$.
\medskip

\noindent
\underline{Sub-case C2}:  $(\mathbf J_{\zeta,1} \cap \mathbf J_{\zeta,2})
\cap \mathbf J^*_{\zeta,\alpha}$ is not dense in $((D^*_\zeta)^+,\supseteq)$.

Let $A_* \in (D^*_\zeta)^+$ be such that $A \subseteq A_* \wedge A \in
(D^*_\zeta)^+ \Rightarrow A \notin (\mathbf J_{\zeta,1} \cup \mathbf
J_{\zeta,2}) \cap \mathbf J^*_{\zeta,\alpha}$.  By $(*)_2$ \wilog \, 
for some $\ell \in \{1,2\}$ we have $A_* \in \mathbf J_{\zeta,\ell}$.

As in the proof of \ref{c13}, necessarily $\alpha$ is a limit ordinal
of cofinality $\ge \theta$.  Now as in Sub-Case C1 we define $\bar
f^1_{\zeta,\alpha} = \langle f^1_{\zeta,\alpha,A}:A \in (\mathbf
J_{\zeta,1} \cup \mathbf J_{\zeta,2}) \cap \mathbf
J^*_{\zeta,\alpha}\rangle$ satisfying: $f^1_{\zeta,\alpha,A}$ is a
$<_{D^*_\zeta + A}$-upper bound of $\langle f_\beta:\beta <
\alpha\rangle$.  Let $f^* \in \Pi \bar\delta$ be defined by
\mn
\begin{enumerate}
\item[$\bullet$]  $f_*(s)$ the supremum below if it is $< \delta_s$
  and is zero otherwise, where $\sup\{f^1_{\zeta,\alpha,A}(s)+1:A \in
  (\mathbf J_{\zeta,1} \cup \mathbf J_{\zeta,2}) \cap
  J^*_{\zeta,\alpha}\}$.
\end{enumerate}
\mn
As in the proof of \ref{c13} there is $\beta < \alpha$ such that
$\gamma \in [\beta,\alpha) \Rightarrow \{s \in Y:f_\gamma(s) <
f_*(s)\} = \{s \in y:f_\beta(s) < f_*(s)\} \mod D^*_\zeta$.

Let $\beta_*$ be the minimal such $\beta$.  Lastly, let $A_\zeta = \{s
\in y:f_{\beta_*}(s) \ge f_*(s)\}$ and
\mn
\begin{enumerate}  
\item[$\bullet$]  $E_\varepsilon = E_\zeta + A_\zeta$
\sn
\item[$\bullet$]  $D^*_\varepsilon = D^*_\zeta + (Y \backslash A_\zeta)$
\sn
\item[$\bullet$]  $\alpha^*_\varepsilon = \alpha^*_\varepsilon + \alpha$
\sn
\item[$\bullet$]  $g_\beta = f_\beta$ for $\beta \in
  (\alpha^*_\zeta,\alpha^*_\varepsilon)$
\sn
\item[$\bullet$]  $g_{\alpha^*_\varepsilon} = f_*$.
\end{enumerate}
\bigskip

\noindent
\underline{Case D}:  None of the above.

So $Y \in D^*_\varepsilon$ and we are done.
\end{PROOF}

\begin{discussion}
\label{c22}
In the results above, is 
$\langle \cf(\alpha^*_{\varepsilon +1}):\varepsilon <
\varepsilon(*)\rangle$ without repetitions?  Certainly this is not
obviously so and it seems we can manuever $\bar\delta$ and the closure
operation to be otherwise.  But can we replace $\bar\alpha^*$ and
$\bar g$ to take care of this?  Clearly if $\cU \subseteq
\alpha^*_{\varepsilon(*)}$ satisfies $\varepsilon < \varepsilon(*)
\Rightarrow \alpha^*_{\varepsilon +1} = \sup(\cU \cap
\alpha^*_{\varepsilon +1})$ then we can replace $\bar g$ by $\bar g
\rest \cU$ so by renaming get $\bar\alpha' = \langle\text{otp}(\cU
\cap \alpha^*_\varepsilon):\varepsilon \le \varepsilon(*)\rangle$.  So
cf$(\alpha^*_\varepsilon) = \text{ cf}(\alpha^*_\zeta) \Leftrightarrow
\text{ cf}(\alpha'_\varepsilon) = \text{ cf}(\alpha'_\zeta)$ and if we
have cf$(\alpha'_\varepsilon) = \text{ cf}(\alpha'_\zeta) \Rightarrow
\alpha'_{\varepsilon +1} - \alpha'_\varepsilon = \alpha'_{\zeta +1} -
\alpha'_\zeta$ we can change $\bar g$ to get desired implication.  So
if AC$_{\varepsilon(*)}$ holds we are done but we are not assuming
it.  In this case we also get $\langle \alpha'_{\varepsilon +1} \backslash
   \alpha'_\varepsilon:\varepsilon < \varepsilon(*)\rangle$ is a
   sequence of regular cardinals.
\end{discussion}
\newpage

\section {More on the pcf theorem} \label{2} 
\bigskip

\subsection {When the Cofinalities are Smaller} \label{2A} \
\bigskip

\begin{definition}
\label{d2}
1) We say $\mathbf x = (Y,\bar\delta,\theta,
\varepsilon(*),\bar\alpha^*,\bar D^*,\bar E^*,\bar f) = 
(Y_{\mathbf x},\bar\delta_{\mathbf x},\theta_{\mathbf x},\varepsilon_{\mathbf x},
\bar\alpha_{\mathbf x},\bar D_{\mathbf x},\bar E_{\mathbf x},\bar f_{\mathbf x})$ is 
a pcf-system or a pcf-system for $\bar\delta$ or for 
$(\Pi \bar\delta,<_D)$ \when \, they are
as in $(B)^+$ of \ref{c17}, with $\bar f$ here standing for $\bar g$ there; so 
$\bar\delta = \langle \delta_s:s \in Y\rangle,\delta_s$ a
limit ordinal; now \ref{c24} below apply,
we will use $\bar D_{\mathbf x} =
   \langle D^{\mathbf x}_\varepsilon:\varepsilon < \varepsilon_{\mathbf
   x}\rangle = \langle D_{\mathbf x,\varepsilon}:\varepsilon <
   \varepsilon_{\mathbf x}\rangle$, similarly for $\bar f,D_{\mathbf x} =
   D^{\mathbf x}_0$; let $\varepsilon(\mathbf x) = \varepsilon_{\mathbf x}$.

\noindent
2) Above we say is ``almost a $\pcf$-system" if we demand $\bar f \rest
   [\alpha_{\mathbf x,\varepsilon},\alpha_{\mathbf x,\varepsilon +1})$ is
   only $\le_{D_{\mathbf x,\varepsilon}}$-increasing (still cofinal) so
using $(B)^+(i)'$ instead of $(B)^+(i)$, see \ref{c17},\ref{c19}(7).

\noindent
3) Above we say $\mathbf x$ is ``weakly a $\pcf$-system" \when \, in
\ref{c17}(B)$^+$ - we weaken clause (i) as in part (2) and 
we omit $\bar E^*$, i.e. omit clauses (j),(k) but retain (l) which
means: if $X_0 \in D^*_{\varepsilon +1} \backslash D^*_\varepsilon,
X_1 = Y \backslash X_0$ then $\bar f \rest 
[\alpha^*_\varepsilon,\alpha^*_\varepsilon)$ is
$\le_{D^*_\varepsilon}$-increasing and cofinal in 
$(\Pi \bar\delta,<_{D^*_\varepsilon +
   X_1})$ and $\bar f$ is $\le_{D^*_\varepsilon}$-increasing.
\end{definition}

\begin{observation}
\label{d5}
1) If $\theta,Y,D$ and $\bar\delta = \langle \delta_s:s \in Y\rangle$
satisfies clause (A) of \ref{c17}, \then\, there is a pcf-sytem $\mathbf
x$ for $(\Pi \bar\delta,<_D)$ with $\theta_{\mathbf x} = \theta$.

\noindent
2) We can above use $D = \cf-\fil_{< \theta}(\bar\delta)$.
\end{observation}

\begin{PROOF}{\ref{d5}}
By \ref{c17}.
\end{PROOF}

\begin{observation}
\label{c24}
Let $\mathbf x = (Y,\bar\delta,\theta,\varepsilon(*),
\bar\alpha^*,\bar D^*,\bar E^*,\bar f)$
 be as in \ref{c17} (with $\bar f$ instead of $\bar g$)
or Definition \ref{d2}(2).

\noindent
1) $(\Pi \bar\delta,<_{D_{\mathbf x}})$ has a cofinal well 
orderable subset, in fact, of
cardinality $|\alpha^*_{\varepsilon(*)}|$.

\noindent
2) Assume $f \in \Pi \bar\delta$ and for $\varepsilon <
\varepsilon(*)$ we let $\beta_\varepsilon = \min\{\beta:\beta \in
[\alpha^*_\varepsilon,\alpha^*_{\varepsilon +1})$ satisfy $f <
f_\beta \mod(E^*_{\varepsilon +1})\}$, then:
\mn
\begin{enumerate}
\item[$(a)$]  $\beta_\varepsilon \in
[\alpha^*_\varepsilon,\alpha^*_{\varepsilon +1})$ is well defined
hence $\langle \beta_\varepsilon:\varepsilon < 
\varepsilon(*)\rangle$ is well defined
\sn
\item[$(b)$]  for some finite $u \subseteq \varepsilon(*)$ we have $f <
\sup\{f_{\beta_\varepsilon}:\varepsilon \in u\}$
\sn
\item[$(b)^+$]  moreover $\langle f_{\beta_\varepsilon}:\varepsilon
  \in u\rangle$ is $\partial$-uniformly definable from $f$ and
  $\bar\delta$ and $D^*_0$ (equivalently, $f$ and $\mathbf x$).
\end{enumerate}
\end{observation}

\begin{PROOF}{\ref{c24}}
1) By (2).

\noindent
2) Easy; e.g..

\noindent
\medskip

\noindent
\underline{Clause (b)}:

Let $\varepsilon \le \varepsilon_{\mathbf x}$ be minimal such that
\mn
\begin{enumerate}
\item[$(*)$]  $\varepsilon = \varepsilon_*$ for some finite $u
  \subseteq [\varepsilon,\varepsilon_{\mathbf x})$ we have $f <
  \max\{f_{\beta_\zeta}:\zeta \in u\} \mod D_{\mathbf x,\varepsilon}$.
\end{enumerate}
\mn
Now $\varepsilon$ is well defined because $\varepsilon_{\mathbf x}$ is a
successor ordinal and $\langle f_\beta:\beta <
\alpha^*_{\varepsilon(\mathbf x)}\rangle$ is cofinal in $(\Pi
\bar\delta,<_{D_{\mathbf x,\varepsilon(\mathbf x)-1}})$ and so $u =
\{\beta_{\varepsilon(\mathbf x)-1}\}$ is as required.

If $\varepsilon = \zeta +1 < \varepsilon_{\mathbf x}$ and $u$ are as in
$(*)$ the set $Z = \{s \in Y:f(s) < \max\{f_{\beta_\zeta}(s):\zeta \in
u\}$ is $= \emptyset \mod E_{\zeta +1}$ and repeat the argument for
$\varepsilon = \varepsilon_{\mathbf x} -1$.

If $\varepsilon$ is a limit ordinal, this leads to contradiction as
$D_{\mathbf x,\varepsilon} = \cup\{D_{\mathbf x,\zeta}:\zeta <
\varepsilon\}$.

Lastly, if $\varepsilon = 0$ then we are done.
\end{PROOF}

\begin{discussion}
\label{c29}
1) In \ref{c24}, we may restrict ourselves to $\aleph_1$-complete
filters only, so
replace $\varepsilon_*$ by $\{\varepsilon <
\varepsilon_*:E^*_\varepsilon$ is $\aleph_1$-complete$\}$ but use
countable $u$. 

\noindent
2) Similarly for $\theta$-complete.

\noindent
3) Recall that with choice or just AC$_Y$, the ideal $\cf-\id_{<
\theta}(\bar\delta)$ is degenerate: if, for transparency, 
$\theta$ is regular, then
cf-id$_{< \theta}(\bar \delta) = \{X \subseteq Y:(\forall s \in
X)[\cf(\delta_s) < \theta]$ and $|X| < \theta\}$.

We have dealt with $(\prod\limits_{s} \delta_s,<_D)$ when $D
\supseteq \cf-\fil_{< \theta}(\bar\delta)$ and $\theta \ge
\hrtg(\Fil^4_{\aleph_1}(Y))$; we try to lower the restriction on the
cardinal $\theta$ with some price.
\end{discussion}

\begin{definition}
\label{c41}
Assume $D$ is a filter on $Y,\alpha(*)$ an ordinal and $\bar f =
\langle f_\alpha:\alpha < \alpha(*)\rangle$ is a $\le_D$-increasing
sequence of members of ${}^Y\Ord$ and $f \in {}^Y\Ord$ is not 
$<_D$-below any $f_\alpha$.  We define

\begin{equation*}
\begin{array}{clcr}
\id(f,\bar f,D) = \{Z \subseteq Y: &\text{ there is } \alpha < \alpha(*)
\text{ such that} \\
  &Z \subseteq \{s \in Y:f(s) < f_\alpha(s)\} \mod D\}.
\end{array}
\end{equation*}
\end{definition}

\begin{claim}
\label{c44}
For $Y,D,\bar f,f$ as in Definition \ref{c41} above.

\noindent
1) $\id(f,\bar f,D)$ is an ideal on $Y$ extending $\dual(D)$.

\noindent
2) $f$ is a $\le_{\id(f,\bar f,D)}$-upper bound of $\bar f$.

\noindent
3) For $A \in D^+$ we have: $\cP(A) \cap \id(f,\bar f,D) \subseteq
\dual(D)$ \Iff \, $f$ is a $\le_{D+A}$-upper bound of $\bar f$.

\noindent
4) If $A \in D^+ \cap \id(f,\bar f,D)$ \then \, for every $\alpha <
   \alpha(*)$ large enough, $f < f_\alpha \mod(D+A)$.

\noindent
5) $\id(f,\bar f,D) = \id(f',\bar f,D)$ when $f' \in {}^Y\Ord$ and $f'
   =_D f$.
\end{claim}

\begin{PROOF}{\ref{c44}}
Straightforward.
\end{PROOF}

\begin{notation}
\label{c45}
1) Given $\bar\delta = \langle \delta_s:s \in Y\rangle$ and set $u$ of
ordinals let $h_{[u,\bar\delta]}$ be the function $h$ with domain $Y$
such that: $h(s)$ is $\sup(u \cap \delta_s)$ \when \, it is $<
\delta_s$, is $0$ \when \, otherwise.

\noindent
2) For $\bar u = \langle u_s:s \in Y\rangle$ we define $h_{[\bar
   u,\bar\delta]}$ similarly.
\end{notation}

\begin{claim}
\label{c47}
If we assume $\oplus$ below and (A) + (B) \then \, (C) where:
\mn
\begin{enumerate}
\item[$\oplus$]  $(a) \quad \Ax_{4,\theta} \wedge |Y| \le
\aleph_0$ 
\sn
\item[${{}}$]  $(b)_{\kappa,\theta} \quad$ the union of any sequence
  of length $\le \kappa$ of sets of ordinals

\hskip25pt  each of cardinality $< \theta$ is of cardinality $< \theta$
\sn
\item[${{}}$]  $(c) \quad \kappa \le \theta$
\sn
\item[$(A)$]  $(a) \quad \bar \delta = \langle \delta_s:s \in
Y\rangle$ is a sequence of limit ordinals
\sn
\item[${{}}$]  $(b) \quad D$ is a filter on $Y$
\sn
\item[${{}}$]  $(c) \quad D \supseteq \cf-\fil_{<
\theta}(\bar\delta)$
\sn
\item[${{}}$]  $(d) \quad \mu = \cup\{\delta_s:s \in Y\}$
\sn
\item[$(B)$]  $\delta_*$ is an ordinal and
\sn
\item[${{}}$]  $(a) \quad f_\alpha \in \prod\limits_{s \in Y} \delta_s$
for $\alpha < \delta_*$
\sn
\item[${{}}$]  $(b) \quad$ if $\alpha < \beta < \delta_*$ then
$f_\alpha < f_\beta \mod D$
\sn
\item[${{}}$]  $(c) \quad \bar f = \langle f_\alpha:\alpha <
  \delta_*\rangle$ is not cofinal in $(\prod\limits_{s
\in Y} \delta_s,<_D)$
\sn
\item[${{}}$]  $(d) \quad \cf(\delta_*) > \kappa$
\sn
\item[$(C)$]  we can $\theta$-uniformly define (or
 $(\theta,\kappa)$-uniformly define) $g$ such that:
\sn
\item[${{}}$]  $(a) \quad g \in \prod\limits_{s \in Y} \delta_s$ is
not $<_D$-below any $f_\alpha$
\sn
\item[${{}}$]  $(b) \quad$ if $g \le_D g' \in \prod\limits_{s \in Y}
\delta_s$ then $\id(g',\bar f,D) = \id(g,\bar f,D)$.
\end{enumerate}
\end{claim}

\begin{remark}
1) See more in \ref{c59}.

\noindent
2) Do we uniformly have the parallel of: some stationary $S \subseteq
   S^\lambda_{\kappa^+}$ belongs to $\check I_{\kappa^+}[\lambda]$?
   See later.

\noindent
3) We can weaken \ref{c47} $\oplus(a)$ to 
$\Ax_{4,\mu,\theta,\kappa} \wedge \hrtg(Y) \le \kappa$, (see
\ref{z4}(3)) the proof is written for this.
\end{remark}

\begin{PROOF}{\ref{c47}}

\noindent
\underline{Stage A}:  

Let $(\cS_*,<_*)$ witness $\Ax_{4,\mu,\theta,\kappa}$.

We try to choose $g_\varepsilon,u_\varepsilon,Y_\varepsilon$ by
induction on $\varepsilon < \kappa$ such that:
\mn
\begin{enumerate}
\item[$\boxplus$]  $(a) \quad g_\varepsilon \in \prod\limits_{s \in Y}
\delta_s$
\sn
\item[${{}}$]  $(b) \quad u_\varepsilon \subseteq \mu$ has cardinality
$< \theta$ and $\zeta < \varepsilon \Rightarrow u_\zeta \subseteq
u_\varepsilon$
\sn
\item[${{}}$]  $(c) \quad Y_\varepsilon = \{s \in Y:\delta_s =
\sup(\delta_\varepsilon \cap u_\varepsilon)\} = \emptyset \mod D$ 
\sn
\item[${{}}$]  $(d) \quad$ if $s \in Y \backslash Y_\varepsilon$ and
$\zeta < \varepsilon$ then $g_\zeta(s) < g_\varepsilon(s)$ 
\sn
\item[${{}}$]  $(e) \quad g_\varepsilon = h_{[u_\varepsilon,\bar\delta]}$, see
Definition \ref{c45}
\sn
\item[${{}}$]  $(f) \quad$ if $\varepsilon$ is a limit ordinal \then \,:
\sn
\begin{enumerate}
\item[${{}}$]  $\bullet \quad u_\varepsilon = 
\cup\{u_\zeta:\zeta < \varepsilon\}$ 
\sn
\item[${{}}$]  $\bullet \quad g_\varepsilon(s)$ is
$\cup\{g_\zeta(s):\zeta < \varepsilon\}$ \when \, it is $< \delta_s$

\hskip35pt is $0$ \when \, otherwise
\end{enumerate}
\sn
\item[${{}}$]  $(g) \quad$ if $\varepsilon = \zeta +1$ then
\sn
\begin{enumerate}
\item[${{}}$]  $(\alpha) \quad g_\zeta$ is not as required on $g$ in
clause (C)
\sn
\item[${{}}$]  $(\beta) \quad u_\varepsilon$ is the $<_*$-first $u \in
\cS_*$ extending $u_\zeta$ such that if we define $g_\varepsilon$ 

\hskip35pt  as $h_{[u,\bar\delta]}$ \then \, it 
is a counterexample like $g'$ there
\end{enumerate}
\sn
\item[${{}}$]  $(h) \quad$ if $\varepsilon = 0,g_\varepsilon$ is
defined from $u_\varepsilon$ similarly.
\end{enumerate}
\mn
Now we shall finish by proving in stages B,C below that:
\mn
\begin{enumerate}
\item[$(*)_1$]  if we have defined $g_\varepsilon$ but $g_\varepsilon$
is as required on $g$ in clause (C)(b), then we are done; this is
obvious
\sn
\item[$(*)_2$]   we can choose $g_\varepsilon$ if $\varepsilon =0$ 
\sn
\item[$(*)_3$]  if $\langle g_\zeta:\zeta < \varepsilon \rangle$ was
defined we can define $g_\varepsilon$ if $\varepsilon$ is a limit
ordinal $< \kappa$
\sn
\item[$(*)_4$]  if $\varepsilon = \zeta +1$ and $\langle g_\xi:\xi \le
\zeta\rangle$ has been defined and $g_\zeta$ fail (C), then we can
define $g_\varepsilon$
\sn
\item[$(*)_5$]  we cannot succeed to choose $\langle
g_\varepsilon:\varepsilon < \kappa \rangle$.
\end{enumerate}
\bigskip

\noindent
\underline{Stage B}:

\noindent
\underline{Proof of $(*)_5$}:

Toward contradiction assume $\langle g_\varepsilon:\varepsilon <
\kappa\rangle$ is well defined.

For $\varepsilon < \kappa$ and $\alpha < \delta_*$ let
$Z_{\varepsilon,\alpha} = \{s \in Y:g_\varepsilon(s) \ge
f_\alpha(s)\}$ and let $Y_\varepsilon = \{s \in Y:\sup(u_\varepsilon
\cap \delta_s) = \delta_s\}$, it belongs.
By clauses (b),(c),(e) of $\boxplus$ we have $Z_{\varepsilon_1,\alpha}
 \backslash Y_{\varepsilon_1} \subseteq Z_{\varepsilon_2,\alpha}
 \backslash Y_{\varepsilon_2}$ for $\varepsilon_1 <
 \varepsilon_2 < \kappa,\alpha < \delta_*$.

Now by clause $(g)(\beta)$ of $\boxplus$, if
$\varepsilon = \zeta +1$ then for some $\alpha < \delta_*,
Z_{\varepsilon,\alpha} \notin \id(g_\zeta,\bar f,D)$ and
let $\alpha_\zeta$ be the minimal such $\alpha$.  As $\cf(\delta_*) >
\kappa$ by Clause (B)(d) of the assumption, $\gamma :=
\cup\{\alpha_\zeta:\zeta < \kappa\}$ is $< \delta_*$.

Now the sequence $\langle Y_\varepsilon:\varepsilon < \kappa
\rangle$ is $\subseteq$-increasing sequence of subsets of $Y$ because
$\langle u_\varepsilon:\varepsilon < \kappa\rangle$ is by
$\boxplus(b)$ and the choice of $Y_\varepsilon$.  By 
$\oplus(a)$ we have $\hrtg(Y) \le \kappa$.

Also clearly
\mn
\begin{enumerate}
\item[$\bullet_2$]  $Z_{\varepsilon +1} \nsubseteq Z_{\varepsilon +1}
  \mod D$ and $Y_\varepsilon$.
\end{enumerate}
\mn
Together $\langle Z_{\varepsilon +1,\gamma} \backslash 
Z_{\varepsilon,\gamma} \backslash
Y_\varepsilon:\varepsilon < \kappa \rangle$ is a sequence pairwise
distinct non-empty of subsets of $Y$, so recalling
$\hrtg(Y) \le \kappa$, this is contradiction to the first paragraph.
\bigskip

\noindent
\underline{Stage C}:

Obviously $(*)_1$ holds.

\noindent
\underline{Proof of $(*)_2$}:  we can choose $g_\varepsilon$ for
$\varepsilon =0$
\mn
\begin{enumerate}
\item[$\bullet_1$]  there is $g'' \in \prod\limits_{s \in Y} \delta_s$
  such that $\alpha < \delta_* \Rightarrow g'' \nleq f_\alpha \mod D$.
\end{enumerate}
\mn
[Why?  By clause (B)(c) of the claim.  For such a $g''$ there is 
$u \in \cS_*$ such that $\Rang(g'')
\subseteq u$ because $\hrtg(Y) \le \kappa$ and $\cS_*$ witness
$\Ax_{4,\mu,\theta,\kappa}$.  
We choose $u \in \cS_*$ as the $<_*$-first such $u \in \cS_*$ and
choose $g \in \prod\limits_{s \in Y} \delta_s$ as $h_{[u,\delta]}$.] 

So
\mn
\begin{enumerate}
\item[$\bullet_2$]   $g \in \prod\limits_{s \in Y} \delta_s$
\sn
\item[$\bullet_3$]  $g'' \le g \mod D$.
\end{enumerate}
\mn
[Why?  Recall $\cf-\fil_{< \theta}(\bar\delta) \subseteq D$ by the
assumption (A)(c), hence $\{s \in Y:\sup(u \cap \delta_2) \le g(s)\}$
as $|u| < \theta$ being a membre of $\cS_*$.  So as $(\forall s \in
Y)(g''(s) \in \delta_s \cap u)$ we have $g'' \le g \mod D$ by the
choice of $u$.]
\mn
\begin{enumerate}
\item[$\bullet_4$]   $\alpha < \delta_* \Rightarrow g \nleq f_\alpha
 \mod D$.
\end{enumerate}
\mn
[Why?  By $\bullet_3$ and by the choice of $g''$ in $\bullet_1$.]
\medskip

\noindent
\underline{Proof of $(*)_3$}:  limit $\varepsilon$

We define $g_\varepsilon$ as in $\boxplus(f)$, as it is as required
because $D \supseteq \cf-\fil_{< \theta}(\bar\delta)$ by clause (A)(c)
of the assumption recalling $\oplus(b)_{\kappa,\theta}$ of the assumption.
\medskip

\noindent
\underline{Proof of $(*)_4$}:  

So we are assuming $g_\zeta$ is well defined but fail (C)(b) as
exemplified by $g$, let $u \in \cS_*$ be $<_*$-minimal such that
$\Rang(g) \subseteq u$ and let $h = h^*_{[u,\bar\delta]} +1$, that is
$s \in Y \Rightarrow h(s) = h_{[u,\bar\delta]}(s)+1 < \delta_s$
 hence $g <_J h_{[u,\bar\delta]} \mod D$ and we can finish easily as
 in the proof of $(*)_2$.
\end{PROOF}

\begin{observation}
\label{c51}
$\cf(\alpha(*)) \ge \theta$ \when \,
\mn
\begin{enumerate}
\item[$(a)$]  $D$ is a filter on $Y$
\sn
\item[$(b)$]  $\bar\delta = \langle \delta_s:s \in Y\rangle$ is a
  sequence of limit ordinals
\sn
\item[$(c)$]  $D \supseteq \cf-\fil_{< \theta}(\bar\delta)$
\sn
\item[$(d)$]  $\bar f = \langle f_\alpha:\alpha <
\alpha(*)\rangle$ is $<_D$-increasing sequence of members of
$\prod\limits_{s \in Y} \delta_s$
\sn
\item[$(e)$]  $\bar f$ has no $<_D$-upper bound in
$\prod\limits_{s \in Y} \delta_s$.
\end{enumerate}
\end{observation}

\begin{PROOF}{\ref{c51}}
The proof splits into cases proving the existence of a $<_D$-upper bound
$g \in \prod\limits_{s \in Y} \delta_s$.
\bigskip

\noindent
\underline{Case 1}:  $\alpha(*) = 0$

The constantly zero function $g:Y \rightarrow \{0\}$ can serve.
\bigskip

\noindent
\underline{Case 2}:  $\alpha(*)$ is a successor ordinal

Let $\alpha(*) = \beta +1$ and $g$ be defined by $g(s) = f_\beta(s)
+1$.  As each $\delta_s$ is a limit ordinal, $g \in \prod\limits_{s \in Y}
\delta_s$.
\bigskip

\noindent
\underline{Case 3}:  $\cf(\alpha(*)) \in [\aleph_0,\theta)$

Let $w \subseteq \alpha(*)$ be cofinal of order type $\cf(\alpha(*))$,
let $u_s = \{f_\alpha(s):\alpha \in w\}$ for $s \in Y$ so $\bar u := \langle
u_s:s \in Y\rangle$ is well defined and $s \in Y \Rightarrow |u_s| < \theta$, 
hence $g=h_{[\bar u,\bar\delta]}$ is as required.
\end{PROOF}

\begin{claim}
\label{c54}
If $\boxplus$ below holds then $\oplus_1 \Rightarrow
\oplus_2 \Rightarrow \oplus_3$ \underline{where}
\mn
\begin{enumerate}
\item[$\oplus_1$]  $\Ax_{4,\mu,\theta,\kappa}$
\sn
\item[$\oplus_2$]  there is a well orderable set cofinal in 
$(\Pi \bar \delta,<_D)$, defined $(\mu,\theta,\kappa)$-uniformly
\sn
\item[$\oplus_3$]  we can $(\theta,\kappa)$-uniformly define a 
$<_D$-increasing sequence $\bar f = \langle f_\alpha:\alpha < 
\alpha(*)\rangle$ in $(\prod\limits_{s \in Y} \delta_s,<_D)$ with 
no upper bound
\end{enumerate}
\mn
where
\mn
\begin{enumerate}
\item[$\boxplus$]  $(a) \quad D$ a filter on $Y$
\sn
\item[${{}}$]  $(b) \quad \bar\delta = \langle \delta_s:s \in Y \rangle$ is a
sequence of limit ordinals
\sn
\item[${{}}$]  $(c) \quad D \supseteq \cf-\fil_{< \theta}(\bar\delta)$
\sn
\item[${{}}$]  $(d) \quad \hrtg(Y) \le \kappa \le \theta$
\sn
\item[${{}}$]  $(e) \quad \mu =\sup\{\delta_s:s \in Y\}$.
\end{enumerate}
\end{claim}

\begin{PROOF}{\ref{c54}}
\smallskip

\noindent
\underline{$\oplus_1 \Rightarrow \oplus_2$}

Let $(\cS_*,<_*)$ witness $\Ax_{4,\mu,\theta,\kappa}$.

For every $g \in \Pi \bar\delta,\Rang(g)$ is a subset of
$\sup\{\delta_s:s \in Y\} = \mu$ of cardinality $< \hrtg(Y) \le
\kappa$ hence there is $u \in \cS_*$ such that
$\Rang(g) \subseteq u$, so $|u| < \theta$ 
hence easily $g \le h_{[u,\bar\delta]} \mod D$, see Definition \ref{c45}.  
Hence $\cF = \{h_{[u,\bar\delta]}:u \in \cS_*\}$ is a cofinal
subset of $(\Pi \bar\delta,<_D)$ and being $\le_{\qu} \cS_*$ it is well
orderable.  Recall $h_{[u,\bar\delta]} \in \Pi \bar\delta$ is defined
by:
 $h_{[u,\bar\delta]}(s)$ is $\sup(\delta_s \cap u)$ if $\sup(\delta_2
\cap u) < \delta_s$ and is zero otherwise.

Now $\cF \subseteq \Pi \bar\delta$ being cofinal in $(\Pi
\bar\delta,<_D)$ follows from $D \supseteq \cf-\fil_{<
  \theta}(\bar\delta)$ that is $\boxplus(c)$.
\smallskip

\noindent
\underline{$\oplus_2 \Rightarrow \oplus_3$}

Let $\cF \subseteq \Pi \delta$ be cofinal in $(\Pi \bar\delta,<_D)$ and
$<_*$ well order $\cF$.
We try to choose $f_\alpha$ by induction on the ordinal 
$\alpha$.  If $\bar f^\alpha =
\langle f_\beta:\beta < \alpha\rangle$ has no $<_D$-upper bound we are
done so assume $g \in \prod\limits_{s \in Y} \delta_s$ is a $<_D$-upper
bound of $\bar f^\alpha$ so there is $h \in \cF$ such that $g <_D h$, 
so $h$ is a $<_D$-lub of $\bar f$ and let $f_\alpha \in \cF$ be
the $<_*$-minimal such $h$.  Necessarily for some $\alpha$ we cannot
continue so $\bar f^\alpha$ is as promised.
\end{PROOF}

\begin{conclusion}
\label{c56}
In clause (C) of \ref{c47} letting 
\mn
\begin{enumerate}
\item[$\bullet$]  $Z_\alpha = \{s \in Y:g(s) < f_\alpha(s)\} 
\text{ for } \alpha < \delta_*$
\sn
\item[$\bullet$]  $\cW = \{\alpha < \delta_*:Z_\beta \ne Z_\alpha 
\mod D \text{ for every } \beta < \alpha\}$
\sn
\item[$\bullet$]  $D_\alpha = D+ Z_\alpha$ for $\alpha < \delta_*$
\sn
\item[$\bullet$]  $\alpha_* = \min\{\alpha \le \delta_*$: if $\alpha 
< \delta_*$ then $Z_\alpha \in D^+\}$
\end{enumerate}
\mn
and assuming (B)(e) $\bar f$ has no $\le_D$-ub in $\Pi \bar\delta$ we can add:
\mn
\begin{enumerate}
\item[$(c)$]  $\langle Z_\alpha/D:\alpha \in \cW\rangle$ is 
$\subseteq$-increasing and $\alpha_* < \delta_*$
\sn
\item[$(d)$]  for $\alpha \in \cW,\alpha \ge \alpha_*,D_\alpha$ is a filter
on $Y$ and $\langle f_{\alpha + \gamma}:\gamma < \delta_* - \alpha$ and
$\alpha + \gamma \in \cW\rangle$ is
$<_{D_\alpha}$-increasing and cofinal in $\Pi \bar\delta$
\sn
\item[$(e)$]  $\langle D_\alpha:\alpha \in \cW \backslash \alpha_*\rangle$ is a
strictly $\subseteq$-increasing sequence of filters of $Y$ and 
$0 \in \cW$ 
\sn
\item[$(f)$]  $\bar f$ is $<_{D_\alpha}$-increasing and
  $<_{D_\alpha}$-cofinal in $\Pi \bar\delta$ if $\alpha \in \cW \backslash 
\alpha_*$
\sn
\item[$(g)$]  if $\cf(\delta_*) \ge \hrtg(\cP(Y))$ \then \, $\cW$ has
  a last member.
\end{enumerate}
\end{conclusion}

\begin{PROOF}{\ref{c56}}
Easy or see \cite[Ch.II,\S2]{Sh:g}; but we elaborate.
\medskip

\noindent
\underline{Clause (c)}:  First, the sequence is $\subseteq$-increasing as
$\bar f$ is $<_D$-increasing.  Second, $\alpha_* < \delta_*$ as
otherwise we have $\alpha < \delta \Rightarrow f_\alpha \le g \mod D$
but we are assuming $\bar f$ has no $\le_D$-ub in $\Pi \bar\delta$.
\medskip

\noindent
\underline{Clause (d)}:  $D_\alpha$ is a filter as by clause (c),
$\alpha \ge \alpha_* \Rightarrow Z_\alpha \in D^+$ and obviously
$Z_\alpha \in D^+ \Rightarrow$ ($D_\alpha$ is a filter). 
\medskip

\noindent
\underline{Clause (e)}:  By the definition of $\cW$.
\medskip

\noindent
\underline{Clause (f)}:  By (C)(a),(b) and clause (d).
\medskip

\noindent
\underline{Clause (g)}:  Obvious.
\end{PROOF}

\begin{theorem}
\label{c59}
Assume $\boxplus(a)-(e)$ of \ref{c54}.

\noindent
1) If $\cf(\theta) \ge \hrtg(\cP(Y))$ and 
$\Ax_{4,\mu,\theta,\kappa}$, \then \,
the conclusion $(B)^+$ of Theorem \ref{c17} holds, i.e. there is a
$\pcf$-system $\mathbf x$ such that $Y_{\mathbf x} = Y,\bar\delta_{\mathbf
  x} = \bar\delta,\theta_{\mathbf x} = \theta$. 

\noindent
2) Without the extra assumption $\cf(\theta) \ge \hrtg(\cP(Y))$, we
get only a weakly $pcf$-system (see \ref{d2}(3)) $\mathbf x$ with
$\theta = \hrtg(\cP(Y))$.

\noindent
3) If there is a weak $\pcf$-system $\mathbf x$ for $\bar\delta$ 
\then \, $\Pi \bar\delta$ has a subset which is a well-orderable and is cofinal
   in $(\Pi \bar\delta,<_{D_{\mathbf x}})$.

\noindent
4) If $(\Pi \bar\delta,<_D)$ has a well-orderable cofinal subset and
$\hrtg(\cP(Y)) \le \theta$ \then \, there is a $\pcf$-system $\mathbf x$
for $\bar\delta$ with $D_{\mathbf x} = D$.

\noindent
5) If $(\Pi \bar\delta,<_D)$ has a well-ordered cofinal subset and
$\theta \ge \hrtg(Y)$ \then \, there is a $\pcf$-system $\mathbf x$ for
$\bar\delta$ with $D_{\mathbf x} = D,\alpha_{\mathbf x,\varepsilon +1} -
\alpha_{\mathbf x,\varepsilon}$ increasing.
\end{theorem}

\begin{remark}
Note that later parts of \ref{c59} supercede earlier ones.  One
reason for this is that it may be better to avoid using inner models,
developing the set theory of $ZF + DC + Ax_4$ per se.
\end{remark}

\begin{PROOF}{\ref{c59}}
1) We repeat the proof of \ref{c17}, but using \ref{c47}, \ref{c51},
\ref{c56}, i.e. in case (c) after $(*)_3$ we use \cite{Sh:E62}.      
But a simpler argument is that by
\ref{c54} we know that there is a $<_D$-cofinal subset $\cF$ of
$\Pi\bar\delta$ which is well orderable, say by $<_*$.

\noindent
2) Like part (1).

\noindent
3) Let $\mathbf x$ be a weak pcf-system for $(\Pi \bar\delta,<_D)$,
   clearly $\{f_{\mathbf x,\alpha}:\alpha < \alpha_{\mathbf
   x,\varepsilon(\mathbf x)}\}$ is a well orderable subset of 
$\Pi \bar\delta$ and so is $\cF = \{\max\{f_{\mathbf x,\alpha_\ell}:\ell < n\}:
\bar\alpha = \langle \alpha_\ell:\ell < n\rangle$ is
   a finite sequence of ordinals $\langle 
\alpha_{\mathbf x,\varepsilon(\mathbf x))}\rangle\}$.  
Hence it suffices to prove that the set $\cF$ 
is cofinal in $(\Pi\bar\delta,<_{D_{\mathbf x}})$. 

This means to show that
\mn
\begin{enumerate}
\item[$(*)$]  for every $g \in \Pi \bar\delta$ there are $n$ and
  $\alpha_\ell < \alpha_{\mathbf x,\varepsilon(\mathbf x)}$ for $\ell < n$
  such that $g < \max\{f_{\mathbf x,\alpha_\ell}(s):\ell < n\} \mod D_{\mathbf x}$.
\end{enumerate}
\mn 
For this we prove by induction on $\varepsilon \le \varepsilon_{\mathbf x}$ that
\mn
\begin{enumerate}
\item[$(*)_\varepsilon$]  if $X \in D_{\mathbf x,\varepsilon}$ and 
$g \in \Pi \bar\delta$ \then \, we can find $Z \in D_{\mathbf x}$ and
  $n$ and $\alpha_\ell < \alpha_{\mathbf x,\varepsilon}$ for $\ell < n$ 
such that $s \in Z 
\backslash X \Rightarrow g(s) < \max\{f_{\mathbf x,\alpha_\ell}(s): \ell < n\}$.
\end{enumerate}
\mn 
This suffices as for $\varepsilon = \varepsilon_{\mathbf x}$ we can use
$X = \emptyset$.

For $\varepsilon = 0$ necessarily $Z := X$ is as required because 
$X \in D_{\mathbf x,\varepsilon} = D_{\mathbf x}$.

For $\varepsilon$ a limit ordinal, if $X \in D_{\mathbf x,\varepsilon}$
then for some $\zeta < \varepsilon,X \in D_{\mathbf x,\zeta}$ and use
the induction hypothesis for $\zeta$.

For $\varepsilon = \zeta +1$, we are given $X \in D_{\mathbf x,\varepsilon}$
and $g \in \Pi \bar\delta$.  By clause $(B)^+(\ell)$ of \ref{c17} if 
$X \in D_{\mathbf x,\zeta}$ use the
induction hypothesis so \wilog \, $X \notin D_{\mathbf x,\zeta}$ hence
$D_{\mathbf x,\zeta} + (Y \backslash X)$ is a filter on $Y_{\mathbf x}$
and it is $\supseteq E_{\mathbf x,\zeta}$.  So by clause $(B)^+(l)$ of
Theorem \ref{c17} there is $\alpha \in [\alpha_{\mathbf x,\zeta},
\alpha_{\mathbf x,\zeta +1})$ such that $g < f_{\mathbf x,\alpha} 
\mod (D_{\mathbf x,\zeta} + (Y_{\mathbf x} \backslash X))$.

Let $X_1 = \{s \in Y:s \notin X$ and $g(s) < f_{\mathbf x,\alpha}(s)\}$,
so $X_1 \in D_{\mathbf x,\zeta} + (Y_{\mathbf x} \backslash X)$ hence $X_2
:= X \cup X_1 \in D_{\mathbf x,\zeta}$ so by the induction hypothesis
there are $n_1$ and $\beta_\ell < \alpha_{\mathbf x,\zeta}$ for $\ell <
n_1$ and $Z \in D_{\mathbf x}$ such that $s \in Z \backslash X_2
\Rightarrow g(s) < \max\{f_{\mathbf x,\beta_\ell}(s):\ell < n_1\}$.  Let
$n=n_1 +1$ and let $\alpha_\ell$ be $\beta_\ell$ if $\ell < n_1,\alpha_\ell$ be
$\alpha$ if $\ell=n_1$, so $Z,\langle \alpha_\ell:\ell < n\rangle$
witness the desired conclusion in $(*)_\varepsilon$.  So we can carry
the induction and as said above this suffices.

\noindent
4) Let $\cF \subseteq \Pi \bar\delta$ be well orderable $<_D$-cofinal subset so
let $\bar g = \langle g_\alpha:\alpha < \alpha(*) \rangle$ list $\cF$.
\medskip

\noindent
\underline{Case 1}:  $Y \subseteq \Ord$

Let $\mathbf V_1 = \mathbf L[\bar g]$ and $\mathbf V_2 =
\mathbf V_1[D]$, using $D$ as a predicate so $\mathbf
V_1,\mathbf V_2$ are transitive models of ZFC and let $D_2 = D \cap
\mathbf V_2 \in \mathbf V_2$, of course, also $\mathbf V_2 \models
``\theta$ a cardinal $> |Y|"$.

In $\mathbf V_2$ we let $\bar\lambda = \langle \lambda_s:s \in Y\rangle$
be defined by $\lambda_s = \cf(\delta_s)^{\mathbf V_2}$.  Now if $u \in
\mathbf V_2$ is a set of ordinals of cardinality $< \theta$
then the set $\{s:\delta_s > \sup(u \cap \delta_s)\}$ belongs to
$D$ hence to $D \cap \mathbf V_2$; this implies that
$Y_* = \{s \in Y:\lambda_s \ge \theta\}$ belong to $D$.  Now
apply the pcf theorem in $\mathbf V_2$ on $\langle \lambda_s:s \in
Y_*\rangle$ getting $\langle J_{< \mu},Y_\mu:\mu \in \gb\rangle$ and
$\langle g_{\lambda,\alpha}:\lambda \in \gb,\alpha < \lambda\rangle$
where $\ga = \{\lambda_s:s \in Y_*\},\gb = \pcf(\ga)^{\mathbf V_2}$, in
particular such that:
\mn
\begin{enumerate}
\item[$\bullet$]  $\gb = \pcf\{\lambda_s:s \in Y\}$
\sn
\item[$\bullet$]  $Y_\mu \subseteq Y$
\sn
\item[$\bullet$]  $J_{<\mu}$ is the ideal on $Y$ generated by
  $\{Y_\lambda:\lambda \in b \cap \mu\}$
\sn
\item[$\bullet$]  $\langle g_{\lambda,\alpha}:\alpha < \lambda\rangle$
is a sequence of members of $\prod\limits_{s \in Y_*}
\lambda_s,<_{J^+_{<\mu}(Y_* \backslash Y_\mu)}$-increasing and cofinal.
\end{enumerate}
\mn
We can translate this to get a pcf-system for
$(\Pi\bar\delta,<_D)$ in $\mathbf V_2$ hence in $\mathbf V$. 
\medskip

\noindent
\underline{Case 2}:  $Y \nsubseteq \Ord$

We shall show that it essentially suffices to deal with $\bar\delta$
without repetitions.  Note that each $f \in \cF$ or just $f$ a
function from $Y$ into $\Ord$ induces an
equivalence relation $\eq_f$ on $Y_{\mathbf x}:s_1(\eq_f)s_2
\Leftrightarrow f(s_1) = f(s_2) \wedge \delta_{s_1} = \delta_{s_2}$.
For any such equivalence relation $e$ on $Y_{\mathbf x}$, the set $\cF_e =
\{f \in \cF:\eq_f =e\}$ can be translated to one as in Case 1, and if
for some such $e,\cF_e$ is cofinal in $(\Pi\bar\delta,<_{D_{\mathbf x}})$ 
then we are done, but in general this is not clear.  
\Wilog \, $\cE = \{e_f:f \in \cF\}$ is closed under intersection and
assume there is no $e$ as above.  We can define a function $F$ from
$\cE$ into $\alpha(*)$ by $F(e) = \min\{\alpha$: there is no $f \in
\cF$ such that $e_f = e \wedge g_\alpha \le f\}$, it is well defined
by the present assumption and let $u = \Rang(F)$, so $|u| < \hrtg(\cE)
\le \hrtg(\cP(Y \times Y)) \le \theta$, and we can finish easily.

\noindent
5) Let $u_\alpha := \Rang(g_\alpha),v := \{\delta_s:s \in Y\}$ so all
of them are subsets of $\mu$ of cardinality $< \hrtg(Y)$, so $\bar u = \langle
u_\alpha:\alpha <\alpha(*)\rangle$ is well defined and let $\mathbf V'_1
= \mathbf L[\bar u,v]$; it is a well defined universe, a model of ZFC.  In
$\mathbf V'_1$ we define $\bar\delta'$, listing $v$ in increasing
order and $\bar g' = \langle g'_\alpha:\alpha < alpha(*)\rangle$ where
$g'_\alpha = h_{[u_\alpha,\bar\delta']}$.  
In $\mathbf V$ define $\bar f''_\alpha =
\langle g''_\alpha:\alpha < \alpha(*)\rangle$ where $g''_\alpha =
h_{[u_\alpha,\bar\delta]}$.  As $\theta \ge \hrtg(Y)$ clearly $g_\alpha
\le g''_\alpha \mod \cf-\fil_{< \theta}(\bar\delta)$ hence $g_\alpha
\le g''_\alpha \mod D$ hence \wilog \, $\bar g' = \bar g$.  As there
is no real difference between $\bar\delta$ and $\bar\delta'$ and we
can deal with $\bar g',\bar\delta'$ via $\mathbf L[\bar g',\bar\delta']$
as in Case 1 of the proof of part (4) and finish easily.
\end{PROOF}

\begin{discussion}
Alternate proof: suppose we can uniformly choose $\bar f = \langle
f_\alpha:\alpha < \delta_*\rangle$ which is $<_D$-increasing and
cofinal in $(\Pi\bar\delta,<_D)$.

We define an equivalence relation $E$ on $|\cF|$ by: $\alpha E \beta$ \Iff
\, $e_{g_\alpha} = e_{g_\beta}$; let $\bar\beta = \langle \beta_\zeta =
\beta(\zeta):\zeta < \zeta(*)\rangle$ list $\{\alpha < |\cF|:\alpha =
\min(\alpha/E)\}$ in increasing order and let $\zeta:|\cF| \rightarrow
\zeta(*)$ be $\zeta(\alpha) = \min\{\zeta:\alpha \in \beta_\zeta/E\}$.

Let $\bar\xi^* = \langle \xi^*_\zeta:\zeta < \zeta(*)\rangle$
where\footnote{recall the one-to-one function from $\Ord \times \Ord$
  onto $\Ord$ such that $(\alpha_1 \in \alpha \wedge \beta_1 \le
  \beta) \Rightarrow \pr(\alpha_1,\beta_1) \le \pr(\alpha,\beta)$.}
$\xi^*_\zeta = \pr(\otp(\Rang(f_{\alpha_\zeta})),\otp(\Rang(\bar\delta))$ 
and for $\alpha < |\cF|$ let $\hat g_\alpha$ be the function from
$\xi^*_{\zeta(\alpha)}$ to $\Ord$ defined by $\hat g_\alpha(\xi) =
\gamma$ iff for some $s \in Y_{\mathbf x}$ we have $f_\alpha(s) = \gamma
\wedge \xi = \pr(\otp(\Rang(g_{\beta_{\zeta(\alpha)}}) \cap
g_\alpha(s)),\otp(\Rang(\bar\delta) \cap \delta_s))$.

Lastly, let $R = \{(\zeta_1,\zeta_2,\xi_1,\xi_2)$: for some $s \in Y$
for $\ell = 1,2$ we have $\zeta_\ell < \zeta(*),\xi_\ell <
\xi^*_{\zeta_\ell},\xi_\ell = \pr(\otp(\Rang(g_{\alpha_{\zeta_\ell}})
\cap g_{\alpha_{\zeta_\ell}}(s),\otp(\Rang(\bar\delta) \cap \delta_s))\}$.
Now we use $\mathbf V_1 = \mathbf L[\bar\delta,\bar g,E,R,\bar\xi^*]$ let
$\bar D = \langle D_\zeta:\zeta < \zeta(*)\rangle,D_\zeta = D_{\mathbf
  x}(e_{g_{\alpha_\zeta}}),\mathbf V_2 = \mathbf V_1[D_{\mathbf x}]$ and for
$\zeta < \zeta(*)$ let $\bar\lambda_\zeta = \langle
\lambda_{\zeta,\xi}:\xi < \xi_\zeta\rangle,\lambda_{\zeta,\xi} =
\cf(\delta_s)$ when $\xi = \pr(\xi,s)$ for some appropriate
$\varepsilon$.

Clearly $\zeta < \zeta(*) \Rightarrow \xi_\zeta < \theta$, as before
\wilog \, $\lambda_{\zeta,f} = \cf(\lambda_{\zeta,\xi}) \ge \theta$
and $\theta > \hrtg(Y)$ by an assumption hence the pcf analysis in
$\mathbf V_2$ of $\Pi \bar\lambda_\zeta$ is O.K.; moreover and
$\{\lambda_{\eta,\xi}:\xi < \xi_\zeta\}$ does not depend on.

Now the analysis for $\bar\lambda_0$ recalling $\eq_{\bar\delta} =
e_{g_0} = e_{g_{\alpha_0}}$ is enough.
\end{discussion}

\begin{claim}
\label{d12}
If $\mathbf x$ is a pcf-system \then \, there is $\bar Y$ defined
uniformly from $\mathbf x$ such that (so may write $\bar Y^{\mathbf x} = \langle
Y^{\mathbf x}_\varepsilon:\varepsilon < \varepsilon_{\mathbf x}\rangle$):
\mn
\begin{enumerate}
\item[$(a)$]  $\bar Y = \langle Y_\varepsilon:\varepsilon <
\varepsilon(*)\rangle$ 
\sn
\item[$(b)$]  $Y_\varepsilon \subseteq Y_{\mathbf x}$
\sn
\item[$(c)$]  $D^{\mathbf x}_{\varepsilon +1} = D^{\mathbf x}_\varepsilon 
 + Y_\varepsilon$.
\end{enumerate}
\end{claim}

\begin{PROOF}{\ref{d12}}
Fix $\varepsilon < \varepsilon_{\mathbf x}$, if $\varepsilon_{\mathbf x} =
\varepsilon +1$ let $Y_\varepsilon = Y$, hence assume $\varepsilon +1 <
\varepsilon_{\mathbf x}$.  So for some $Y \subseteq Y_{\mathbf x}$ we have
$D_{\mathbf x,\varepsilon +1} = D_{\mathbf x,\alpha} + Y$ hence $E_{\mathbf
  x,\varepsilon} = D_{\mathbf x,\varepsilon} + (Y_{\mathbf x} 
\backslash Y)$; and $f_{\mathbf x,\alpha_{\mathbf x,\varepsilon}}$ 
is a $<_{D_{\mathbf x,\varepsilon +1}}$-upper bound of $\bar f_{\mathbf x} \rest
[\alpha_{\mathbf x,\varepsilon},\alpha_{\mathbf x,\varepsilon +1})$.  But
  $\bar f_{\mathbf x} \rest [\alpha_{\mathbf x,\varepsilon},
\alpha_{\mathbf x,\varepsilon +1})$ is cofinal in $(\Pi \bar
\delta_{\mathbf x},<_{E_{\mathbf x,\varepsilon}})$ hence we can find
$\beta \in [\alpha_{\mathbf x,\varepsilon},\alpha_{\mathbf x,\varepsilon
    +1})$ such that $f_{\mathbf x,\alpha_{\mathbf x,\varepsilon +1}}
< f_{\mathbf x,\beta} \mod E_{\mathbf x,\varepsilon}$.

Let $\beta_*$ be the minimal such $\beta$ and 
easily $Y_\varepsilon := \{s \in Y_{\mathbf x}:f_{\mathbf x,\beta_*}(s) <
g_{\mathbf x},\alpha_{\mathbf x,\varepsilon +1}(s)\}$ is as required.
\end{PROOF}
\bigskip

\subsection {Elaborations} \label{elaborations} \label{2B}\
\bigskip

\begin{claim}
\label{d17}
Assume $\Ax_{4,\lambda,\partial}$.

For any $\lambda$ we can $\partial$-uniformly define the following.

\noindent
1) For $\delta < \lambda$ of cofinality $\aleph_0$, an unbounded
   subset $e_\delta$ of $\delta$ of order type $< \partial$.

\noindent
2) For $\theta = \hrtg(Y),\bar\delta = \langle \delta_s:s 
\in Y\rangle$ a sequence of limit ordinals $< \lambda$ of 
uncountable cofinality satisfying $Y \in \cf-\id_{<
  \theta}(\bar\delta)$, (see \ref{c2}) a closed $u_* \subseteq 
\sup\{\delta_s:s \in Y\}$, unbounded in each $\delta_s$
   of cardinality $< \hrtg([\theta_1]^{< \partial})$ where 
\mn
\begin{enumerate}
\item[$\bullet$]  $\theta_1 = \min\{|u|:(\forall s)[s \in Y 
\rightarrow \delta_s = \sup(u \cap \delta_s)\}$ is necessarily $< \theta$.  
\end{enumerate}
\mn
3) For $\delta < \lambda$, an unbounded subset $e_\delta$ of
   cardinality $< \hrtg([\cf(\delta)]^{\aleph_0})$.
\end{claim}

\begin{PROOF}{\ref{d17}}
1) See \cite{Sh:835} or as 
in the proof of $(*)_4$ inside the proof of \ref{c13}.

\noindent
2) Let $\mathbf U_{\bar\delta} = \{u:u \subseteq \sup\{\delta_s:s \in Y\}$ of
 cardinality $< \theta$ and $u \cap \delta_s$ an unbounded subset of
$\delta_s$ for every $s \in Y\}$.  By the assumption ``$Y \in \cf-\id_{< \theta}
(\bar\delta)$" clearly $\mathbf U_{\bar\delta} \ne \emptyset$, 
hence $\mathbf U'_{\bar\delta} = \{u \in \mathbf
   U_{\bar\delta}:u$ is closed$\}$ is non-empty .  Using $c \ell$ from
   \ref{z8}, the set $u_* = \cap\{c \ell(u):u \in \mathbf
   U'_{\bar\delta}\}$ has cardinality $< \hrtg([\min\{|u|:u \in \mathbf
U_\delta\}]^{< \partial})$.

Now
\mn
\begin{enumerate}
\item[$\bullet_1$]  if $u_n \in \mathbf U'_{\bar\delta}$ for $n <
  \omega$ \then \, $u := \cap \{u_n:n < \omega\}$ belongs to $\mathbf
  U'_\delta$.
\end{enumerate}
\mn
[Why?  Clearly it is a subset of $\mu$ of cardinality $< \theta$,
being $\subseteq u_0$ and it is closed because each $u_n$ is.  But for
any $s \in Y$, why is $u$ unbounded in $\delta_s$?  Because $\delta_s$
has uncountable cofinality
\mn
\begin{enumerate}
\item[$\bullet_2$]  for some $u \in \mathbf U'_{\bar\delta},|u| \le
  \theta_1$ and \wilog \, $u$ is closed, so $|u_*| \le |c
  \ell(u)| \le \hrtg([\theta_1]^{\le \aleph_0})$ as promised.
\end{enumerate}
\mn
By $\bullet_1 + \bullet_2$ we are done.

\noindent
3) By the proof of $(*)_4$ inside the proof of \ref{c13}.
\end{PROOF}

\noindent
We give a sufficient condition for $<_D$-eub existence, try to write
such that we get the trichotomy.
\begin{claim}
\label{d19}
\underline{The eub-existence claim}:

Assume $\Ax_{4,\partial}$ or just $\Ax_{4,\hrtg({}^Y \mu),\partial}$.
The sequence $\bar f$ has a $<_D$-eub (see Definition \ref{z13}(5)), even one
$\partial$-uniformly definable from $(Y,D,\bar f)$ \when \,:
\mn
\begin{enumerate}
\item[$\boxplus$]  $(a) \quad (\theta,Y)$ satisfies clauses
  $(A)(c)(\beta),(\gamma),(\delta)$ of \ref{c13}
\sn
\item[${{}}$]  $(b) \quad D$ is a filter on $Y$, so not
  necessarily $\aleph_1$-complete
\sn
\item[${{}}$]  $(c) \quad \bar f = \langle f_\alpha:\alpha <
  \delta\rangle$
\sn
\item[${{}}$]  $(d) \quad f_\alpha \in {}^Y\Ord$ is $\le_D$-increasing
\sn
\item[${{}}$]  $(e) \quad \cf(\delta) \ge \theta$ and 
$\cf(\delta) \ge \hrtg(\prod\limits_{s \in Y}
  \zeta_s)$ when $\zeta_s < \hrtg(\cP(Y))$ for $s \in Y$.
\end{enumerate}
\end{claim}

\begin{PROOF}{\ref{d19}}
Toward contradiction assume that the desired conclusion fails.
Let $\alpha^*_s = \cup\{f_\alpha(s):\alpha < \delta\}$ for $s \in Y$ and
$\alpha_* = \sup\{\alpha^*_s +1:s \in Y\}$.

We try to choose $g_\zeta$ and $\beta_\zeta <
\delta$ by induction on $\zeta < \hrtg(\cP(Y)/D) \le \hrtg(\cP(Y))$
such that:
\mn
\begin{enumerate}
\item[$\oplus$]  $(a) \quad g_\zeta \in \prod\limits_{s \in
  Y}(\alpha^*_s +1)$
\sn
\item[${{}}$]  $(b) \quad$ if $\alpha < \delta$ then $f_\alpha <
  g_\zeta \mod D$
\sn
\item[${{}}$]  $(c) \quad$ if $\varepsilon < \zeta$ then $g_\zeta \le
  g_\varepsilon \mod D$ and $g_\zeta/D \ne g_\varepsilon/D$
\sn
\item[${{}}$]  $(d) \quad g_\zeta$ and $\beta_\zeta < \delta$ are 
defined as below.
\end{enumerate}
\mn
Clearly impossible as $\cf(\delta) \ge \hrtg(\cP(Y))$ by assumption
$\boxplus(d)$,  so we shall get stuck somewhere.  If 
$\bar g^\zeta = \langle g_\varepsilon:\varepsilon < \zeta\rangle$
is well defined, we let $\bar u_\zeta = \langle u_{\zeta,s}:s \in
Y\rangle$ be defined by $u_{\zeta,s} = \{\gamma$: for some
$\varepsilon < \zeta$ and $n$ we have $\gamma +n = g_\varp(s)$ or $\gamma
+n = \alpha^*_s\}$, so $u_{\zeta,\alpha} \subseteq \alpha^*_s +1$ and
$|u_{\zeta,\alpha}| \le \aleph_0 + |\zeta|$ even uniformly.  Next for
$\alpha < \delta$ we let $f^{\zeta,1}_\alpha \in \prod\limits_{s \in Y}
(\alpha_s +1)$ be defined by $f^{\zeta,1}_\alpha(s) =
\min(u_{\zeta,s}  \backslash f_\alpha(s))$, clearly well defined and
belongs to $\prod\limits_{s \in Y} (\alpha^*_s +1)$ and is
$\le_D$-increasing.  Now
$\{f^{\zeta,1}_\alpha:\alpha < \delta\} \subseteq \prod\limits_{s \in
  Y} u_{\zeta,s}$ so as $\cf(\delta) \ge \hrtg({}^Y(1 +\zeta)) \ge
\hrtg(\prod\limits_{s} u_{\zeta,s})$, necessarily $\langle
f^{\zeta,1}_\alpha/D:\alpha < \delta\rangle$ is eventually constant.
Let $\beta_{\zeta,1} = \min\{\beta < \delta$: if $\alpha \in
(\beta,\delta)$ then $f^{\zeta,1}_\alpha = f^{\zeta,1}_\beta \mod D\}$ so
$\alpha < \delta \Rightarrow f_\alpha \le_D f^{\zeta,1}_{\beta_{\zeta,1}} \mod
D$ and let $g_{\zeta,1} = f^{\zeta,1}_{\beta_{\zeta,1}}$.  If $g_{\zeta,1}$ is a
$<_D$-eub of $\bar f$ we are done, otherwise the construction will split
to cases. 

Let $Y_0 = \{s \in Y:f^{\zeta,1}_{\beta_{\zeta,1}}(s) = 0\},Y_1 = \{s \in
Y:f^{\zeta,1}_{\beta_{\zeta,1}}(s)$ is a successor ordinal$\}$ and $Y_2 =
\{s \in Y:f^{\zeta,1}_{\beta_{\zeta,1}}(s)$ is a limit ordinal of
cofinality $< \theta\}$ and $Y_3 = \{s \in
Y:f^{\zeta,1}_{\beta_{\zeta,1}}(s)$ is a limit ordinal of cofinality
$\ge \theta\}$, so $\langle Y_0,Y_1,Y_2,Y_3\rangle$ is a partition of $Y$
\mn
\begin{enumerate}
\item[$(*)$]  \wilog \, $Y_\ell \in D,g_{\zeta,1}$ is not an lub and
  even $Y_\ell=Y$ from some $\ell < 4$.
\end{enumerate}
\mn
[Why?  For each $\ell < 4$ such that $Y_\ell \in D^+$, clearly we can
replace $D$ by $D+Y_\ell$ hence (by the present assumption) a
$<_{D+Y_\ell}$-eub $g'_\ell$ exists; if $Y_\ell \notin D^+$ let
$g_\ell$ be constantly zero.  Lastly, $\cup\{g^*_\ell \rest
Y_\ell:\ell < 4\}$ is as required.]
\medskip

\noindent
\underline{Case 0}:  $Y_0 \in D$ so $Y_0 = Y$

Trivial.
\medskip

\noindent
\underline{Case 1}: $Y_1 \in D$ so $Y_1 = Y$

Define $g_\zeta \in \prod\limits_{s \in y}(\alpha_s +1)$ by:
$g_\zeta(s) = g_{\zeta,1}(s)-1$.  Clearly it is still a $\le_D$-upper
bound of $\bar f$ as $\bar f$ is $<_D$-increasing, and $g_\zeta <
g_\varepsilon \mod D$ for every $\varepsilon < \zeta$.  Lastly, let
$\beta_\zeta = \beta_{\zeta,1}$. 
\medskip

\noindent
\underline{Case 2}:  $Y_2 \in D$

Let $\langle e_\alpha:\alpha < \alpha_*\rangle$ be as in
\ref{d17}(1),(3) for $\alpha < \delta$, then we
define $f^{\zeta,2}_\alpha \in \prod\limits_{s \in Y_2} (\alpha_s +1)$
by $f^{\zeta,2}_\alpha(s) = \min(e_{g_{\zeta,1}(s)} \backslash 
f_\alpha(s))$ and let $\zeta_s = \otp(e_{g_{\zeta,1}(s)}) < \theta$,
this holds by \ref{c13}(A)(c)$(\beta)$ which in turn holds by
$\boxplus(a)$ of the assumption of the claim.

Now as $\cf(\delta) \ge \hrtg(\prod\limits_{s \in Y_2} \zeta_s) 
= \hrtg(\prod\limits_{s \in Y_2} e_{g_{\zeta,1}(s)})$ 
clearly $\langle f^{\zeta,2}_\alpha/D:\alpha <
  \delta\rangle$ is eventually constant, so $\beta_{\zeta,2} =
  \min\{\beta < \delta$: if $\alpha \in (\beta,\delta)$ then
  $f^{\zeta,2}_\alpha/D = f^{\zeta,2}_\beta/D\}$ is well defined.  Let
  $\beta_\zeta = \sup(\{\beta_{\zeta,1},\beta_{\zeta,2}\} \cup
  \{\beta_\varepsilon +1:\varepsilon < \zeta\})$ it is $<
  \delta,\cf(\delta) > |\zeta|$ and let $g_\zeta =
  f^{\zeta,2}_{\beta_\zeta}$.  Clearly $\varepsilon < \zeta
  \Rightarrow g_\zeta = f^{\zeta,7}_{\beta_\zeta} <
  f^{\zeta,1}_{\beta_{\zeta,1}} \le g_\varepsilon \mod D$, so
  $(g_\zeta,\beta_\zeta)$ are as required.
\medskip

\noindent
\underline{Case 3}:  $Y_3 = Y$

Let $\bar f' = \langle f'_\alpha:\alpha < \delta\rangle,f'_\alpha \in
\prod\limits_{s} g_{\zeta,1}(s)$ defined as $f_\alpha(s)$ if $<
g_{\zeta,1}(s)$, zero otherwise.

Now $g_{\zeta,1}$ is not a $<_D$-eub of $\bar f$ hence there is $h \in
{}^Y\Ord$ such that $h < g_{\zeta,1} \mod D$ and for no $\alpha <
\delta$ do we have $h < f_\alpha \mod D$.  But $h$ was not canonically
chosen.  Clearly the assumption of \ref{d5}, i.e. \ref{c17} holds with
$Y,\theta,g_{\zeta,1},\bar f'$ here standing for
$Y,\theta,\bar\delta,\bar f$ here.  So there is a $\pcf$-system $\mathbf
x$ with $Y_{\mathbf x} = Y,\theta_{\mathbf x} = \theta,D_{\mathbf x} =
D,\bar f_{\mathbf x} = \bar f'$ and $\bar\delta_{\mathbf x} =
g_{\zeta,1}$.

Hence by \ref{c24}(1) we can define a pair
$(\cF,<_*)$ such that $\cF \subseteq \prod\limits_{s \in Y_2}
g_{\zeta,1}(s)$ is cofinal and $<_*$ a well ordering of $\cF$.

So as $g_{\zeta,1}$ is not a $<_D$-eub of $\bar f$ there is $h
\in \cF$ witnessing this and let $h_* \in \cF$ be the $<_*$-first one.

Let

\begin{equation*}
\begin{array}{clcr}
\beta_{\zeta,3} = \min\{\beta <\alpha: &\text{ if } \alpha \in
(\beta,\delta) \text{ then } \{s \in Y:f_\alpha(s) \le g_\zeta(h_*((s))\} = \\
  &\{s \in Y:f_\beta(s) \le h_*(s)\} \mod D + Y_2\},
\end{array}
\end{equation*}

\mn
well defined as before.  Lastly, let $g_\zeta \in {}^Y\Ord$ be defined
as follows: $g_\zeta(s)$ is
\mn
\begin{enumerate}
\item[$\bullet$]  $h_*(s)$ if $f_{\beta_{\zeta,3}}(s) \le h_*(s)$
\sn
\item[$\bullet$]  $f_{\beta_{\zeta,3}}(s)$ if $f_{\beta_{\zeta,3}}(s)> h_*(s)$.
\end{enumerate}
\end{PROOF}
\bigskip

\centerline {$* \qquad * \qquad *$}
\bigskip

Now we give a version of the main theorem of \cite[\S1]{Sh:835}.  From
this we may try to understand better ${}^\kappa \lambda$ and use it in
constructions, i.e. to diagonalize.
\begin{theorem}
\label{d29}
$[\Ax_{4,\lambda,\partial}]$

For $\kappa < \lambda$ letting $X_\kappa = {}^\omega(\Fil^4_{\aleph_1}
(\kappa))$, we can $\partial$-uniformly define $\langle (\cS_t,<_t):t \in
X_\kappa \rangle$ such that:
\mn
\begin{enumerate}
\item[$(a)$]  $\cup\{\cS_t:t \in X_\kappa\} = {}^\kappa \lambda$
\sn
\item[$(b)$]  $<_t$ is a well ordering of $\cS_t$
\sn
\item[$(c)$]  there is an equivalence relation $E$ on ${}^\kappa
\lambda$ such that:
\mn
\begin{enumerate}
\item[$(\alpha)$]   $({}^\kappa \lambda)/E$ is well ordered
\sn
\item[$(\beta)$]  each equivalence class is of power $\le X_\kappa$
\end{enumerate}
\sn
\item[$(d)$]  moreover for some $\bar g = \langle g_{\bar{\gy},\alpha}:
\bar{\gy} \in  X_\kappa,\alpha \in S_{\bar{\gy}}\rangle$ and $\bar S =
\langle S_{\bar{\gy}}:\bar{\gy} \in X_\kappa\rangle$ and 
$\bar{\cF} = \langle \cF_\beta:\beta < \beta(*)\rangle$ we have
\sn
\begin{enumerate}
\item[$(\alpha)$]  $\beta(*) < \hrtg(\alpha(*)]^{\aleph_0}$ where
  $\alpha(*) = \sup\{\rk_D(\lambda):D \in \Fil^1_{\aleph_1}(Y)\}$
\sn
\item[$(\beta)$]  $\beta(*) = \cup\{S_{\bar{\gy}}:\bar{\gy} \in
  X_\kappa\}$
\sn
\item[$(\gamma)$]  $\{g_{\bar{\gy},\alpha}:\bar{\gy} \in X_\kappa,\alpha
  \in S_{\bar{\gy}}\}$ is equal to ${}^\kappa \lambda$
\sn
\item[$(\delta)$]  $g_{\bar{\gy}_1,\alpha_1} =
 g_{\bar{\gy}_2,\alpha_2}$ implies $\alpha _1 = \alpha_2$
\sn
\item[$(\varepsilon)$]  $\bar{\cF} = \langle \cF_\beta:\beta <
  \beta(*)\rangle$ is a partition of ${}^\kappa \lambda$
\sn
\item[$(\zeta)$]  $|\cF_\beta| \le_{\qu} |X_\kappa|$.
\end{enumerate}
\end{enumerate}
\end{theorem}

\begin{remark}
1) We may compare with \cite[\S1]{Sh:835}.

\noindent
2) Recall \ref{z24}(2).
\end{remark}

\begin{PROOF}{\ref{d29}}
Fix a witness $c \ell$ of $\Ax_{4,\lambda,\partial}$. 
For every $\gy \in \Fil^4_{\aleph_1}(Y)$ and ordinal
$\alpha$ there is at most one $f \in {}^Y(\lambda +1)$ such that 
$f$ satisfies $\alpha = \rk_D(f)$ and so $f \rest 
(Y \backslash Z_{\gy})$ is constantly zero
and $D^{\gy}_2 = \dual(J[f,D^{\gy}_1])$, see \ref{z15}, \ref{z21};
if in this case call it $f_{\gy,\alpha}$
and let $S_{\gy,\lambda}$ be a set of $\alpha$ such that
$f_{\gy,\alpha}$ is well defined.

So $\langle f_{\gy,\alpha}:\gy \in \Fil^4_{\aleph_1}(Y),\alpha \in
S_{\gy,\alpha}\rangle$ is well defined.  For every $f \in \lambda$ 
and $\aleph_1$-complete filter $D_1$ on $Y$ for some
$\gy \in \Fil^4_{\aleph_1}(Y)$ satisfying $D_{\gy,1} = D_1$ 
and ordinal $\alpha$ we have $f = f_{\gy,\alpha} \mod D_{\gy,2}$ 
(in fact $\alpha = \rk_{D_1}(f) < \rk_D(\lambda) \le
\alpha(*),\alpha(*)$ from $(d)(\alpha)$ of the Theorem).

Now
\mn
\begin{enumerate}
\item[$(*)_1$]  for every $f \in {}^Y(\lambda +1)$ 
there is a countable set $\gY \subseteq \Fil^4_{\aleph_1}(Y)$ such
that
\sn
\begin{enumerate}
\item[$(\alpha)$]   $f$ semi-satisfies each $\gy \in \gY$ 
\sn
\item[$(\beta)$]   $Y = \cup\{Z_{\gy}:\gy \in \gY\}$
\sn
\item[$(\gamma)$]  for each $\gy \in \gY$, for some $\alpha$ we have
  $f \rest Z_{\gy} = f_{\gy,\alpha} \rest Z_{\gy}$.
\end{enumerate}
\end{enumerate}
\mn
[Why?  Let $\cZ = \{Z_{\gy}:\gy \in \Fil^4_{\aleph_1}(Y)$
  and for some $\alpha \in S_{\gy,\lambda}$ we have $f \rest Z_{\gy} =
  f_{\gy,\alpha} \rest Z_{\gy}\}$.  If $Y$ is the union of a countable
  subset of $\cZ$ then recalling $\AC_{\aleph_0}$ we have
$Y = \cup\{Z_{\gy_n}:n\}$ for some $\{\gy_n:n <
  \omega\} \subseteq \Fil^4_{\aleph_1}(Y)$ and we are easily done.  If not,
  $D_1 := \{Z \subseteq Y:Z$ includes $(Y \backslash
  \bigcup\limits_{n} Z_{\gy}$: for some $\langle \gy_n:n <
\omega\rangle \in {}^\omega(\Fil^4_{\aleph_0}(Y))$ satisfying
$Z_{\gy_n} \in \cZ$ for $n < \omega\}$ is an $\aleph_1$-complete 
filter and we easily get a contradiction.]

Recall $S_{\gy,\lambda} = \{\alpha < \alpha(*):f_{\gy,\alpha}$ well defined$\}$
and by $\Ax_4$ we can find a list $\langle \eta_\beta:\beta <
\beta(*)\rangle$ of $\{\eta:\eta \in {}^\omega \alpha(*)\},\beta(*) <
\hrtg({}^\omega \beta(*))$ and even $\beta(*) = |\beta(*)|^{\aleph_0}\}$.

Now for every $\bar{\gy} \in X_\kappa := {}^\omega(\Fil^4_{\aleph_1}(Y))$, 
let $W_{\bar{\gy}} = \{\beta < \beta(*):\eta_\beta(n) 
\in S_{\gy_n,\lambda}$ for each $n$
and $\cup\{f_{\gy_n,\eta_n(\alpha)} \rest Z_{\gy_n}:
n < \omega\}$ is a function, in
fact one from $Y$ to $\lambda +1\}$.  For $\beta \in W_{\bar{\gy}}$ let
$g_{\bar{\gy},\beta}$ be $\cup\{f_{\gy_n,\eta_\beta(n)}:n < \omega\}$
and let $S_{\bar{\gy}} = \{\beta \in W_{\bar{\gy}}:g_{\bar{\gy},\beta} \notin
\{g_{\bar{\gz},\gamma}:\bar{\gz} \in X_\kappa$ and $\gamma <
\beta\}\}$.

Note that
\mn
\begin{enumerate}
\item[$(*)_2$]  $(a) \quad \langle S_{\bar{\gy}}:\bar{\gy} \in
  X_\kappa\rangle$ exist
\sn
\item[${{}}$]  $(b) \quad \bigcup\limits_{\gy} S_{\bar{\gy}} \subseteq
  \beta(*)$
\sn
\item[${{}}$]  $(c) \quad \langle S_{\bar{\gy}}:\bar{\gy} \in
  X_\kappa\rangle$ exists and $\cup\{S_{\bar{\gy}}:\bar{\gy} \in
  X_\kappa\} = \beta(*)$.
\end{enumerate}
\mn
Note also that clause (d) of the theorem implies clauses (a),(b); (let
$\cS_{\bar{\gy}} = \{g_{\bar{\gy},\alpha}:\alpha \in S_{\bar{\gy}}\}$
and $<_{\bar{\gy}} =
\{(g_{\bar{\gy},\alpha},g_{\bar{\gy},\alpha},g_{\bar{\gy},\beta}):\alpha
< \beta$ are from the set $S_{\bar{\gy}}$ of ordinals).

Also clause (d) implies clause (c) letting $E =
\{(g_{\bar{\gy}_1,\alpha_1},g_{\bar{\gy}_2,\alpha_2}):\bar{\gy}_\ell
\in X_\kappa,\alpha_\ell \notin S_{\bar{\gy}_\ell}$ for $\ell=1,2$ and
$\alpha_1 = \alpha_1\}$ recalling $(d)(\delta)$.

So it is enough to prove clause (d).

Now
\mn
\begin{enumerate}
\item[$\bullet$]  clause $(d)(\alpha)$ holds by the choices of
  $\alpha(*),\beta(*)$
\sn
\item[$\bullet$]  clause $(d)(\beta)$:  we have
  only $\beta(*) \supseteq \cup\{S_{\bar{\gy}}:\bar{\gy} \in
  X_\kappa\}$, but we can replace $\beta(*)$ by
  $\otp(\cup\{S_{\bar{\gy}}:\bar{\gy} \in X_\kappa\}$
\sn
\item[$\bullet$]  clause $(d)(\gamma)$: $g_{\bar{\gy},\beta} \in
  {}^\kappa \lambda$ are defined above but why ${}^\kappa \lambda =
\{g_{\bar{\gy},\alpha}:\bar{\gy} \in X_\kappa,\alpha \in
  S_{\bar{\gy}}\}$?  As said above, if $f \in {}^\kappa \lambda$ by
  $(*)_1$ there is a countable $\gY \subseteq
\FIL^4_{\aleph_1}(Y)$ as there, hence for some sequence $\langle
  (\gy_n,\alpha_n):n < \omega\rangle$ we have $\gY = \{\gy_n:n <
  \omega\}$ and $f \rest Z_{\gy_n} = f_{\gy_n,\alpha_n} \rest
  Z_{\gy}$.  Hence $\bar{\gy} := \langle \gy_n:n < \omega\rangle \in
  X_\kappa$ and for some $\gamma < \beta(*)$ we have $\eta_\gamma =
  \langle \alpha_n:n < \omega\rangle$.  So $f =
  \cup\{f_{\gy_n,\eta_\gamma(n)} \rest Z_{\gy_n}:n < \omega\} =
  g_{\bar{\gy},\beta}$ so $f \in W_{\bar{\gy}}$, and $f =
  g_{\bar{\gy},\gamma}$, hence by the choice of $S_{\bar{\gy}}$ there
  are $\bar{\gz} \in X_\kappa$ and $\beta^* \le \gamma$ such that
  $\beta \in w'_{\bar{\gz}}$ and $f = g_{\bar z,\beta}$, so we are
  done
\sn
\item[$\bullet$]  clause $(d)(\delta)$: look again at the choice of
  $S_{\bar{\gy}}$
\sn
\item[$\bullet$]  clause $(d)(\varp)$: let $\cF_\beta =
  \{g_{\bar{\gy},\beta}:\bar{\gy} \in X_\kappa$ and $\beta$ belongs to
  $S_{\bar{\gy}}\}$
\sn
\item[$\bullet$]  clause $(d)(\zeta)$: check.
\end{enumerate}
\end{PROOF}

\begin{conclusion}
\label{d31}
Assume $\Ax_{4,\partial}$.  If $\partial \le \kappa < \mu$ and
$\hrtg(\Fil^4_{\aleph_1}(\kappa)) < \mu$.  \Then \, the following
   cardinals are almost equal (as in \cite[\S(3A)]{Sh:955}:
\mn
\begin{enumerate}
\item[$(a)$]  $\hrtg({}^\kappa \mu)$
\sn
\item[$(b)$]  $\wlor({}^\kappa \mu)$
\sn
\item[$(c)$]  o-Depth$^+_\kappa({}^\kappa \mu) = 
\sup\{{\rm o}-\Depth^+_D(\mu):D$ a filter$\}$.
\end{enumerate}
\end{conclusion}

\begin{PROOF}{\ref{d31}}
By \ref{d29}.
\end{PROOF}

A drawback of the pcf theorem is the demand
$\theta \ge \hrtg(\Fil^4_{\aleph_1}(Y))$ rather than just $\theta \ge
\hrtg(\cP(Y))$ or even $\theta \ge \hrtg(Y)$.  Note: in
\cite[Ch.XII,\S5]{Sh:b} we work to assume just the parallel of $\theta
\ge \hrtg(\cP(Y))$, i.e. $\Min(\ga) > 2^{|\ga|}$ rather than the
parallel of $\theta \ge \hrtg(\cP(\cP(Y))$, i.e. $\Min(\ga) > 2^{2^{|\ga|}}$
and only in \cite{Sh:345} we succeed to use just
the parallel of $\theta \ge \hrtg(Y)$.

We may try to analyze not $\Pi \bar\delta,\bar\delta = \langle
\delta_s:s \in Y\rangle$ but rather all $\Pi(\bar\delta \rest Z),Z \in
\cA$ simultaneously where $\cA \subseteq \cP(Y)$, demanding $Z \in
\cA \Rightarrow \theta \ge \hrtg(\Fil^4_{\aleph_1}(Z))$ but less on
$|Y|$; hopefully see \cite{Sh:F1303}.

\noindent
We may consider
\begin{definition}
\label{d34}
Let $\Ax_{5,F}$ say: if $Y = \kappa \in \Card$ then
$\Ax_{5,\kappa,F(\kappa)}$ where $\Ax_{5,Y,\theta}$ means that: if
$\bar\delta = \langle \delta_s:s \in Y\rangle$ is a sequence of limit
ordinals and $D = \cf-\fil_{< \theta}(\bar\delta)$ \then \, there is a
pcf-system $\mathbf x_{\bar\delta}$ for $(\Pi \bar\delta,<_D)$, see
\ref{c59}.  Moreover, the choice of $\mathbf x_\delta$ is $\partial$-uniform.
\end{definition}

\begin{definition}
\label{s27}
1) We say $\mathbf p$ is a pcf-problem \when \, it consists of:
\mn
\begin{enumerate}
\item[$(a)$]  $\bar\delta = \langle \delta_s:s \in Y\rangle$
  and $\mu = \sup\{\delta_s:s \in Y\}$ and $\cA \subseteq \cP(Y)$
\sn
\item[$(b)$]  $D_* = D_{\mathbf p}$ is a filter on $Y$, it may be $\{Y\}$
\sn
\item[$(c)$]  $\theta = \theta[Y,\bar\delta,D_*] =
  \theta[Y,\delta,D_*,\partial]$ is any cardinal satisfying:
\sn
\begin{enumerate}
\item[$(\alpha)$]  $\cf-\id_{< \theta}(\bar\delta)
  \subseteq \dual(D_*)$, note that this
holds when each $\delta_s$ is an
ordinal $\le \mu$ of cofinality $\ge \theta$, see below
\sn
\item[$(\beta)$]  $\alpha < \theta \Rightarrow
\hrtg([\alpha]^{\aleph_0} \times \partial) \le \theta$ so 
$\partial < \theta$ and so if $\Ax_4$ then 
the demand is equivalent to ``$\partial < \theta$ and
$\alpha < \theta \Rightarrow |\alpha|^{\aleph_0} < \theta$"
\sn
\item[$(\gamma)$]  $\hrtg(\Fil^4_{\aleph_1}(Z))\le \theta$ for every
  $Z \in \cA$.
\end{enumerate}
\end{enumerate}
\mn
2) For $\mathbf p$ a pcf-problem let $\bar\delta_{\mathbf p} =
\delta,\delta_{\mathbf p,s} = \delta_s$, etc., if clear from the context
$\mathbf p$ is omitted.

\noindent
3) For $D$ a filter on $Y_{\mathbf p}$ extending $D_{\mathbf p}$ let $c
\ell_{\mathbf p}(D) = c \ell(D,\mathbf p) = 
\{A \subseteq Y_{\mathbf p}$: if $Z \in \cA_{\mathbf p}$ 
then $A \cup (Y_{\mathbf p} \backslash Z) \in D$.

\noindent
4) $\mathbf p$ is nice if $\hrtg(\cP(Y)) \le \theta_{\mathbf p}$.
\end{definition}

\begin{definition}
\label{s35}
We say $\mathbf x$ is a wide pcf system \when \, $\mathbf x$ consists of
(if we omit $(f),(g)(\alpha),(\beta)$, (as in \ref{d2}(3)) we say
``almost wide"):
\mn
\begin{enumerate}
\item[$(a)$]  $\mathbf p$, a pcf-problem let $D_{\mathbf x} = D_{\mathbf p},
\theta = \theta_{\mathbf p}$, etc.
\sn
\item[$(b)$]   an ordinal $\varepsilon_{\mathbf x} = \varepsilon(\mathbf x)$
\sn
\item[$(c)$]  $\bar\alpha^* = \langle \alpha^*_\varepsilon:\varepsilon
  \le \varepsilon_{\mathbf x}\rangle$ is increasing continuous
\sn
\item[$(d)$]  $(\alpha) \quad \bar D = \langle D_\varepsilon:\varepsilon \le
  \varepsilon_{\mathbf x}\rangle$ is a continuous sequence
of filters on $Y$ except that 

\hskip25pt  possibly $D_{\varepsilon_{\mathbf x}} = \cP(Y)$
\sn
\item[${{}}$]  $(\beta) \quad D_\varepsilon = c \ell_{\mathbf
  p}(D_\varepsilon)$
\sn
\item[${{}}$]  $(\gamma) \quad$ for limit $\varepsilon,D_\varepsilon 
= c \ell_{\mathbf p}(\bigcup\limits_{\zeta < \varepsilon} D_\zeta)$
\sn
\item[$(e)$]  $D_0 = D_{\mathbf x}$ is $\cf-\fil_\theta(\bar\delta)$
\sn
\item[$(f)$]  $\bar E = \langle E_\varepsilon:\varepsilon <
  \varepsilon_{\mathbf x}\rangle$
\sn
\item[$(g)$]  for each $\varepsilon <
  \varepsilon_{\mathbf x} < \theta$ there is $A_\varepsilon \in
  D^+_\varepsilon$ such that
\sn
\item[${{}}$]  $(\alpha) \quad D_{\varepsilon +1} = D_\varepsilon +
  A_\varepsilon$
\sn
\item[${{}}$]  $(\beta) \quad E_\varepsilon = D_\varepsilon + 
(u \backslash A_\varepsilon)$
\sn
\item[${{}}$]  $(\gamma) \quad$ there are $a_\varepsilon \subseteq
  \kappa$ and $h_\varepsilon \in \prod\limits_{i \in \varepsilon} u_i$
  such that $\{(i,h_\varepsilon(i)):i \in a_\varepsilon\} \notin D_\varepsilon$
\sn
\item[${{}}$]  $\hskip10pt \bullet \quad$ but $A_\varepsilon$ is not
  necessarily unique, only $A_\varepsilon/D_\varepsilon$ is, and of
  course, also 

\hskip25pt $a_\varepsilon,h_\varepsilon$ are not necessarily unique
\sn
\item[${{}}$]  $(\delta) \quad$ there is $Z \in \cA$ such that $Z \in
  \dual(D_{\varepsilon +1}) \backslash \dual(D_\varepsilon)$
\sn
\item[$(h)$]  $\bar f = \langle f_\alpha:\alpha < \varepsilon_{\mathbf
  x}\rangle,f_\alpha \in \Pi \bar\delta$
\sn
\item[$(i)$]  $\bar f \rest \alpha_{\varepsilon +1}$ is
  $\le_{D_\varepsilon}$-increasing 
\sn
\item[$(j)$]  $\bar f \rest [\alpha_\varepsilon,\alpha_{\varepsilon
    +1})$ is $<_{E_\varepsilon +Z}$-cofinal for some $Z \in D^+_\varepsilon$.
\end{enumerate}
\end{definition}

\begin{theorem}
\label{s38}
Assume $\Ax_{4,\partial}$.
Assume $\mathbf p$ is a pcf-problem and $\hrtg(\cA_{\mathbf p}) \le
\theta_{\mathbf p},\partial < \theta_{\mathbf p}$.  
\Then \, there is a wide pcf-system $\mathbf x$ such
that $\mathbf p_{\mathbf x} = \mathbf p$.
\end{theorem}

\begin{PROOF}{\ref{s38}}
As in \S1 we try to choose $\alpha_\varepsilon$ and $\langle
f_\alpha:\alpha \le
\alpha_\varepsilon\rangle,D_\varepsilon,E_\varepsilon$ by induction on
$\varepsilon$ satisfying the relevant demands.  The main point is
having chosen $\langle \alpha_\xi,D_\xi:\xi \le \zeta \rangle,
\langle f_\alpha:\alpha \le \alpha_\zeta \rangle$, we try to choose
for $\varepsilon = \zeta +1$.  So we try to choose $f_\alpha$
for $\alpha > \alpha_\zeta$ by induction on $\alpha$ satisfying the relevant
conditions.  Arriving to limit $\alpha$ let
$\cA^1_\alpha := \{Z \in \cA:Z \notin \dual(D_\varepsilon)\}$ and
$\cA^2_\alpha = \{Z \in \cA^1_\alpha:\langle f_\beta:\beta <
\alpha\rangle$ has a $<_{D_\varepsilon +Z}$-upper bound in
$\Pi\bar\delta\}$.  If $\cA^1_\alpha = \emptyset$ we are done.  If
$\cA^2_\alpha \ne \emptyset$ by \S1 we can define $\langle
f_{\alpha,Z}:Z \in \cA^2_\alpha\rangle$ such that $f_{\alpha,Z} \in
\Pi \bar\delta$ is an $<_{D_\varepsilon +Z}$-upper bound of $\langle
f_\beta:\beta < \alpha\rangle$ and let $f_\alpha \in \Pi\bar\delta$ be
defined by $f_\alpha(s) = \sup\{f_{\alpha,Z}(s):Z \in \cA^2_\alpha\}$
if $< \delta_s$ and zero otherwise.  As $\theta \ge 
\hrtg(\cA_{\mathbf p}) \ge \hrtg(\cA^2_\alpha)$, clearly 
$\beta < \alpha \wedge Z \in
\cA^2_\alpha \Rightarrow f_\beta < f_\alpha \mod (D_\varepsilon +Z)$.
If $\cA^2_\alpha = \cA^1_\alpha \ne \emptyset$, then $f_\alpha$ is as
required as we are assuming $D_\varepsilon = c \ell_{\mathbf
  p}(D_\varepsilon)$.  If $\cA^2_\alpha \ne \cA^1_\alpha$, let
$\alpha_{\varepsilon +1} = \alpha$ and $f_\alpha$ is as required.
\end{PROOF}
\bigskip

\subsection {True successor cardinals} \label{2C}\
\bigskip

Contrary to our ZFC intuition, without full choice successor
cardinals, may be singular.  On history we may start with Levy proving
$\ZF + ``\aleph_1$ is singular" is consistent and end with Gitik
proving $\ZF + (\forall \lambda),\cf(\lambda) = \aleph_0$ is
consistent, using suitable large cardinals.  Note: ``two successive cardinals
are singular" has quite high consistency strength.

A major open question is whether $\ZF + \DC + (\forall
\lambda)(\cf(\lambda) \le \aleph_1)$ is consistent. But when $\ZF +
\DC + \Ax_4$ holds the situation is very different.  Also contrary to
our ZFC intuition, successor cardinals may be measurable.

For a cardinal to be a true successor is saying it fits our ZFC
intuition.  In particular, it avoids the two anomalities mentioned
above, and eventually itwill enable us to carry various constructions;
all this motivates Question \ref{p18}.

We continue the investigation in \cite{Sh:835} of successor of
singulars, not relying on \cite{Sh:835}.
\begin{definition}
\label{p16}
1) We say $\lambda$ is a true successor cardinal \when \, for some
cardinal $\mu,\lambda = \mu^+$ and we have a witness $\bar f$,
   which means $\bar f = \langle f_\alpha:\alpha \in
   [\mu,\lambda)\rangle$ and $f_\alpha$ is a one-to-one function from
   $\alpha$ into $\mu$.

\noindent
1A) We say $\bar f$ is an onto-witness when each $f_\alpha$ is onto
$\mu$, see \ref{p21}(1) below.

\noindent
2) We say a set 
$\cU \subseteq \Ord$ is a smooth set \when \, there is a witness
   $\bar f$ which means that $\bar f = \langle f_\alpha:\alpha \in
   \cU\rangle,f_\alpha$ is a one-to-one function 
from $\alpha$ onto $|\alpha|$. 
\end{definition}

\noindent
We may naturally ask
\begin{question}
\label{p18}
Assume, e.g. $\ZF + \DC + \Ax_4$.

\noindent
1) Is there a class of successor of regular cardinals which are true
successor cardinal?  See \ref{p21}(2).

\noindent
2) Assume $\mu$ is strong limit (i.e. $\alpha < \mu \Rightarrow
   \hrtg(\cP(\mu)) < \mu)$ of cofinality $\aleph_1$, so
   $\mu^+$ is regular, but assume in addition that
$\mu^{++}$ is regular $< \pp(\mu)$, see\footnote{generality with weak
  choice there is a choice to be made, but assuming $\Ax_4$ or so and
  $\cf(\mu) = \aleph_0$, there is no problem} \cite[Ch.II]{Sh:g}.  
Is $\mu^{++}$ truely successor?

\noindent
3) Assume $\mu$ is strong limit of cofinality $\aleph_0$ and $\mu^{+2}$
is singular, is $\mu^{+3}$ a true successor cardinal?
\end{question}

\begin{claim}
\label{p21}
1) If $\lambda$ is a true successor, \then \, $\lambda$ is regular and has
   an onto-witness (computed uniformly from a witness).

\noindent
2) [$\Ax^4_{\mu^+}$ or just $\Ax_{4,\mu^+,\partial}$]
Assume $\mu$ is singular and $(\forall \alpha < \mu)
(\hrtg([\alpha]^{\aleph_0} \times \partial) < \mu)$.  \Then \, 
$\mu^+$ is a true successor cardinal.

\noindent
3) [$\Ax_{4,\lambda}$ or just $\Ax_{4,\lambda,\partial}$] 
The set $\cU$ of ordinals $\alpha < \lambda$ such that $|\alpha|$ is 
singular and $(\forall \beta < |\alpha|)
[\hrtg([\beta]^{\aleph_0} \times \partial) \le |\alpha|]\}$ 
is a smooth set of ordinals.

\noindent
4) For every ordinal $\alpha_*,\alpha_* \in
\cf-\id_{\langle(\hrtg([\cf(\alpha)]^{\aleph_0} \times \partial):\alpha <
\alpha_*\rangle}(\langle \alpha:\alpha < \alpha_*\rangle)$. 
\end{claim}

\begin{PROOF}{\ref{p21}}
Let $\pr$ be the classical one-to-one function from $\Ord \times \Ord$
onto $\Ord$ such that $\pr(\alpha,\beta) < (\max\{\alpha,\beta\})^2$
and $\pr_\mu = \pr \rest (\mu \times \mu)$.

\noindent
1) Let $\bar f = \langle f_\alpha:\alpha \in [\mu,\mu^+)\rangle$ witness
$\lambda$ is truely a successor.  First define, for $\alpha \in
  [\mu,\mu^+)$ a fucntion $f'_\alpha:\alpha \rightarrow \mu$ by
    $f'_\alpha(\beta) = \otp(\Rang(f_\alpha) \cap f_\alpha(\beta))$;
    obviously it is a one-to-one function from $\alpha$ into $\mu$
    with range an initial segment; but $|\Rang(f'_\alpha)| = |\alpha|
    = \mu$ so $\Range(f'_\alpha) = \mu,\langle f'_\alpha:\alpha \in
    [\mu,\mu^+)\rangle$ is as promised.

Second proving $\lambda$ is regular, 
toward contradiction let $\cU$ be such that $\cU \subseteq
   \lambda = \sup(\cU),\cU \cap \mu = \emptyset$ 
and $\otp(\cU) < \lambda$, so \wilog \, $\le
   \mu$.  Now we shall combine $\langle f_\alpha:\alpha \in \cU \rangle$ to get
$|\lambda| \le \mu$ by getting a one to one function $f$ from $\lambda$
   into $\mu \times \mu$; for $i < \lambda$ let $\alpha_i =
   \min\{\alpha \in \cU:\alpha > i\}$ and define $f(i) = \pr(\otp(\cU \cap
\alpha_i),f_{\alpha_i}(\alpha))$.  So $f$ exemplifies $|\lambda| \le
   |\mu \times \mu|$ but the latter is $\mu$, contradiction.

\noindent
2) By part (3) applied to $\cU = [\mu,\mu^+)$.

\noindent
3) Let $\cS \subseteq [\lambda]^{< \partial}$ witness
$\Ax_{4,\lambda,\partial}$ and $<_*$ a well ordering of $\cS$. 
Let $\alpha_* = \cup\{\alpha +1:\alpha \in \cU\}$ let 
$c \ell:[\alpha_*]^{\aleph_0} \rightarrow \alpha_*$ be as in
   \ref{z8}, let $<_*$ be a well order $\cS$ and let
   $u_\beta$ for $\beta < \alpha_*$ be defined by
\mn
\begin{enumerate}
\item[$\bullet$]  if $\beta=0$ then $u_0 = \emptyset$
\sn
\item[$\bullet$]  if $\beta= \gamma +1$ then $u_\beta = \{\gamma\}$
\sn
\item[$\bullet$]  if $\cf(\beta) > \aleph_0$ then $u_\beta =
\cap\{\cup\{c \ell(v):v \in [u]^{\aleph_0}\}:u$ a club of $\beta\}$
\sn
\item[$\bullet$]  if $\cf(\beta) = \aleph_0$ the $u_\beta = v_\beta
  \cap \beta$ where $v_\beta$ is the $<_*$-first $v \in \cS$ such that
  $\beta = \sup(v \cap \beta)$.
\end{enumerate}
\mn
Now choose $f_\alpha$ for $\alpha \in \cU$ by induction
on $\alpha$ using $\pr_{|\alpha|}$ as in the proof of part (2).

\noindent
4) By $(*)_4$ in the proof of \ref{c13}, in particular, $(c)_2$ there.
\end{PROOF}

\noindent
Recalling $\cf-\id_{< \gamma}(\bar\delta)$ from Definition \ref{c2}.
\begin{claim}
\label{p23}
1) If $\lambda = \mu^+$ \then \, $\lambda$ is a true successor iff
$\lambda \in \cf-\id_{<(\mu +1)}(\lambda)$, (which means $\lambda \in
\cf-\id_{< (\mu+1)}(\langle \alpha:\alpha < \lambda\rangle))$ iff
$\lambda \in \cf-\id_{< \gamma}(\langle \alpha:\alpha < \lambda\rangle)$
for some $\gamma < \lambda$.

\noindent
2) When $\mu$ is singular, we can add: iff
$\lambda \in \cf_{< \mu}(\langle \alpha:\alpha < \lambda\rangle)$.
\end{claim}

\begin{PROOF}{\ref{p23}}

\noindent
1) \underline{First condition implies second condition}:

So assume $\lambda$ is a true successor, let $\langle f_\alpha:\alpha
\in [\mu,\mu^+)\rangle$ witness it.  For each $\alpha < \mu^+ =
  \lambda$ we choose $u_\alpha$ as follows:
\medskip

\noindent
\underline{Case 1}:  $u_\alpha = \alpha$ if $\alpha < \mu$
\medskip

\noindent
\underline{Case 2}:  $\alpha \ge \mu$

For any $j < \mu$ let $\cU_{\alpha,j} = \{\beta < \alpha:f_\alpha(\beta)
< j\}$, so $\langle \cU_{\alpha,j}:j < \mu\rangle$ is
$\subseteq$-increasing with union $\alpha$ and $|\cU_{\alpha,j}| \le
|j| <\mu$.  If for some $j$ the set $\cU_{\alpha,j}$ is unbounded in $\alpha$
let $j(\alpha)$ be the minimal such $j$ and $u_\alpha =
\cU_{\alpha,j(\alpha)}$.

If for every $j,\cU_{\alpha,j}$ is bounded in $\alpha$ let $u_\alpha =
\{\sup(\cU_{\alpha,j}):j < \mu\}$, so easily $\otp(u_\alpha) \le \mu$.  
So $\langle u_\alpha:\alpha
<\lambda\rangle$ witness $\lambda \in \cf-\id_{<(\mu +1)}(\lambda)$,
i.e. the second condition holds.
\bigskip

\noindent
\underline{Second condition implies third condition}:

Trivial.
\bigskip

\noindent
\underline{Third condition implies first condition}:

Let $\gamma < \lambda$ and let $\bar u = \langle u_\alpha:\alpha <
\lambda\rangle$ witness $\lambda \in \cf-\id_{< \gamma}(\langle
\alpha:\alpha < \lambda\rangle)$; let $f_*:\gamma \rightarrow \mu$ be
one-to-one.  Defined a one-to-one function $f_\alpha:\alpha \rightarrow
\mu$ by induction on $\alpha \in [\mu,\lambda)$, the induction step as
 in the proof of \ref{p21}(1).

\noindent
2) Lastly, assume $\mu$ is singular; obviously the fourth condition
implies the third.
\bigskip

\noindent
\underline{Second condition implies the fourth condition}:

Let $\langle u_\alpha:\alpha < \lambda\rangle$ witness $\lambda \in
\cf-\id_{<(\mu+1)}(\langle \alpha:\alpha < \lambda\rangle)$, let
$f_\alpha$ be the unique order preserving function from $u_\alpha$
onto $\otp(u_\alpha)$.  Let $u \subseteq \mu = \sup(u)$ has order type
$\cf(\mu)$ or just $< \mu$.
Let $u'_\alpha$ be $u_\alpha$ if $\otp(u_\alpha) < \mu$ and be $\{\beta \in
u_\alpha:f_\alpha(\beta) \in u\}$ if $\otp(u_\alpha) = u$.
\end{PROOF}

\noindent
The next claim says that quite many partial squares on $\lambda =
\mu^+$ exists.
\begin{claim}
\label{p36}
[$\Ax_{4,\partial}$]  Assume $\lambda$ is the true successor of
$\mu,\theta \le \kappa = \cf(\mu),\theta \le \theta_1 < \mu,
\partial < \theta$ and $\alpha < \mu \Rightarrow \hrtg({}^{\theta
  >}\alpha) < \mu$ and $\alpha < \theta = \hrtg([\alpha]^{< \partial})
< \theta_1$.  

\Then \, we can find $\bar C = \langle C_{\varepsilon,\alpha}:\varepsilon
< \mu,\alpha \in S_\varepsilon\rangle$ such that:
\mn
\begin{enumerate}
\item[$(a)$]  $S_\varepsilon \subseteq S^\lambda_{< \theta_1} :=
  \{\delta < \lambda:\cf(\delta) < \theta_1\}$
\sn
\item[$(b)$]  $S^\lambda_{< \theta} \subseteq \cup
  \{S_\varepsilon:\varepsilon < \mu\}$
\sn
\item[$(c)$]  $C_{\varepsilon,\alpha} \subseteq \alpha$ and
  $C_{\varepsilon,\alpha}$ is closed unbounded in $\alpha$
\sn
\item[$(d)$]  $\beta \in C_{\varepsilon,\alpha} \Rightarrow
  C_{\varepsilon,\beta} = C_{\varepsilon,\alpha} \cap \beta$
\sn
\item[$(e)$]  $\otp(C_{\varepsilon,\alpha}) < \theta_1$.
\end{enumerate}
\end{claim}

\begin{PROOF}{\ref{p36}}
Let $X \subseteq \lambda$ code:
\mn
\begin{enumerate}
\item[$\bullet$]   a witness to ``$\lambda$ is the true successor 
of $\mu$"
\sn
\item[$\bullet$]  the set $S^*_0 := S^\lambda_{< \theta},S^*_1
= S^\lambda_{<\theta_1}$
\sn
\item[$\bullet$]  a witness to $\cf(\mu) = \kappa$
\sn
\item[$\bullet$]  $\langle e_\alpha:e < \lambda\rangle$ as in 
$(*)_4$ of the proof of \ref{c13} so $\alpha \in S^*_0 \Rightarrow
  |e_\alpha| < \theta_1$.
\end{enumerate}
\mn
So $\mathbf L[X] \models ``\lambda = \mu^+,\cf(\mu) = \kappa \ge
\theta"$ and $\chi < \mu \Rightarrow \chi^{< \theta} < \mu$.  If
$\mathbf L[X] \models ``\mu$ is regular", by \cite[\S4]{Sh:351} and if
$\mathbf L[X] \models ``\mu$ is singular" by 
Dzamonja-Shelah \cite{Sh:562} we get the result in $\mathbf L[X]$ and
the same $\bar C$ works in $\mathbf V$.
\end{PROOF}

\noindent
For more on successor, see \cite[\S(3A)]{Sh:955} and in
\cite[0x=Ls3]{Sh:F1303}.
\bigskip

\subsection {Covering number} \label{2D} \
\bigskip

\begin{definition}
\label{s2}
1) Let $\cov(\lambda,\theta,\le Y,\sigma)$ be the minimal cardinal
   $\chi$ such that (if no such $\chi$ exists, it is $\infty$ (or not
   well defined)): there is a set $\cP$ of cardinality $\chi$ such
   that:
\mn
\begin{enumerate}
\item[$(a)$]  $\cP \subseteq [\lambda]^{< \theta}$
\sn
\item[$(b)$]  if $f \in {}^Y \lambda$ \then \, there is $\cP'
  \subseteq \cP$ of cardinality $< \sigma$ such that $\Rang(f)
  \subseteq \cup \{u:u \in \cP'\}$.
\end{enumerate}
\mn
1A) Writing $\kappa$ instead ``$\le Y$" means $f \in
\bigcup\limits_{\alpha < \kappa} {}^\alpha \lambda$.

\noindent
2) If $\sigma = 2$ we may omit it.

\noindent
3)  Writing ``$\le \theta$" instead of $\theta$ means $\theta^+$,
i.e.  $\cP \in [\lambda]^{\le \theta}$.
\end{definition}

\begin{definition}
\label{s5}
1) We say $([\gamma]^\theta,\subseteq)$ strongly\footnote{without
 ``strongly" we have only $f_\alpha:\gamma_\alpha \rightarrow
\mu$ where $\gamma_\alpha < \theta^+$} has cofinality $\le \chi$
\when \, there is $\bar f = \langle f_\alpha:\alpha < \alpha_*\rangle$
such that $|\alpha_*| = \chi$ and $f_\alpha:\theta \rightarrow \mu$
and for every $u \in [\gamma]^\theta$ there is $\alpha$ such that $u
\subseteq \Rang(f_\alpha)$.

\noindent
2) We replace $``\le \chi$" by $``\chi"$ \when \, in addition
   $([\gamma]^\theta,\subseteq)$ has cofinality $\chi$.
\end{definition}

\begin{claim}
\label{s8}
If $([\gamma]^\theta,\subseteq)$ has cofinality $\chi$ and
$\theta^+$ is a truely successor \then \,
$([\gamma]^\theta,\subseteq)$ strongly has cofinality $\chi$.
\end{claim}

\begin{PROOF}{\ref{s8}}
Easy.
\end{PROOF}

\begin{theorem}
\label{s14}
Assume $\Ax_{4,\partial},\partial < \theta_*,\langle \theta_Y = \theta(Y):Y \in
\theta_*\rangle$ is such that $(\theta_Y,Y)$ satisfies the
demands on $(\theta,Y)$ in \ref{c13} and $\theta_Y < \theta_*$ and so
$\theta_*$ is strong limit in the sense that $Y \in \theta_* \Rightarrow
\hrtg(\Fil^4_{\aleph_1}(Y)) < \theta_*$, equivalently $\kappa <
\theta_* \Rightarrow \hrtg(\cP(\cP(\kappa)) < \theta_*$
 (and $\theta_* > \partial$; see \ref{z24}).

\noindent
1) For all cardinals $\lambda \ge \theta_*$ we have 
$\cov(\lambda,\le \theta_*,< \theta_*,2)$ is well defined (i.e. $<
\infty$).

\noindent
2) Even $\partial$-uniformly and in some inner model $\mathbf L[X],X
   \subseteq \Ord$ we have witness for those covering numbers. 
\end{theorem}

\begin{PROOF}{\ref{s14}}
Let $\lambda_* = \cup\{\hrtg({}^\kappa \lambda):\kappa < \theta_*\}$
\mn
\begin{enumerate}
\item[$\boxplus_1$]  $(a) \quad$ let $(\cS_*,<_*)$ be such that $\cS_*
  \subseteq [\lambda_*]^{< \partial}$ satisfy

\hskip25pt $(\forall u \in [\lambda_*]^{\aleph_0})
(\exists v \in \cS_*)[u \subseteq v]$ and
 $<_*$ is a well ordering of $\cS_*$
\sn
\item[${{}}$]  $(b) \quad$ we define $c \ell$ and 
$\cS_{\lambda_*,\kappa} \subseteq
[\lambda_*]^{< \partial},<_{\lambda_*,\kappa},\langle w^*_{\kappa,i},i
< \otp(\cS_{\lambda_*},<_*)\rangle,\Omega_\kappa,\bar e_\kappa$ 

\hskip25pt as in $(*)_1 - (*)_4$ in the proof of \ref{c13} with $\kappa$
  here standing for $Y$ 

\hskip25pt there, from $(\cS_*,<_*)$.
\end{enumerate}
\mn
So we can choose $\bar F = \langle F^1_\kappa:\kappa < \theta_*\rangle$ where
\mn
\begin{enumerate}
\item[$\boxplus_2$]  $(a) \quad F^1_\kappa$ is a function
\sn
\item[${{}}$]  $(b) \quad \Dom(F^1_\kappa) = \{f:f \in 
{}^\kappa(\lambda +1)$ and
  $i < \kappa \Rightarrow \cf(f(i) \ge \theta_\kappa\}$
\sn
\item[${{}}$]  $(c) \quad F^1_\kappa(f)$ is a pair
  $(\cF^1_f,<^1_f)$ such that
\sn
\begin{enumerate}
\item[${{}}$]  $(\alpha) \quad \cF^1_f \subseteq \prod\limits_{i <
  \kappa} f(i)$ is cofinal, i.e. modulo the filter $\{\kappa\}$
\sn
\item[${{}}$]  $(\beta) \quad <^1_f$ is a well ordering of $\cF^1_f$.
\end{enumerate}
\end{enumerate}
\mn
[Why possible?  By \ref{d5} and \ref{c24}(2).]

Let $(\theta_{n+1}(\kappa))$ exist and is $< \theta_*$, see
\cite[0.14]{Sh:835} where
\mn
\begin{enumerate}
\item[$\boxplus_3$]  for $\kappa < \theta_*$, let 
$\theta_0(\kappa) = \theta_\kappa$ and
  $\theta_{n+1}(\kappa) := \min\{\sigma$: if $\langle u_i:i <
  \kappa\rangle$ is a sequence of sets of ordinals each of cardinality
  $< \theta_n(\kappa)$ then $\sigma > |\bigcup\limits_{i < \kappa} u_i|\}$.
\end{enumerate}
\mn
Choose $\langle (\cF^2_{\kappa,n},<^2_{\kappa,n}):\kappa <
\theta_*\rangle$ by induction on $n$, so
$\langle (\cF^2_{\kappa,n},<^2_{\kappa,n}):n < \omega$ and ordinal $\kappa <
\theta_*\rangle$ exists, such that:
\mn
\begin{enumerate}
\item[$\boxplus_4$]  $(a) \quad$ if $n=0$ then $\cF^2_{\kappa,n} =
  \{f^2_*\},f^2_* \in {}^\kappa(\lambda +1)$ is constantly $\lambda$
\sn
\item[${{}}$]  $(b) \quad$ if $f \in \cF^2_{\kappa,n}$ then $f$ is a
  function from $\kappa$ into $\{u \subseteq \lambda +1:|u| \le
\theta_n(\kappa)\}$
\sn
\item[${{}}$]  $(c) \quad <^2_{\kappa,n}$ well orders
  $\cF^2_{\kappa,n}$
\sn
\item[${{}}$]  $(d) \quad$ if $f \in \cF^2_{\kappa,n}$ then for $\ell
  < 4$ we let $g^\ell_f$ be the following function; 

\hskip25pt  its domain is $\kappa$ and for $i < \kappa$ we let:

\underline{$\ell=0$}:  $g^\ell_f(i) = \{\alpha \in f(i):\alpha = 0\}$

\underline{$\ell=1$}:  $g^\ell_f(i) = \{\alpha \in f(i):\alpha$ is a
successor ordinal$\}$

\underline{$\ell=2$}:  $g^\ell_f(i) = \{\alpha \in f(i):\alpha$ is a
limit ordinal of cofinality $< \theta_\kappa\}$

\underline{$\ell=3$}:  $g^\ell_f(i) = \{\alpha \in f(i):\cf(\alpha)
\ge \theta_\kappa\}$
\sn
\item[${{}}$]  $(d)(\alpha) \quad$ if $f_1 \in \cF^2_{\kappa,n}$ \then \, for
  some $f_2 \in \cF^2_{\kappa,n+1},f_2(i) =$

\hskip35pt $\{\beta:\beta +1 \in g^1_{f_1}(i)\}$
\sn
\item[${{}}$]  $\,\,\,\,\,\,(\beta) \quad$ if $f_1 \in \cF^2_{\kappa,n}$ 
then for some $f_2 \in \cF^2_{\kappa,n+1}$ we have $f_2(i) =
\cup\{e_{\kappa,\alpha}$:

\hskip25pt $\alpha \in g^2_{f_1}(i)$ and $\cf(\alpha) < \theta\}$, 
\sn
\item[${{}}$]  $\,\,\,\,\,\,(\gamma) \quad$ if $f_1 \in \cF^2_{\kappa,n}$ 
letting $u := \otp(\cup\{g^3_{f_1}(i):i < \kappa\})$, i.e. $\zeta =
\zeta_f =$

\hskip25pt $\otp(u) < \theta_*,\bar\delta_{f_1} = \langle 
\delta_{f_1,\iota}:\iota < \zeta \rangle$ increasing $\delta_{f_1,\iota}
\in u$ and 

\hskip25pt $\otp(\delta_{f_1,\iota} \cap u) = \iota$ then 
$F^1_{\otp(u)}(\bar\delta_{f_1}) \subseteq \cF^2$
\sn
\item[${{}}$]  $(e)(\alpha) \quad \cF^2_{\kappa,n+1}$ is minimal under the
  conditions above
\sn
\item[${{}}$]  $\,\,\,\,\,\,(\beta) \quad <^2_{\kappa,n+1}$ is chosen
  naturally.
\end{enumerate}
\mn
We can choose a set $X_2$ of ordinals such that $\langle
\cF^2_{\kappa,n}:\kappa \in \theta_*,n < \omega\rangle$ belongs 
to $\mathbf L[X_2]$ hence a list $\langle w^*_\alpha:
\alpha < \alpha_2(*)\rangle \in \mathbf L[X_2]$ of $\{\Rang(f):f 
\in \cF^2_{\kappa,n}$ for some $\kappa < \theta_*,n < \omega\}$ and a 
list $\bar u = \langle u_\alpha:\alpha <
\alpha_3(*)\rangle$ of a cofinal subset of
$[\alpha_2(*)]^{\aleph_0}$ and $X_3$ such that
$X_2,\bar u \in \mathbf L[X_3]$.

Now for any ordinal $\kappa <  \theta_*$ and $f \in {}^\kappa 
\lambda$ we can choose
finite $v_n \subseteq \alpha_2(*)$ by induction on $n$ such that:
\mn
\begin{enumerate}
\item[$(*)_n$]  $(a) \quad \lambda \in \cup\{w^*_\alpha:\alpha \in
  v_n\}$ for $n=0$
\sn
\item[${{}}$]  $(b) \quad$ if $i < \kappa,f(i) \notin
  \cup\{w^*_\alpha:\alpha \in v_n\}$ then
  $\min(\bigcup\limits_{\alpha \in v_n} w^*_\alpha \backslash f(i)) >
  \min(\bigcup\limits_{\alpha \in v_n} w^*_\alpha \backslash f(i))$.
\end{enumerate}
\mn
So $\langle v_n:n < \omega\rangle$ exists hence $v =
\bigcup\limits_{n} v_n \in \mathbf L[X_3]$, hence $w =
\bigcup\limits_{\alpha \in v} w^*_\alpha \in \mathbf L[X_3]$ has
cardinality $\le \theta_*$ and includes $\Rang(f)$ because if $i <
\kappa \wedge f(i) \notin \cup\{w^*_\alpha:\alpha \in v\}$ then
$\langle \min(\bigcup\limits_{\alpha \in v_n} w^*_\alpha \backslash
f(i)):n < \omega\rangle$ is a strictly decreasing sequence of ordinals.
So we should just
let $\cP = \{u \subseteq \lambda:u \in \mathbf L[X_3]$ and $\mathbf L[X_3]
\models ``|u| \le \theta_*\}$ witness the desired conclusion.
\end{PROOF}

\noindent
Now (like \cite[\S(3A)]{Sh:955} see definitions there)
\begin{conclusion}
\label{s21}
Assume $\Ax_4$.  If $\mu$ is a singular cardinal such that $\kappa <
\mu \Rightarrow \theta_\kappa := \hrtg(\cP(\cP(\kappa))^+ < \mu$ and
$\lambda \le \kappa$ \then \, for some $\kappa < \mu$ we have:
$\cov(\lambda,\mu,\mu,\kappa) = \lambda$. 
\end{conclusion}

\begin{PROOF}{\ref{s21}}
Use \cite{Sh:460} in $\mathbf L[X]$ where $X \subseteq \Ord$ is as in
\ref{s14}(2). 
\end{PROOF}

\begin{discussion}
\label{s24}
0) From \ref{s14}, \ref{s21} we can get also smooth closed generating
sequence (see \cite[\S6]{Sh:430}, \cite{Sh:E69} (an earlier version is 
\cite{Sh:E29}).

\noindent
1) We would like to get better bounds.  A natural way is to fix
   $\kappa$, consider $\theta_1 > \kappa$ and $\mathbf f:\kappa
   \rightarrow [\lambda]^{< \theta_1}$ and ask for $\cF \subseteq
   \{f:\kappa \rightarrow [\lambda]^{< \theta_2}\}$ such that for
   every $g \in \prod\limits_{i < \kappa} (\mathbf f(i) \cup \{1\}
   \backslash \{0\})$ and $g_i \in \prod\limits_{i < \kappa} g_*(i)$
   there is $f \in \cF$ such that $(\forall i < \kappa)(f(i) \cap
   [g_1(i),g_*(i)) \ne \emptyset)$.

\noindent
2) We can get also strong covering, see \cite[Ch.VII]{Sh:g}.

\noindent
3) Can we get something better on $\mu$ singular strong limit?  a BB?,
(BB means black box, see \cite{Sh:309} and in \S3, possibly see more in 
\cite{Sh:F1200}.

\noindent
4) We like to improve \ref{s14}, in particular \S(2C), for this we have
to improve \S(2A). We would like to replace $\Fil^4_{\aleph_1}(Y)$,
i.e. $\hrtg(\Fil^4_{\aleph_1}(Y))$ by $\hrtg(\cP(Y))$ and even
$\hrtg(Y)$, as done in ZFC in \cite{Sh:345}.  We do not know to do
this but we try a more modest aim: suppose we deal only with $[Y]^{\le
  \kappa}$ or so.  So hopefully in \cite{Sh:F1303}, we still have
$\hrtg(\Fil^4_{\aleph_1}(\kappa))$ but $\hrtg(\cP(Y))$ only. 
\end{discussion}
\newpage

\section {Black Boxes} \label{3}

There are many proofs in ZFC using diagonalization of various kinds so
they seem to depend heavily on choice.  Using $\Ax_4$ we succeed to
generalize one such method - one of the black boxes from
\cite{Sh:309}, it seems particularly helpful in constructing abelian
groups and modules; see on applications in the books 
Eklof-Mekler \cite{EM02} and G\"obel-Trlifaj \cite{GbTl12}.

The proof specifically uses countable models and $\Ax_4$.  Naturally we
would like to assume we have only $\Ax_{4,\partial}$.  But existing
versions implies $\cP(\bbN)$ is well ordered and more, whereas
$\Ax_{4,\partial}$ does not imply this.
\bigskip

\subsection {Existence proof} \label{3A} \
\bigskip

\begin{hypothesis}
\label{gy6}
$\ZF + \DC + \Ax_4$ [so $\partial = \aleph_1$]

The following is like \cite[3.24(3)]{Sh:309}, the relevant cardinals
provably exists but may be less common than there: 
conceivably true successor are only successor of singular 
strong limit cardinals.
\end{hypothesis}

\begin{theorem}
\label{g2}
If (A) then (B) \underline{where}:
\mn
\begin{enumerate}
\item[$(A)$]  $(a) \quad \lambda = \mu^+$ is a true successor
\sn
\item[${{}}$]  $(b) \quad \mu = \mu^{\aleph_0}$
\sn
\item[${{}}$]  $(c) \quad S = \{\delta < \lambda:\cf(\delta) =
  \aleph_0$ and $\mu$ divides $\delta\}$ or just $S$ is a stationary 

\hskip25pt subset of $\lambda$ 
such that $\delta \in S \Rightarrow \cf(\delta) =
\aleph_0 \wedge \mu < \delta \wedge (\mu|\delta)$
\sn
\item[${{}}$]  $(d) \quad \bar\gamma^* = \langle
  \bar\gamma^*_\delta:\delta \in S\rangle$ with $\bar\gamma^*_\delta =
  \langle \gamma^*_{\delta,n}:n < \omega\rangle$ an increasing
  $\omega$-sequence

\hskip25pt  of ordinals with limit $\delta$
\sn
\item[$(B)$]  we can find $\mathbf w = (\alpha,\mathbf W,
\dot\zeta,h,\bar{\mathbf k}) 
= (\alpha_{\mathbf w},\mathbf W_{\mathbf w},\dot\zeta_{\mathbf w},h_{\mathbf
  w},\bar{\mathbf k}_{\mathbf w})$ such that (we may denote 
$\alpha_{\mathbf w}$ by $\ell g(\mathbf w)$ and may omit it):
\sn
\begin{enumerate}
\item[$(a)$]  $(\alpha) \quad \mathbf W = \langle \bar N_\alpha:
\alpha < \alpha_{\mathbf w}\rangle$
\sn
\item[${{}}$]  $(\beta) \quad \bar N_\alpha = \langle N_{\alpha,n}:n <
  \omega\rangle$ is $\prec$-increasing sequence of models
\sn
\item[${{}}$]  $(\gamma) \quad \tau(N_{\alpha,n}) \subseteq
  \cH(\aleph_0)$ and $\tau(N_{\alpha,n}) \subseteq \tau(N_{\alpha,n+1})$
\sn
\item[${{}}$]  $(\delta) \quad \mathbf k = \langle \bar k_\alpha:
\alpha < \alpha_{\mathbf w}\rangle,\bar k_\alpha = \langle
 k_{\alpha,n}:n < \omega\rangle$ is increasing, 

\hskip25pt let $k_{\mathbf w}(\alpha,n) = k(\alpha,n) = k_{\alpha,n}$
\sn
\item[${{}}$]  $(\varepsilon) \quad |N_{\alpha,n}| \subseteq
  |N_{\alpha,n+1}| \subseteq \lambda$
 but $N_{\alpha,n} \ne N_{\alpha,n+1}$
\sn
\item[${{}}$]  $(\zeta) \quad$ let $N_\alpha =N_{\alpha,\omega} =
 \lim(\bar N_\alpha)$, that is, $\tau(N_{\alpha,\omega}) =$

\hskip25pt $\cup\{\tau(N_{\alpha,n}):n <\omega\}$
  and $(N_{\alpha,\omega} \rest \tau(N_{\alpha,n})) \supseteq N_{\alpha,n}$
\sn
\item[${{}}$]  $(\eta) \quad$ the universe of $N_{\alpha,n}$ is a
 countable subset of $\lambda$
\sn
\item[$(b)$]  $(\alpha) \quad \dot\zeta$ is a function from
  $\alpha_{\mathbf w}$ into $S$, non-decreasing
\sn
\item[${{}}$]  $(\beta) \quad$ if $\dot\zeta(\alpha) = \delta$ then
  $\delta = \sup\{\gamma^*_{\delta,n}:n < \omega\} = \sup(N_\alpha)$
\sn
\item[${{}}$]  $(\gamma) \quad$ if $\alpha < \alpha_{\mathbf w}$ and
  $\dot\zeta(\alpha) = \delta \in S$ and $n < \omega$ then
  $N_{\alpha,n+1} \backslash N_{\alpha,n}$

\hskip25pt $\subseteq (\gamma^*_{\delta,k(\alpha,n)},
\gamma^*_{\delta,k(\alpha,n)+1})$ and $|N_{\alpha,n}| \subseteq
\gamma^*_{\delta,k(\alpha,n)}$ 
\sn
\item[$(c)$]  if $M$ is a model with universe $\lambda$ and vocabulary
  $\subseteq \cH(\aleph_0)$ \then \, for stationarily many $\delta \in S$,
  there is $\alpha$ such that $\dot\zeta(\alpha) = \delta,N_\alpha \prec M$.
\sn
\item[$(d)$]   $(\alpha) \quad$ if $\dot\zeta(\alpha) = \delta =
  \dot\zeta(\beta)$ then $N_\alpha \cong N_\beta,|N_\alpha| \cap \mu =
  |N_\beta| \cap \mu,\bar k_\alpha = \bar k_\beta$; 

\hskip25pt moreover, $\otp(|N_\alpha|) = 
\otp(|N_\beta|)$ and the unique order preserving

\hskip25pt  mapping is an isomorphism from $N_{\alpha,n}$ onto $N_{\beta,n}$
  for every $n$ 

\hskip25pt and is the identity on $|N_\alpha| \cap \mu$ and on
$N_\alpha \cap N_\beta$ and so maps

\hskip25pt $N_\alpha \cap \gamma^*_{\delta,k(\alpha,n)}$ onto $N_\beta
\cap \gamma^*_{\delta,k(\beta,n)}$
\sn
\item[${{}}$]  $(\beta) \quad$ if $\dot\zeta(\alpha) = \delta =
 \dot\zeta(\beta)$ but $\alpha \ne \beta$ then

\hskip35pt $\bullet \quad N_\alpha \cap N_\beta$ is an initial segment
 of both $N_\alpha$ and of $N_\beta$ 

\hskip35pt $\bullet \quad N_\alpha \cap N_\beta \subseteq
 N_{\alpha,n+1} \cap N_{\beta,n+1}$ and $N_\alpha \cap N_\beta
 \supseteq N_{\alpha,n} = N_{\beta,n}$

\hskip45pt  for some $n$.
\end{enumerate}
\end{enumerate}
\end{theorem}

\begin{remark}
\label{g5}
1) The existence proof is uniform (that is, $\mathbf w$ can be defined
from $(<_*,\bar f)$ where: $<_*$ is a well ordering of
$[\chi]^{\aleph_0}$ for $\chi$ large enough and $\bar f$ is a witness for
$\lambda$ being a true successor.  Moreover, also $\bar\gamma^*$ can
be chosen uniformly (as well as the witness for 
$\lambda$-being a true successor.

\noindent
2)  We would like to add (A)(e) to the assumption 
and add (B)(e) to the conclusion of \ref{g2} where:
\mn
\begin{enumerate}
\item[$(A)(e)$]  $(\alpha) \quad \bar C 
= \langle C_\delta:\delta \in S\rangle$
\sn
\item[${{}}$]  $(\beta) \quad C_\delta \subseteq \delta =\sup(C_\delta)$
\sn
\item[${{}}$]  $(\gamma) \quad \otp(C_\delta) = \omega$ and let
  $\bar\gamma^*_\delta = \langle \gamma^*_{\delta,n}:n <
  \omega\rangle$ list $C_\delta$ in increasing order
\sn
\item[${{}}$]  $(\delta) \quad \bar C$ weakly guess clubs, i.e. for
  every club $E$ of $\lambda$ for stationarily many 

\hskip25pt $\delta \in S$ we
  have $(\forall n)(E \cap (\gamma^*_{\delta,n},\gamma^*_{\delta,n+1})
  \ne \emptyset)$, moreover
\sn
\item[${{}}$]  $(\varepsilon) \quad \langle S_\varepsilon:\varepsilon <
  \lambda\rangle$ is a partition of $S$ such that $\bar C \rest
S_\varepsilon$ weakly guess clubs 

\hskip25pt for each $\varepsilon$
\sn
\item[$(B)$\,\,\,]  $(e) \quad N_{\alpha,n+1} \backslash N_{\alpha,n}$ is
included in $[\gamma^*_{\delta,n},\gamma^*_{\delta,n+1})$, that is
 $k_{\mathbf w}(\alpha,n) = n$.
\end{enumerate}
\mn
But not clear if (A) is provable in our context.  Still, repeating the ZFC 
proof works in $\ZF + \DC_{\aleph_1}$ and gives even ``$\bar C$ guess
clubs", i.e. ``$\{\gamma_{\delta,n}:n < \omega\} \subseteq
C_\delta$".  But we ask only 
for ``weakly guess", see \ref{g5}(2), $(A)(e)(\delta)$ so using $\Ax_4$ just
adding $\AC_{\cP(\bbN)}$ suffice\footnote{That is, having $\bar S =
  \langle S_\varepsilon:\varepsilon < \mu\rangle$ for each
  $\varepsilon$ choose the first increasing function $f \in {}^\omega
  \omega$ such that $\langle \gamma^*_{\delta,f(n)}:\delta \in
  S_\varepsilon\rangle$ weakly guess clubs.}.
However, clause $(B)(d)(\beta)$ is a reasonable substitute.

\noindent
2) We may strengthen clause (B)(d) by adding:
\mn
\begin{enumerate}
\item[$(\gamma)$]  if $\dot\zeta(\alpha) = \delta = \zeta(\beta)$ then
  $|N_\alpha| \cap \gamma(\delta,0) = |N_\beta| \cap \gamma(\delta,0)$
  call it $u_\delta$.
\end{enumerate}
\mn
For this in $(*)_6$ the partition should be $\langle
S_\varepsilon:\varepsilon < \lambda\rangle$ as $\varepsilon$ should
determine also $N_\delta$, etc.

\noindent
3) The use of $\kappa$ possibly $> \aleph_1$ in \ref{g17} is not
   necessary for \ref{g2}.

\noindent
4) Note that in proof we need $\mu = \mu^{\aleph_0}$ for proving $(*)_3$.
Note that for $(*)_6(a),(b),(c)$ we need just ``$\lambda$ is a true
successor of $\mu$".  To get clause (d) too, it suffices to have
   $\mu = \mu^{\aleph_0}$.

\noindent
5) We may prove also \ref{g23} inside the proof of \ref{g2}.
\end{remark}

\begin{PROOF}{\ref{g2}}
Now
\mn
\begin{enumerate}
\item[$\boxplus_1$]  there are $g^0,g^1$ such that
\sn
\begin{enumerate}
\item[$(a)$]  $g^0,g^1$ are two-place functions from $\lambda$ to
  $\lambda$ which are zero on $\mu$
\sn
\item[$(b)$]  $(\alpha) \quad$ if $\alpha \in [\mu,\lambda)$ then $\langle
  g^0(\alpha,i):i < \mu\rangle$ enumerate 

\hskip25pt $\{j:j < \alpha\}$ without repetitions
\sn
\item[${{}}$]  $(\beta) \quad$ if $\alpha,i < \lambda$ and $\alpha <
  \mu \vee i \ge \mu$ then $g^0(\alpha,i) = 0$
\sn
\item[$(c)$]  $(\alpha) \quad g^1(\alpha,g^0(\alpha,i)) = i$ \when \, 
$i < \mu \le \alpha < \lambda$
\sn 
\item[${{}}$]  $(\beta) \quad$ if $\alpha < \mu$ and $i < \lambda$
 then $g^1(\alpha,i) = 0$
\sn
\item[${{}}$]  $(\gamma) \quad$ if $\alpha \le i < \lambda$
then $g^1(\alpha,i) = 0$
\sn
\item[$(d)$]  there is $\gamma_* \in (\mu,\lambda)$ such that
for every countable $u \subseteq \lambda$ closed under 
$g^0,g^1$ there is $v$ such that:
\sn
\item[${{}}$]  $(\alpha) \quad v \subseteq \gamma_*$ is countable
\sn
\item[${{}}$]  $(\beta) \quad \otp(v) = \otp(u)$
\sn
\item[${{}}$]  $(\gamma) \quad v \cap \mu = u \cap \mu$
\sn
\item[${{}}$]   $(\delta) \quad v$ is closed under $g^0,g^1$
\sn
\item[${{}}$]   $(\varepsilon) \quad$ the (unique) order preserving
  function from $u$ onto $v$ commute 

\hskip25pt with $g^0,g^1$
\sn
\item[${{}}$]  $(\zeta) \quad$ we can arrange that $\gamma_* = \mu + \mu$.
\end{enumerate}
\end{enumerate}
\mn
[Why?  As $\lambda$ is truely successor there is no problem to choose
$g^0,g^1$ satisfying clauses (a),(b),(c).  On 
$\cU = \{u \subseteq \mu^+:u$ countable closed under
  $g^0,g^1\}$ we define an equivalence relation $E$ by
$(d)(\beta),(\gamma),(\varepsilon)$.   Now as
  $\mu = \mu^{\aleph_0},\cU/E$ has cardinality $\mu$ hence recalling 
$\lambda$ is regular we can prove that there is
$\gamma_*$ as required in $(d)(\alpha)-(\varepsilon)$ exists.  In
  fact, $\partial$-uniformly we have a well ordering $<_{\cU}$ of
  $\cU$; \wilog \, $u_1 <_{\cU} u_2 \Rightarrow \sup(u_1) \le
  \sup(u_2)$.

To have $\gamma_* = \mu + \mu$, let $\tau_*$ be the vocabulary
$\{F_0,F_1\}$ with $F_0,F_2$ binary function and let $\mathbf M = \{M:M$
is a $\tau_*$-model with universe $|M|$ a countable subset of $\mu +
\mu$ such that $\alpha,\beta \in M \cap \mu \Rightarrow
F_0(\alpha,\beta) = 0 = F_1(\alpha,\beta)$ and the functions
$F^M_0,F^M_1$ satisfies the relevant cases of the demands $(a),(b),(c)$
on $(g^0,g^1)\}$.

Clearly $\mathbf M$ has cardinality $\mu$ and moreover we can (uniformly)
define a list $\langle M_\varepsilon:\varepsilon < \mu\rangle$ of
$\mathbf M$.

Let $i_\varepsilon = \otp(|M_\varepsilon| \backslash \mu)$ and by
induction on $\varepsilon < \mu$ we choose
$(h_\varepsilon,\gamma_\varepsilon)$ such that:
\mn
\begin{enumerate}
\item[$\boxplus_{1.2}$]  $(a) \quad \gamma_0 = \mu$
\sn
\item[${{}}$]  $(b) \quad \langle \gamma_\zeta:\zeta \le
  \varepsilon\rangle$ is increasing continuous
\sn
\item[${{}}$]  $(c) \quad h_\varepsilon$ is an order preserving
  function from $|M_\varepsilon| \backslash \mu$ onto
  $[\gamma_\varepsilon,\gamma_{\varepsilon +1})$.
\end{enumerate}
\mn
Next let $N_\varepsilon \in \mathbf M$ be such that $h_\varepsilon \cup
\id_{|M_\varepsilon| \cap \mu}$ is an isomorphism from $M_\varepsilon$
onto $N_\varepsilon$. 

Now we define the two-place function $g^*_0,g^*_1$ from $\lambda$ to
$\lambda$ as follows
\mn
\begin{enumerate}
\item[$\boxplus_{1.3}$]  $(a) \quad$ if $\varepsilon < \mu$ and
$\gamma_\varepsilon \le \alpha < \gamma_{\varepsilon +1}$ then
\sn
\begin{enumerate}
\item[$\bullet$]  if $i \in N_\varepsilon \cap \mu$ then
  $g^*_0(\alpha,i) = F^{N_\varepsilon}_0(\alpha,i)$
\sn
\item[$\bullet$]  $\langle g^*_0(\alpha,i):i \in \mu \backslash
  N_\varepsilon\rangle$ lists $\alpha \backslash N_\varepsilon$ without
repetition and is derived from $\langle g^0(\alpha,i):i < \mu\rangle$
  and $N_\varepsilon$ as in the proof of the Cantor-Bendixon theorem
  (that $|A| \le |B| \wedge |B| \le |A| \Rightarrow |A| = |B|$):
\end{enumerate}
\sn
\item[${{}}$]  $(b) \quad$ if $\alpha \in [\mu + \mu,\lambda)$ then $i
  < \mu \Rightarrow g^*_0(\alpha,i) = g^0(\alpha,i)$
\sn
\item[${{}}$]  $(c) \quad$ if $\alpha \in [\mu,\lambda)$ and $j <
  \alpha$ then $g^*_1(\alpha,j)$ is defined as the unique $i < \mu$

\hskip25pt   such that $g^*_0(\alpha,i) = j$
\sn
\item[${{}}$]  $(d) \quad$ in all other cases the value is zero.
\end{enumerate}
\mn
Now $g^*_0,g^*_1$ are well defined, just recall $\boxplus_1(a),(b),(c)$.
So $\boxplus_1$ holds indeed.]

Clearly
\mn
\begin{enumerate}
\item[$(*)_1$]  if $u_1,u_2 \subseteq \lambda$ are closed under
  $g^0,g^1$ and $u_1 \cap \mu = \mu_2 \cap \mu$ \then \, $u_1 \cap
  u_2$ is an initial segment of $u_1$ and of $u_2$.
\end{enumerate}
\mn
Let $\mathbf N$ be the set of tuples $(\bar N,\bar\gamma)$ satisfying
\mn
\begin{enumerate}
\item[$(*)_2$]  $(a) \quad \bar N = \langle N_n:n < \omega\rangle$
\sn
\item[${{}}$]  $(b) \quad N_n$ is a model with vocabulary $\tau(N_n)
  \subseteq \cH(\aleph_0)$
\sn
\item[${{}}$]  $(c) \quad N := \cup\{N_n:n < \omega\}$ 
is countable with universe $\subseteq \gamma_*$
\sn
\item[${{}}$]  $(d) \quad \tau(N_n) \subseteq \tau(N_{n+1})$ with $N_n
  \subseteq N_{n+1} \rest \tau_n$
\sn
\item[${{}}$]  $(e) \quad \bar\gamma = \langle \gamma_n:n <
  \omega\rangle$ is an increasing sequence of ordinals satisfying 

\hskip25pt $\cup\{\gamma_n:n < \omega\} = \cup\{\alpha +1:\alpha \in
\cup\{N_n:n <\omega\}\} < \gamma_*$
\sn
\item[${{}}$]  $(f) \quad N_n = (N_{n+1} \rest \tau(N_n)) \rest
  \gamma_n$
\sn
\item[${{}}$]  $(g) \quad \sup(N_n) < \gamma_n = \min(N_{n+1}
  \backslash N_n)$
\sn
\item[${{}}$]  $(h) \quad N_n$ is closed under $g_0,g_1$.
\end{enumerate}
\mn
Recalling $\cH_{< \aleph_1}(\gamma) = \{u:u$ a countable set such that
$u \cap \Ord \subseteq \gamma$ and $y \in u \backslash \gamma
\Rightarrow |y| < \aleph_1$. 
Clearly $\mathbf N \subseteq \cH_{< \aleph_1}(\gamma_*)$ so as 
$\mu^{\aleph_0} = \mu = |\gamma_*|$, clearly 
$\mathbf N$ is well orderable so (and
using parameter witnessing, $\Ax^4_\lambda + ``\lambda$ is a
true successor cardinal" to uniformize) let
\mn
\begin{enumerate}
\item[$(*)_3$]  $(a) \quad \langle (\bar
  N_\varepsilon,\bar\gamma_\varepsilon):\varepsilon < \mu\rangle$ list
  $\mathbf N$
\sn
\item[${{}}$]  $(b) \quad \bar N_\varepsilon = \langle
  N_{\varepsilon,n}:n < \omega\rangle,\bar\gamma_\varepsilon = \langle
  \gamma_{\varepsilon,n}:n < \omega\rangle$
\sn
\item[${{}}$]  $(c) \quad N_\varepsilon = N_{\varepsilon,\omega} :=
  \cup\{N_{\varepsilon,n}:n < \omega\}$, i.e. $N_\varepsilon = 
\lim(\bar N_\varepsilon)$.
\end{enumerate}
\mn
Next
\mn
\begin{enumerate}
\item[$(*)_4$]  for each $\varepsilon < \mu$ let $\mathbf N_\varepsilon$
be the set of pairs $(\bar N,\bar\gamma)$ such that:
\sn
\begin{enumerate}
\item[${{}}$]  $(a) \quad \bar N = \langle N_n:n < \omega\rangle$
\sn
\item[${{}}$]  $(b) \quad N = \cup\{N_n:n < \omega\}$ 
is a $\tau(N_\varepsilon)$-model
\sn
\item[${{}}$]  $(c) \quad N_n$ is a $\tau(N_{\varepsilon,n})$-model
  with universe $\subseteq \lambda$
\sn
\item[${{}}$]  $(d) \quad$ there is $h$, an order preserving function from
$N_{\varepsilon,\omega}$ onto $N$

\hskip25pt  commuting with $g^0,g^1$ mapping $N_{\varepsilon,n}$ onto
$N_n$, 

\hskip25pt (i.e. $h \rest N_{\varepsilon,n}$ is an isomorphism from
$N_{\varepsilon,n}$ onto $N_n$) 

\hskip25pt  and being the identity on $N_\varepsilon \cap \mu$ and so
mapping $\gamma_{\varepsilon,n}$

\hskip25pt to $\gamma_n$
\end{enumerate}
\sn
\item[$(*)_5$]  for $\delta \in S$ and $\varepsilon < \mu$ let $\mathbf
  N_{\varepsilon,\delta}$ be the set of pairs $(N,\bar\gamma) \in
  \mathbf N_\varepsilon$ such that $\sup\{\gamma_n:n < \omega\} =
  \delta$ and for clause $(B)(b)(\gamma)$ for every $n$ for some
  $k,N_{n+1} \backslash N_n \subseteq
  (\gamma^*_{\delta,k},\gamma^*_{\delta,k+1})$
\sn
\item[$(*)_6$]   there is a partition $\bar S = \langle
  S_\varepsilon:\varepsilon < \mu\rangle$ of $S$ to stationary sets.
\end{enumerate}
\mn
[Why?  By Larson-Shelah \cite{Sh:925}.]
\mn
\begin{enumerate}
\item[$(*)_7$]  there is $\langle \bar\gamma^*_\delta:\delta \in
  S\rangle$ such that each $\bar\gamma^*_\delta$ is an increasing
  $\omega$-sequence with limit $\delta$.
\end{enumerate}
\mn
[Why?  By $\Ax_4$.]
\mn
\begin{enumerate}
\item[$(*)_8$]  there is, (in fact as in all cases in this proof,
uniformly definable), a sequence $\langle 
(\bar N_\alpha,\bar\gamma_\alpha,u_\alpha):
\alpha <\alpha(*)\rangle$ and function
  $\dot\zeta:\alpha(*) \rightarrow S$ such that:
\sn
\item[${{}}$]  $(a) \quad \dot\zeta$ is non-decreasing
\sn
\item[${{}}$]  $(b) \quad (\bar N_\alpha,\bar\gamma_\alpha) \in 
\mathbf N_{\varepsilon,\dot\zeta(\alpha)}$ when $\dot\zeta(\alpha) \in
S_\varepsilon$, moreover
\sn
\item[${{}}$]  $(b)' \quad$ if $\varepsilon < \mu$ and $\delta \in
  S_\varepsilon$ then $\{(\bar N_\alpha,\bar\gamma_\alpha):\alpha <
  \alpha(*)$ satisfies $\dot\zeta(\alpha) = \delta\}$ 

\hskip25pt list $\mathbf N_{\varepsilon,\delta}$
\sn
\item[$(*)_9$]   let $N_{\alpha,\omega} =
  \cup\{N_{\alpha,n}:n < \omega\}$.
\end{enumerate}
\mn
[Why?  By $(*)_5,(*)_6$ and using a well ordering of $[\lambda]^{\aleph_0}$.]

Now ignoring clause (c), clauses of (B) should be clear.  
Lastly, clause (c) holds by the following Theorem \ref{g17}, in 
our case $\kappa = \aleph_1$.
\end{PROOF}

\begin{theorem}
\label{g17}
If (A) then (B) where:
\mn
\begin{enumerate}
\item[$(A)$]  $(a)(\alpha) \quad \lambda > \kappa$ are regular
uncountable cardinals
\sn
\item[${{}}$]  $\quad (\beta) \quad \alpha < \lambda
  \Rightarrow |\alpha|^{\aleph_0} < \lambda$
\sn
\item[${{}}$]  $(b)(\alpha) \quad$ if $\alpha <\lambda$ then
  $\cf([\lambda]^{< \kappa},\subseteq)$ is $< \lambda$
\sn
\item[${{}}$]  $\quad (\beta) \quad \mathbf U_* \subseteq [\lambda]^{<
  \kappa}$ is well orderable and cofinal (under $\subseteq$)
\sn
\item[${{}}$]  $\quad (\gamma) \quad |\mathbf U_* \cap 
[\alpha]^{< \kappa}| < \lambda$ for $\alpha < \lambda$
\sn
\item[${{}}$]  $(c) \quad M$ is a model with universe $\lambda$ and
  vocabulary $\tau,\tau$ not necessarily

\hskip25pt  well orderable
\sn
\item[${{}}$]  $(d) \quad$ if $\alpha < \kappa$ then $\lambda >
  \hrtg(\{N:N$ a $\tau$-model with universe $\alpha$; may add

\hskip25pt that some order preserving mapping is an elementary
embedding 

\hskip25pt of $N$ into $M\})$
\sn
\item[$(B)$]  there is $\bar N$, uniformly defined from witnesses to
  (A) such that:
\sn
\item[${{}}$]  $(a) \quad \bar N = \langle N_\eta:\eta \in {}^{\omega
  >}\lambda\rangle$
\sn
\item[${{}}$]  $(b) \quad \tau(N_\eta) = \tau$
\sn
\item[${{}}$]  $(c) \quad N_\eta$ has cardinality $< \kappa$ and
  $N_\eta \cap \kappa$ is an ordinal $< \kappa$
\sn
\item[${{}}$]  $(d) \quad N_\eta$ is an elementary submodel of $M$
\sn
\item[${{}}$]  $(e) \quad$ if $\nu \triangleleft \eta$ then $N_\nu$ is
  a (proper) initial segment of $N_\eta$
\sn
\item[${{}}$]  $(f) \quad$ if $n < \omega$ and $\eta,\nu \in {}^n
  \lambda$ \then \, there is an order preserving function 

\hskip25pt from $N_\eta$ onto $N_\nu$ which is an isomorphism
\sn
\item[${{}}$]  $(g) \quad$ if $n < \omega,\eta \in {}^n \lambda$ and
$\gamma < \lambda$ \then \, there is $\nu$ such that $\eta
\triangleleft \nu \in {}^{n+1}\lambda$

\hskip25pt  and $\min(N_\nu \backslash N_\eta) > \gamma$.
\end{enumerate}
\end{theorem}

\begin{remark}
\label{g18}
1) We may consider adding: 
$N_\eta(\eta \in {}^\omega \lambda)$ has $\Sigma_1$-property
   and use: $\hrtg$(the set of expansions of $\bar N^*) < \lambda$.

\noindent
2) The ZFC version of \ref{g17} is from Rubin-Shelah \cite{Sh:117}.

\noindent
3) Note that in \ref{g17} the vocabulary is constant whereas in 
\ref{g2} it is not.  But the difference is not serious as in 
\ref{g2} the vocabulary is $\subseteq \cH(\aleph_0)$ so 
there is one vocabulary which is enough to code any other.

\noindent
4) We may continue in \cite[8.2=Lg19]{Sh:F1303}.
\end{remark}

\begin{PROOF}{\ref{g17}}
Now
\mn
\begin{enumerate}
\item[$(*)_0$]    \wilog \, $\mathbf U_* \subseteq [\lambda]^{< \kappa}$ is
  closed under countable unions and initial segments.
\end{enumerate}
\mn
[Why?  By (A)(a),(b), the point is that the closure retains the
  properties.]
\mn
\begin{enumerate}
\item[$(*)_1$]   let $\mathbf N$ be the set of $\bar N$ such that
\sn
\begin{enumerate}
\item[$(a)$]  $\bar N = \langle N_n:n < \omega\rangle$
\sn
\item[$(b)$]  $(\alpha) \quad N_n \prec M$ has cardinality $< \kappa$
\sn
\item[${{}}$]  $(\beta)\quad$ moreover, $|N_n| \in \mathbf U_*$
\sn
\item[$(c)$]  $|N_n|$ is an initial segment of $|N_{n+1}|$
\sn
\item[$(d)$]  $N_n$ has cardinality $< \kappa$ and $N_0 \cap \kappa$
  is an ordinal $< \kappa$
\sn
\item[$(e)$]  $\tau(N_n) = \tau$
\end{enumerate}
\end{enumerate}
\mn
Now
\mn
\begin{enumerate}
\item[$(*)_2$]  $\mathbf N$ is well orderable
\end{enumerate}
\mn
[Why?  Recall $\mathbf U_*$ is well orderable so let $\langle u^*_\alpha:\alpha <
  \alpha_* \rangle$ list it.  Now $N_n$ is determined by $|N_n|$ (because
$N_n \prec M$) and $|\alpha_*|^{\aleph_0}$ is well orderable so we are done.]
\mn
\begin{enumerate}
\item[$(*)_3$]   let $\langle \bar N_\alpha:\alpha < \alpha_{*} \rangle$
  list $\mathbf N$ and let $\langle u^*_\alpha:\alpha < \alpha_*\rangle$
  list $\mathbf U_*$.
\end{enumerate}
\mn
[Why exists?  By $(*)_2$ and $(A)(b)(\beta)$ of the theorem assumption.]
\mn
\begin{enumerate}
\item[$(*)_4$]   $(a) \quad$ we say $\bar N',\bar N'' \in \mathbf N$ are
  equivalent and write $\bar N' \cE \bar N''$ \when \, 

\hskip25pt  for every $n,\otp(|N'_n|) = \otp(N''_n)$ and the order
  preserving function 

\hskip25pt  from $|N'_n|$ onto $|N''_n|$ is an
isomorphism and $N'_0 = N''_0$
\sn
\item[${{}}$]  $(b) \quad$ let $\mathbf N' = \{\bar N:\bar N =
  \langle N_\ell:\ell \le n\rangle = \bar N' \rest (n+1)$ for some $\bar
  N' \in \mathbf N,n \in \bbN\}$
\sn
\item[${{}}$]  $(c) \quad$ we define the equivalence relation $\cE'$
on $\mathbf N'$ by $\bar N^1 \cE' \bar N^2$ if $\bar N^1,\bar N^2$ has

\hskip25pt  the same length and the parallel of clause (a) holds
\sn
\item[${{}}$]  $(d) \quad \cE$ and $\cE'$ have $\le \mu$ equivalence
  classes.
\end{enumerate}
\mn
[Why?  E.g. clause (d) by clause $(A)(d)$ of the theorem's
  assumption.]
\mn
\begin{enumerate}
\item[$(*)_5$]  $E_1$ is a club of $\lambda$ where $E_1 := \{\delta <
  \lambda:\delta$ is a limit ordinal such that $M \rest \delta \prec
  M$ and if $\bar N \in \mathbf N$ and $\sup(N_0) < \delta$ then there
is $\bar N' \in \mathbf N$ which is $\cE$-equivalent 
to $\bar N$ with $N'_0 = N_0$
  and $\sup(\cup\{N'_n:n < \omega\}) <\delta\}$.
\end{enumerate}
\mn
[Why?  Think, noting that we can consider only $\{\bar N_\alpha:\alpha
  <\alpha_{**}$ and $\bar N_\alpha$ is not $\cE$-equivalent to $\bar
  N_\beta$ when $\beta < \alpha\}$.]
\mn
\begin{enumerate}
\item[$(*)_6$]   for $\bar N^* \in \mathbf N$ and $\bar N \in \mathbf N'$ such
  that $N_0 = N^*_0$ we define $\rk(\bar N,\bar N^*) \in \Ord \cup
\{-1,\infty\}$ by defining when $\rk(\bar N,\bar N^*) \ge \alpha$ by
  induction on the ordinal $\alpha$ as follows:
\sn
\begin{enumerate}
\item[$(a)$]  \underline{$\alpha = 0$}:  $\rk(\bar N,\bar N^*) \ge
  \alpha$ iff $\bar N \cE'(\bar N^* \rest \ell g(\bar N))$
\sn
\item[$(b)$]  \underline{$\alpha$ limit}:  $\rk(\bar N,\bar N^*) \ge
\alpha$ iff $\beta < \alpha \Rightarrow \rk(\bar N,\bar N^*) \ge
\beta$
\sn
\item[$(c)$]  \underline{$\alpha = \beta +1$}:  $\rk(\bar N,\bar N^*) \ge
\alpha$ \Iff \, for every $\gamma < \lambda$ there is $\bar N^+$ such
that
\sn
\item[${{}}$]  \hskip20pt $\bullet \quad \bar N \triangleleft \bar N^+
  \in \mathbf N'$
\sn
\item[${{}}$]  \hskip20pt $\bullet \quad \rk(\bar N^+,\bar N^*) \ge \beta$
\sn
\item[${{}}$]  \hskip20pt $\bullet \quad \ell g(\bar N^+) = \ell
  g(\bar N) +1$
\sn
\item[${{}}$]  \hskip20pt $\bullet \quad$ if 
$n = \ell g(\bar N)$ then $\gamma < \min(N^+_n \backslash N_{n-1})$.
\end{enumerate}
\end{enumerate}
\mn
Consider the statement
\mn
\begin{enumerate}
\item[$\boxtimes$]  for some $\bar N^* \in \mathbf N,\rk(\langle
  N^*_0\rangle,\bar N^*) = \infty$.
\end{enumerate}
\mn
Why enough?  Reflect.

\noindent
Why true?  First
\mn
\begin{enumerate}
\item[$\boxplus_1$]   $E_2$ is a club of $\lambda$ where

\begin{equation*}
\begin{array}{clcr}
E_2 = \{\delta \in E_1: &\text{ if } \bar N^* \in
\mathbf N,\sup(\cup\{N^*_n:n < \omega\} < \delta,\bar N \in \mathbf N', \\
 &\sup(\cup\{\bar N_\ell:\ell < \ell g(\bar N)\} < \delta \text{ and } 
0 \le \rk(\bar N,\bar N^*) < \infty, \text{ \then \, there is no} \\
  &\bar N' \text{ such that } \bar N \triangleleft \bar N' \in \mathbf
N',\rk(\bar N',\bar N^*) = \rk(\bar N,\bar N^*) \text{ and } \ell
g(\bar N') = \ell g(\bar N) +1 \\
  &\text{ such that letting } n = \ell g(\bar N) \text{ we have } 
\min(N'_n \backslash N_{n-1}) \ge \delta\}
\end{array}
\end{equation*}
\end{enumerate}
\mn
[Why?  Reflect.]

Now choose
\mn
\begin{enumerate}
\item[$\boxplus_2$]  there is an increasing sequence $\langle \delta_n:n <
  \omega\rangle$ of members of $E_2$ with limit $\delta \in E_2$ (in
  fact can do this uniformly; e.g. let $\delta_n$ be the $n$-th member
  of $E_2$).
\end{enumerate}
\mn
Lastly, choose $\langle u_{n,\ell}:n < \omega\rangle$ by induction on
$n$ such that
\mn
\begin{enumerate}
\item[$(a)$]  $u_{n,\ell} \in \mathbf U_* \cap [\delta_n]^{< \kappa}$
\sn
\item[$(b)$]  $u_{n,\ell +1}$ is $u^*_\alpha$ for the minimal $\alpha$
 such that $u^*_\alpha \subseteq \delta_n$ and it includes
 $u^*_{n,\ell +1} \cap \delta_n$ where $u^*_{n,\ell +1}$ is
 the $M$-Skolem hull of the set
\sn
\begin{enumerate}
\item[$\bullet$]  $(\cup\{u_{m,k} \cup \{\delta_m\}:
m < \omega,k < \ell\} \cup \{\alpha:\alpha \le 
\sup(u_{n,\ell} \cap \kappa)\}$,
\end{enumerate}
\mn
(the Skolem function are just ``the first example"; note that the
$\sup(u_{n,\ell} \cap \kappa)$ may be zero).
\end{enumerate}
\mn
Let $u_n = \cup\{u_{n,\ell}:\ell < \omega\},N^*_n = M \rest u_n$.  Now
we are done by $(*)_0(a)$ so $\boxplus$ is indeed true and said above
is enough.
\end{PROOF}

\begin{conclusion}
\label{g22}
Assume $\lambda = \mu^+$ is a true successor and $\mu =
\mu^{\aleph_0}$.  \Then \, there is an $\aleph_1$-free Abelian group of
cardinality $\lambda$ such that $\Hom(G,\bbZ) = \{0\}$.
\end{conclusion}

\begin{PROOF}{\ref{g22}}
Straightforward by Theorem \ref{g2} as in \cite{Sh:172} or see in \S(3B).
\end{PROOF}

\begin{theorem}
\label{g23}
1) We can strengthen the conclusion of \ref{g2} by replacing (B)(c) to
\mn
\begin{enumerate}
\item[$(B)$]  $(c)^+ \quad$ if $\langle \bar N'_\eta:\eta \in
  {}^{\omega >} \lambda\rangle$ is as in \ref{g17} (B)(a),(c)-(f) for
  $\kappa = \aleph_1$, replacing

\hskip25pt  (B)(b) by ``$\tau(N_\eta) \subseteq 
\cH(\aleph_0),|N'_\eta| \in [\lambda]^{\le \aleph_0}$" \then \, for
  stationarily 

\hskip25pt  many $\delta \in S$ for some $\alpha < \alpha_{\mathbf w}$ 
and $\eta \in {}^\omega \lambda$ we have $\dot\zeta(\alpha) = \delta$ and

\hskip25pt  $\bar N_\alpha = \langle N'_{\eta \rest n}:n < \omega\rangle$.
\end{enumerate}
\mn
2) In \ref{g2}, if $\kappa < \lambda$ as in \ref{g17} and we can
replace $(\bar N,\bar \gamma)$ by $\langle N_\eta:\eta \in 
{}^{\omega >}\kappa\rangle$.
\end{theorem}

\begin{discussion}
\label{g20}
There is a recent BB helpful in constructing
$\aleph_n$-free abelian groups, (usually is the product of $n$ BB's); in
\cite{Sh:883} it is proved to exist, and using it construct
$\aleph_n$-free Abelian group $G$ such that $\Hom(G,\bbZ) = 0$.  This
is continued, G\"obel-Shelah \cite{Sh:920}, G\"obel-Shelah-Str\"ungman
\cite{Sh:981} use it to deal with modules and in
G\"obel-Herden-Shelah \cite{Sh:970} use it to construct 
$\aleph_n$-free Abelian group with endomorphism ring isomorphic to a
given suitable ring.
Lately is \cite{Sh:1028}.

We try to generalize a version of it but note that we cannot use BB for
$\lambda_{n+1}$ with $\|N_\eta\| = \lambda_n$ as in the ZFC-proof.
But instead we can use \ref{g23}!  See \S(3B) below and maybe more in
\cite{Sh:F1303}.
\end{discussion}
\bigskip

\subsection {Black Boxes with No Choice} \label{3B} \
\bigskip

\begin{context}
\label{k10}
We assume ZF only (for this sub-section).  

Here we try to deal with ZF-proofs.
\end{context}

\noindent
We now define a black box, BB suitable without choice (even weak ones).
\begin{definition}
\label{k11}
1) For a natural number $\mathbf k$ we 
say $\mathbf x$ is a $\mathbf k$-g.c.p. (general 
combinatorial parameter) \when \,
$\mathbf x$ consists of (so $Y = Y_{\mathbf x}$, etc.):
\mn
\begin{enumerate}
\item[$(a)$]  the set $Y$ and the sets
$X_m$ for $m < \mathbf k$ are pairwise disjoint
\sn
\item[$(b)$]  $\Lambda \subseteq \{\bar\eta:\bar\eta = \langle
\eta_m:m < \mathbf k \rangle$ and $\eta_m \in {}^\omega(X_m)$ for $m <
\mathbf k\}$
\sn
\item[$(c)$]  $|Y| \le |X_0|$ and moreover
\sn
\item[$(c)^+$]   $f_0:Y \rightarrow X_0$ is one to one
\sn
\item[$(d)$]  if $m \in (0,\mathbf k)$ then 
$|X_m| \ge {}^{(X_{<m})}Y$ where $X_{<m} =
\prod\limits_{\ell < m} {}^\omega(X_\ell)$, moreover
\sn
\item[$(d)^+$]   $f_m:\{t:t$ a function from $X_{<m}$ to $Y\} 
\rightarrow X_m$ is one to one.
\end{enumerate}
\mn
1A) We say a $\mathbf k$-g.c.p. $\mathbf x$ is standard when $f_{\mathbf x,m}$ is the
identity for every $m < \mathbf k$ and we fix $y_* \in Y$.

\noindent
2)  For $\mathbf x$ a $\mathbf k$-g.c.p. (as above) we say $\mathbf w$ is an $\mathbf
x$-BB, i.e. an $\mathbf x$-black box \when \, $\mathbf w$ 
consists of ($\mathbf x = \mathbf x_{\mathbf w}$ and):
\mn
\begin{enumerate}
\item[$(a)$]   $\Lambda = \Lambda_{\mathbf w} \subseteq \Lambda_{\mathbf x}$; (if
$\Lambda = \Lambda_{\mathbf x}$ we may omit it)
\sn
\item[$(b)$]  $(\alpha) \quad h:\Lambda \rightarrow {}^{(\mathbf k+1) 
\times \omega}Y$, so we write $h(\bar\eta) = \langle
h_{m,n}(\bar \eta):m \le \mathbf k,n < \omega \rangle$ so 

\hskip25pt  $h_{m,n}$ is a function from $\Lambda$ into $Y$
\sn
\item[${{}}$]  $(\beta) \quad$ for every $g:\Omega \rightarrow Y$, see below 
for some $\bar\eta \in \Lambda$ we have 

\hskip25pt $(\forall m < \mathbf k)
(\forall n)(h_{m,n}(\bar\eta) = g(\bar\eta \upharpoonleft (m,n))$
\sn
\item[$(c)$]  \underline{notation}:
\sn
\item[${{}}$]  $(\alpha) \quad$ if $\bar\eta \in \Lambda_{\mathbf x}$ then
$\bar\nu = \bar\eta \upharpoonleft (m,n)$ 
when $\bar\nu = \langle \nu_\ell:\ell < \mathbf k\rangle$ and
  $\nu_\ell$ is $\eta_\ell$ 

\hskip25pt  if $\ell < \mathbf k \wedge \ell \ne m$
and is $\nu_\ell = \eta_\ell \rest n$ if $\ell=m$
\sn
\item[${{}}$]  $(\beta) \quad \Omega_m = \{\bar\eta \upharpoonleft
(m,n):n < \omega$ and $\eta \in \Lambda_{\mathbf w}\}$ so $\Omega_m \subseteq
\{\bar\eta:\bar\eta = \langle \eta_\ell:\ell < \mathbf k\rangle$

\hskip25pt   and for $\ell < \mathbf k,[\ell \ne m \Rightarrow \eta_\ell 
\in {}^\omega(X_\ell)]$ and 
$[\ell = m \Rightarrow \eta_\ell \in {}^{\omega >}X_\ell]\}$
\sn
\item[${{}}$]  $(\gamma) \quad \Omega = \bigcup\limits_{m < \mathbf k} \Omega_m$.
\end{enumerate}
\mn
3) Above $\mathbf k_{\mathbf x} = \mathbf k(\mathbf x) = \mathbf k,\Omega_{\mathbf w} =
\Omega,\Omega_{\mathbf w,m} = \Omega_m$, etc.

\noindent
4) In Claim \ref{k12} below we call 
$\bar z$ \underline{simple} when it has the form $\langle
a_{\bar\eta,n} z:\bar\eta \in \Lambda_{\mathbf x},n < \omega\rangle$
where $a_{\bar\eta,n} \in \bbZ$. 
\end{definition}

\begin{claim}
\label{k11f}
1) For every $Y,y_* \in Y$ and $\mathbf k$ there is, moreover we can define a
 standard $\mathbf k$-g.c.p. $\mathbf x_{\mathbf k}$ (with witnesses 
$f_{\mathbf x,m} =$ identity).

\noindent
2) For every such $\mathbf x_{\mathbf k}$ we can define an 
$\mathbf x$-BB $\mathbf w = \mathbf w_{\mathbf x_{\mathbf k}}$. 
\end{claim}

\begin{remark}
Why we do not choose $\Lambda_{\mathbf w} = \Lambda_{\mathbf x}$?  We can have
$\Lambda_{\mathbf w} = \Lambda_{\mathbf x}$ using a constant value $\in Y$
for the additional cases, so for definability choose a fixed $y_* \in Y$
in \ref{k11}(1), see \ref{k11}(2).
\end{remark}

\begin{PROOF}{\ref{k11f}}
1) By induction $m < \mathbf k$ we define $(X_m,f_m)$ by:
\mn
\begin{enumerate}
\item[$\bullet$]  $X_m = Y$ if $m=0$
\sn
\item[$\bullet$]  $X_m = \{t:t$ is a function from $X_{< m} =
  \prod\limits_{\ell < m} (X_\ell)$ to $Y\}$ if $m > 0$
\sn
\item[$\bullet$]  $f_m = \id_{X_m}$ (so is one to one onto).
\end{enumerate}
\mn
Now check.

\noindent
2) \underline{Case 1}: $\mathbf k = 1$

Let $\Lambda_{\mathbf w} = {}^\omega(\Rang(f_0))$, so $\Omega_{\mathbf w}
= \Omega_{\mathbf w,0} = {}^{\omega >}(\Rang(f_0)),h_{\mathbf w,n}$ or
pedantically $h_{\mathbf w,0,n}$ is a function from $\Omega_{\mathbf w,0}
= \Lambda_{\mathbf w} = \{\langle\eta\rangle:\eta \in
{}^\omega(\Rang(f_0))\}$ and $\Omega_{\mathbf w} = \{\langle \eta
\rangle:\eta \in {}^{\omega >}(\Rang(f_0))\}$ and $\langle \eta\rangle
\in \Lambda_{\mathbf w} \Rightarrow \langle \eta \rangle \upharpoonleft
(0,n) = \langle \eta \rest n \rangle$.

Now for $n < \omega$ we let $h_{\mathbf w,0,n}:\Lambda_{\mathbf w} \rightarrow
Y$ be defined by
\mn
\begin{enumerate}
\item[$\bullet$]  $h_{\mathbf w,0,n}(\langle \eta \rangle) 
= \eta(n) \in Y$ for $\eta \in \Lambda_{\mathbf w}$.
\end{enumerate}
\mn
Obviously clauses (a),(b)$(\alpha)$ from \ref{k11}(2) holds but what
about clause (b)$(\beta)$ of \ref{k11}?

Now for any $g:\Omega_{\mathbf w} \rightarrow Y$ we choose $y_n \in Y$ by
induction on $n$ as follows: $y_n = g(\langle f_0(y_\ell):
\ell < n\rangle) = g(\langle y_\ell:\ell < n\rangle)$.  
So $\eta := \langle y_\ell:\ell < \omega\rangle \in
{}^\omega(\Rang(f_0))$ is as required. 
\medskip

\noindent
\underline{Case 2}:  $\mathbf k > 1$

Let $\Lambda_{\mathbf w} = \{\bar\eta:\bar\eta = \langle \eta_m:m <
\mathbf k\rangle$ and $\eta_m \in {}^\omega(\Rang(f_m))$ for $m < \mathbf
k\}$ hence $\Omega_m = \Omega_{\mathbf w,m}$ and $\Omega_* = \Omega_{\mathbf w}$
are well defined.

We now define $h_{m,n} = h_{\mathbf w,m,n}$ for $m < \mathbf k,n < \omega$
\mn
\begin{enumerate}
\item[$(*)_1$]  for $\bar\eta \in \Lambda_m = \{\bar\eta
  \upharpoonleft (m,n):\bar\eta \in \Lambda_{\mathbf w}$ and $n <
  \omega\}$ we let $h_{m,n}(\bar\eta) =
\big(f^{-1}_m(\eta_m(n))\big)(\bar\eta \rest m)$ if $m>0$ and
  $h_{m,n}(\bar\eta) = f^{-1}_m(\eta_m(n))$ if $m=0$.
\end{enumerate}
\mn
Why well defined and $\in Y$?   Clearly if $m=0$ then
$h_{m,n}(\bar\eta) = f^{-1}_m(\eta_m(n)) \in Y$ as $\eta_m \in
{}^\omega(X_0) = {}^\omega Y$ and if $m>0$ then
$\eta_m(n) \in X_m$ hence $f^{-1}_m(\eta_m(n)) \in {}^{(X_{<m})}Y$ so
is a function from $X_{<m} := \prod\limits_{\ell < m} {}^\omega
X_\ell$ into $Y$
so $\bar\eta \rest m \in X_{<m}$ hence
$\big(f^{-1}_m(\eta_m(n))\big)(\bar\eta \rest m) \in Y$.  So
clause $(b)(\alpha)$ of Definition \ref{k11} is satisfied.  What about
clause $(b)(\beta)$ of Definition \ref{k11}(2)? so let a function
$g:\Omega \rightarrow Y$ be given and we shall prove that there is $\bar\eta
\in \Lambda_{\mathbf w}$ as required, in fact define it.  Toward this we
choose $\eta_m \in {}^\omega(\Rang(f_m)) \subseteq {}^\omega(X_m)$ by
downward induction on $m$, and for each $n$, we shall let $\eta_m =
\langle f_m(t_{m,n}):n < \omega\rangle$, where we choose $t_{m,n} \in
\Dom(f_m)$ by induction on $n < \omega$ as follows:
\mn
\begin{enumerate}
\item[$(*)_2$]  if $m>0$ then $t_{m,n}$ is the 
following function from $\{\bar\eta
\rest m:\bar\eta \in \Lambda_{\mathbf w}\} = X_{<m} = \prod\limits_{\ell
  < \omega} {}^\omega(X_\ell)$ to $Y$: if $\bar\nu =
\langle \nu_\ell:\ell < m\rangle \in \Dom(t_{m,n})$ then $t_{m,n}(\bar
\nu)$ is $g(\bar\rho) \in Y$ where $\bar\rho = \langle \rho_\ell:\ell
< \mathbf k\rangle$ is defined by:
\sn
\begin{enumerate}
\item[$\bullet$]  if $\ell > m$ then $\rho_\ell = \eta_\ell$, is well
defined by the induction hypothesis on $m$
\sn
\item[$\bullet$]  if $\ell=m$ then $\rho_\ell = \langle
f_m(t_{m,0}),\dotsc,f_m(t_{m,n-1})\rangle$, well defined by the
induction hypothesis on $n$
\sn
\item[$\bullet$]  if $\ell < m$ then $\rho_\ell = \nu_\ell$, given
\end{enumerate}
\sn
\item[$(*)_3$]  if $m = 0$ then $t_{m,n} = g(\bar\rho)$ where
  $\bar\rho$ is chosen as above except that there is no $\bar\nu$.
\end{enumerate}
\mn
Now check.
\end{PROOF}

\begin{claim}
\label{k12}
Let $\mathbf x$ be a $\mathbf k$-g.c.p. see \ref{k11}(1) and $\mathbf w$ an
$\mathbf x$-BB, see \ref{k11}(2) and $\Lambda = \Lambda_{\mathbf w},\Omega
= \Omega_{\mathbf w}$, etc.  Then $G \in \cG_{\mathbf x}
\Rightarrow G_{\mathbf x,0} \subseteq G \subseteq_{\purely} G_{\mathbf
x,1}$ where $\subseteq_{\purely}$ is from \ref{k15}(0) and
$G \in \cG_{\mathbf x}$ \Iff \, some $\bar z,G = G_{\mathbf x,\bar z}$,
which means:
\mn
\begin{enumerate}
\item[$(a)$]  $G_0 = G_{\mathbf x,0} =
\oplus\{\bbZ x_\rho:\rho \in \Omega\} \oplus \bbZ z$
\sn
\item[$(b)$] $G_1 = G_{\mathbf x,1} = \oplus\{\bbQ x_\rho:\rho \in \Omega\}
 \oplus \bbQ z \oplus \{\bbQ y_{\bar\eta}:\bar\eta \in \Lambda_{\mathbf x}\}$
\sn
\item[$(c)$]  $\bar z = \langle z_{\bar\eta,n}:\eta \in \Omega_{\mathbf
w}\rangle$ is a sequence of members of $G_{\mathbf x,1}$
\sn
\item[$(d)$]  for $\bar\eta \in \Lambda$ we define $y_{\bar\eta,n} =
y_{\bar z,\bar\eta,n}$ by induction on $n$:
\sn
\begin{enumerate}
\item[$\bullet$]  $y_{\bar\eta,0} = y_{\bar\eta}$,
\sn
\item[$\bullet$]  $n!y_{\bar\eta,n+1} =
y_{\bar\eta,n} - \sum\limits_{m < \mathbf k} x_{\bar\eta \upharpoonleft
 (m,n+1)} - z_{\bar\eta,n}$
\end{enumerate}
\sn
\item[$(e)$]  $G$ is the (Abelian) subgroup of $G_1$
 generated by $\{x_{\bar\eta}:\bar\eta \in \Omega\} \cup
 \{y_{\bar\eta,m}:\bar\eta \in \Lambda,n < \omega\} \cup \{z\}$.
\end{enumerate}
\end{claim}

\begin{PROOF}{\ref{k12}}
Straightforward.
\end{PROOF}

\begin{claim}
\label{k13}
Let $\mathbf k,\mathbf x,\mathbf w,\bar z$ be as in \ref{k11}, \ref{k11}(2),
\ref{k12}.

\noindent
1) $G_{\mathbf x,\bar z}$ is almost $\aleph_{\mathbf k(\mathbf x)}$-free (see below
   Definition \ref{k15} and \ref{k14}) provided that $\bar z$ has the
   form $\langle a_{\bar\eta,n} z:\bar\eta \in \Lambda_{\mathbf x},n <
   \omega\rangle$ where $a_{\bar\eta,n} \in \bbZ$ (or less as in
   \cite{Sh:883}).

\noindent
1A) If $\mathbf k \ge 2$ then $G_{\mathbf x,\bar z}$ is strongly
$\aleph_1$-free.

\noindent
2) In Claim \ref{k12} above,
$G_{\mathbf x,\bar z}$ is definable (in ZF!) from $(\mathbf x,\bar z)$.

\noindent
3) For $\mathbf x$ a $\mathbf k$-g.c.p. and $\mathbf w$ an $\mathbf x$-BB such
that $Z \subseteq Y_{\mathbf x}$ we can define $\bar z = 
\bar z_{\mathbf w}$ such that $G_{\mathbf x,\bar z}$ (is
   well defined and) satisfies 
$h \in \Hom(G_{\mathbf x,\bar z},\bbZ) \Rightarrow h(z)=0$.

\noindent
4) For $\mathbf x$ a $\mathbf k$-g.c.p. and $\mathbf w$ an $\mathbf x$-BB we can
   define an $\aleph_{\mathbf k(\mathbf x)}$-free Abelian group $G$ such
   that $\Hom(G,\bbZ) = \{0\}$. 
\end{claim}

\begin{discussion}
\label{k14}
1) Assume $H \subseteq G = G_{\mathbf x,\bar z}$ is a subgroup of
cardinality $< \aleph_{\mathbf k(\mathbf x)}$.  For each $t \in G$ let
$Y_t$ be the minimal $Y \subseteq Y_{\mathbf x} = \{x_\rho:\rho \in
\Omega_{\mathbf x}\} \cup \{z\} \cup \{y_{\bar\eta}:\eta \in
\Lambda_{\mathbf x}\}$ such that $t \in \oplus\{\bbQ x:x \in Y\}$.  If
$\Omega_{\mathbf x} \cup \Lambda_{\mathbf x}$ is linearly ordered then
$\cup\{Y_t:t \in H\}$ has cardinality $< \aleph_{\mathbf k(\mathbf x)}$
but in general this explains the ``weakly" or ``almost" in \ref{k13}.
However, it
may occur that this holds for the ``wrong" reason say $\aleph_0
\nleq |A|$ in Definition \ref{k15}(2). But the proof of \ref{k11f},
\ref{k11} gives ``many" such subsets of the set $A$.

\noindent
2) For proving \ref{k13}(1) note that in the definition of $\cG_{\mathbf x}$ in
\cite{Sh:883} there is a use of choice: dividing the stationary set
$S_m \subseteq \lambda_m$ to $\lambda_m$ pairwise disjoint sets or just
the choice of $\bar z = \langle z_{\bar\eta}:\bar\eta \in
\Lambda_{\mathbf w}\rangle$.  However, we can just ``glue together"
copies of the $G$ constructed above; i.e. start with $G$ and for every
non-zero pure $z \in G$, add $G_z$ of $h_z:G \rightarrow G_z$ identify
$x_{<>}$ with $z$, etc.
\end{discussion}

\begin{definition}
\label{k15}
Let $G$ be a torsion free Abelian group (the torsion free means $G
\models ``nx=0",n \in \bbZ,x \in G$ implies $n=0 \vee x = 0_G$).

\noindent
0) Recall $H \subseteq G$ means $H$ is a subgroup.  Let 
$H \subseteq_{\purely} G$ mean $H$ is a
pure subgroup of $G$, which means $H \subseteq G$ and $n \in \bbZ \backslash
\{0\},nx \in G,nx \in H \Rightarrow x \in H$.

\noindent
1) We say $G$ is a weakly $\kappa$-free \when \,: there is a set $A$
   such that the pair $(G,A)$ is $\kappa$-free, see part (2).

\noindent
2) We say $(G,A)$ is $\kappa$-free when: $A \subseteq G$ and
$\PC_G(A)=G$ and if $B \subseteq A$ has cardinality $< \kappa$ 
then $\PC_G(B) \subseteq G$ is a free Abelian group 
recalling $\PC_G(A) =$ the minimal pure
subgroup of $G$ which includes $A$.

\noindent
3) We say $G$ is almost $\kappa$-free \when \, there is a set $A$ such
   that the pair $(G,A)$ is almost $\kappa$-free, see part (4).

\noindent
4) The pair $(G,A)$ is almost $\kappa$-free when: $(G,A)$ is 
 $\kappa$-free and $A$ is independent in $G$
(i.e. $\sum\limits_{\ell < n} a_\ell x_\ell = 0 \Rightarrow
   \bigwedge\limits_{\ell < n} a_\ell = 0$ when $x_0,\dotsc,x_n \in A$
  without repetitions).
\end{definition}

\begin{PROOF}{\ref{k13}}

\noindent
\underline{Proof of \ref{k13}}:

\noindent
1) Let $A = \{x_\rho:\rho \in \Omega_{\mathbf x}\} \cup \{z\} \cup
   \{y_{\bar\eta}:\bar\eta \in \Lambda_{\mathbf x}\}$.  It is easy to
check that $A$ is independent in $G$ (see \ref{k15}(4)) and $\PC_G(A)
   = G$ so for any $t \in G$ there is a unique finite $Y_t
   \subseteq A$ such that $t \in \PC_G(Y_t),Y_t$ of minimal
   cardinality.

Now if $B \subseteq A$ has cardinality $< \aleph_{\mathbf k(\mathbf x)}$,
then also $Y_B := \{\rho:x_\rho \in B\} \cup \{\bar\eta
\upharpoonleft (m,n):y_{\bar\eta} \in B,m < \mathbf k(\mathbf x)$ and $n <
\omega\}$ has cardinality $< \aleph_{\mathbf k(\mathbf x)}$.

For some $Y \subseteq \Ord$ in $\mathbf L[Y]$ there is a $\mathbf
k$-c.p. $\mathbf x'_1$ and $\bar z_1$ such that $G_{\mathbf x_1,\bar z_1}
\in \mathbf L[Y]$ is isomorphic (in $\mathbf V$) to $\PC_G(B)$.  So by
\cite{Sh:883} we are done.

\noindent
2) Should be clear.

\noindent
3) We shall define uniformly (in ZF) from $\mathbf k$-g.c.p. 
$\mathbf x$ and $\mathbf w$ an $\mathbf x$-BB a sequence $\bar z$ such that the 
Abelian group $G = G_{\mathbf x,z_{\mathbf w}}$ satisfies $h \in \Hom(G,\bbZ)
\Rightarrow h(z) = 0$.

For each $\bar\eta \in \Lambda$ let $\bar a = \langle a_{\mathbf
w,\bar\eta,n}:n < \omega\rangle \in {}^\omega \bbZ$ be defined by:
\mn
\begin{enumerate}
\item[$(*)$]  $a_{\mathbf w,\bar\eta,n}$ is
\sn
\item[${{}}$]  $\bullet \quad \sum\limits_{m < \mathbf k}
h_{m,n+1}(\bar\eta)$ when $\{h_{m,n}(\bar\eta):m < \mathbf k\} \subseteq
\bbZ$
\sn
\item[${{}}$]  $\bullet \quad 0 \qquad$ when otherwise.
\end{enumerate}
\mn
We shall choose $b_{\mathbf w,\bar\eta,n} \in \bbZ$ for $n < \omega$ such that
\mn
\begin{enumerate}
\item[$(*)$]  if $a_{\mathbf w,\bar\eta,0} \ne 0$ \then \, there are no
$t_n \in \bbZ$ for $n < \omega$ such that for every 
$n$ we have
\newline
$(\eq_n) \, n!t_{n+1} = t_n - a_{\mathbf w,\bar\eta,n+1} = b_{\mathbf
w,\bar\eta,n} \cdot a_{\mathbf w_{\eta,0}}$.
\end{enumerate}
\mn
Why then can we choose?  We choose $b_{\mathbf w,\bar\eta,n} \in \bbN
\subseteq \bbZ$ as minimal such that we cannot find $t_0,\dotsc,t_n
\in \bbZ$ such that $t_0 = \{-n,-n-1,\dotsc,-1,0,1,m\dotsc,n\}$ and
for every $m < n+$ we have $\bbZ \models ``n!t_{m+1} = t_n - a_{\mathbf
  w,\bar\eta,m+1} - b_{\mathbf w,\bar\eta,m} - a_{\mathbf w_{\eta,0}}$.

Now we define
\mn
\begin{enumerate}
\item[$(*)$]  $\bar z = \bar z_{\mathbf w} = \langle b_{\mathbf
  w,\bar\eta,n} \cdot z:\bar\eta \in \Lambda_{\mathbf x},n <
  \omega\rangle$.
\end{enumerate}
\mn
So
\mn
\begin{enumerate}
\item[$(*)$]  $(a) \quad G_{\mathbf x,\bar z}$ is well defined
\sn
\item[${{}}$]  $(b) \quad$ if $g \in H(G_{\mathbf x,\bar z},\bbZ)$ then
  $h(z) = 0_{\bbZ}$.
\end{enumerate}
\mn
[Why?  Clause (a) is obvious.  For clause (b) if $g$ is a
  counterexample by the choice of $\mathbf w$ there is $\bar\eta \in
  \Lambda_{\mathbf w}$ such that $m < \mathbf k \wedge n < \omega
  \Rightarrow g(x_{\bar\eta \upharpoonleft (m,n)}) =
  h_{m,n}(\bar\eta)$ that is $n < \omega \Rightarrow \sum\limits_{m
    < \mathbf k} g(x_{\bar\eta \upharpoonleft (m,n+1)} = a_{\mathbf
    w,\bar\eta,n}$.  Now use the choice of $\langle b_{\mathbf
    w,\bar\eta,n}:n < \omega\rangle$ to get a contradiction.]

\noindent
4) We derive an example from $G_{\mathbf w}$ from part (3).

Let $\Omega' = \Omega'_{\mathbf x} = \{\rho:\rho$ a finite sequence
   of members of $\Omega\}$ and for $\rho \in \Omega'$ let
\mn
\begin{enumerate}
\item[$(*)$]  $(a) \quad X_\rho = X_{\mathbf x,\rho} =
  \{x_{\rho,\bar\eta}:\bar\eta \in \Omega_{\mathbf w}\}$
\sn
\item[${{}}$]  $(b) \quad Y_\rho = Y_{\mathbf x,\rho} =
  \{y_{\rho,\bar\eta}:\bar\eta \in \Lambda_{\mathbf w}\}$
\sn
\item[$(*)$]  $(a) \quad G'_0 = G'_{\mathbf x,0} = G'_{0,0} \otimes
  G'_{0,1}$ where
\sn
\item[${{}}$]  $(b) \quad G'_{0,0} = G'_{\mathbf x,0,0} = 
\oplus\{\bbZ_{\rho,\bar\eta}:\rho \in \Omega'_{\mathbf w},\bar\eta \in
  \Omega_{\mathbf w}\}$
\sn
\item[${{}}$]  $(c) \quad G'_{0,1} = G'_{\mathbf x,0,1} = \bbZ z$
\sn
\item[$(*)$]  $(a) \quad G'_1 = G'_{\mathbf x,1} \oplus G_{\mathbf w,2,1}
  \oplus G_{\mathbf w,2,1} \oplus G_{\mathbf w,1,2}$ where
\sn
\item[${{}}$]  $(b) \quad G'_{1,0} = G'_{\mathbf w,1,0} = \oplus\{\bbQ
  x_{\rho,\bar\eta}:\rho \in \Omega'_{\mathbf x}$ and $\bar\eta \in
  \Omega_{\mathbf x}\} \supseteq G'_{0,1}$
\sn
\item[${{}}$]  $(c) \quad G'_{1,1} = G'_{\mathbf w,1,1} = \bbQ z
  \supseteq G'_{0,1}$
\sn
\item[${{}}$]  $(d) \quad G'_{1,2} = G'_{\mathbf w,1,2} = \oplus\{\bbQ
  y_{\rho,\bar\eta}:\rho \in \Omega'_{\mathbf x}$ and $\bar\eta \in
  \Lambda_{\mathbf x}\}$.
\end{enumerate}
\mn
Let
\mn
\begin{enumerate}
\item[$(*)$]  $(a) \quad z_\rho$ be $z$ if $\rho =\langle \rangle$ and
  $x_{\rho \rest \ell,\rho(\ell)}$ if $\beta \in \Omega'_{\mathbf x}
  \backslash \{<>\}$
\sn
\item[${{}}$]  $(b) \quad$ let $y_{\rho,\bar\eta,0} =
  y_{\rho,\bar\eta}$
\sn
\item[${{}}$]  $(c) \quad$ for $\rho \in \Omega'_{\mathbf w}$ and
  $\bar\eta \in \Lambda_{\mathbf x}$ we define $y_{\rho,\bar\eta,n}$ by
  induction on $n>0$
\sn
\begin{enumerate}
\item[${{}}$]  $\bullet \quad y_{\rho,\bar\eta,n+1} =
  (y_{\rho,\bar\eta,n} + \sum\limits_{m \le k} x_{\rho\bar\eta
  \upharpoonleft (m,n)} + \bar a_{\bar\eta,n} z_{\bar\eta}$ where

\hskip30pt $\langle a_{\bar\eta,n}:n < \omega\rangle \in {}^\omega \bbZ$ was
  defined above using $h(\bar\eta)$
\end{enumerate}
\sn
\item[$(*)$]  $(a) \quad$ for every $t \in G'_1$ let $\supp(x)$ be the
  minimal subset $X_t$ of $X_{\mathbf s} = \{x_{\rho,\bar\eta}$:

\hskip25pt $\rho \in \Omega'_{\mathbf x},\bar\eta \in \Omega_{\mathbf x}\} \cup
  \{y_{\rho,\bar\eta}:\rho \in \Omega'_{\mathbf x}$ and $\bar\eta \in
  \Lambda_{\mathbf w}\}$ such that: 

\hskip25pt  $t \in \Sigma\{\bbQ x:x \in X_*\}$; used in part (2)
\sn
\item[$(*)$]  for $\rho \in \Omega'$ we define an embedding $h_\rho$
  from $G_{\mathbf w}$ into $G'_1$ by (see $\boxplus_4$ below):
\sn
\item[${{}}$]  $(a) \quad h_\rho(z) = z_\rho$
\sn
\item[${{}}$]  $(b) \quad h_\rho(x_{\bar\eta}) = x_{\rho,\bar\eta}$
  for $\bar\eta \in \Omega_{\mathbf w}$
\sn
\item[${{}}$]  $(c) \quad h_\rho(y_{\bar\eta,n}) =
  y_{\rho,\bar\eta,n}$.
\end{enumerate}
\mn
Now
\mn
\begin{enumerate}
\item[$\boxplus_1$]   let $G'_{\mathbf w}$ be the subgroup of
  $G'_{\mathbf w,1}$ generated by $\{X_{\rho,\bar\eta}:\rho \in
  \Omega'_{\mathbf w}$ and $\bar\eta \in \Omega_{\mathbf w}\} \cup \{z\}
\cup \{y_{\rho,\bar\eta,n}:\rho \in \Omega'_{\mathbf w},\bar\eta \in
  \Lambda_{\mathbf w}$ and $n < \omega\}$
\sn
\item[$\boxplus_2$]    $G'_{\mathbf w,0} \subseteq G'_{\mathbf x}$ is
  dense in the $\bbZ$-adic topology.
\end{enumerate}
\mn
[Why?  Just look at each $y_{\rho,\bar\eta,n}$.]
\mn
\begin{enumerate}
\item[$\boxplus_3$]   for $\rho \in \Omega'_{\mathbf x}$
\sn
\item[${{}}$]  $(a) \quad h_\rho$ is a well defined homomorphism
\sn
\item[${{}}$]  $(b) \quad h_\rho$ is indeed an embedding
\sn
\item[${{}}$]  $(c) \quad \Rang(h_\rho) \subseteq G'_{\mathbf x}$
\sn
\item[${{}}$]  $(d) \quad \Rang(h_\rho)$ is a pure subgroup of
  $G'_{\mathbf x}$
\sn
\item[${{}}$]  $(e) \quad h_{<>}$ is ?
\end{enumerate}
\mn
[Why?  For clause (a) note the definition of $y_{\rho,\bar\eta,n}$,
  also the other clauses are obvious.]
\mn
\begin{enumerate}
\item[$\boxplus_4$]   $\Hom(G'_{\mathbf w},\bbZ) = 0$.
\end{enumerate}
\mn
[Why?  Let $g \in \Hom(G'_{\mathbf w},\bbZ)$.  For each $\rho \in
  \Omega'_{\mathbf x}$, the function $g \circ h_\rho$ is a homomorphism
from $G_{\mathbf x}$ into $\bbZ$ hence by the previous claim \ref{k13},
$(G \circ h_\rho)(z)=0$.  This means that $0 = (g \circ h_\rho)(z) =
g(h_\rho(z)) = g(z_\rho)$ hence $g(z)=0$, using $\rho = \langle
\rangle$ and $g(x_{\rho,\bar\eta})=0$ for $\rho \in \Omega'_{\mathbf
  x},\bar\eta \in \Omega_{\mathbf x}$ using $z_{\rho \char 94 \langle
  \bar z \rangle} = X_{\rho,\bar\eta}$.  By the choice of $G'_{\mathbf
  w,0}$ this implies $g \rest G'_{\mathbf x,0}$ is zero and by
$\boxplus_3$ this implies $g \rest G'_{\mathbf w}$ is zero, as
promised.]
\end{PROOF}
\newpage


\bibliographystyle{amsalpha}
\bibliography{shlhetal}

\end{document}